\documentclass{amsart}
\usepackage{amsfonts}
\usepackage{amsmath}
\usepackage{amssymb}
\usepackage{latexsym}
\input cyrfont.def
\input cyracc.def
\newtheorem{theorem}{Theorem}[section]
\newtheorem{lemma}[theorem]{Lemma}

\newtheorem{proposition}[theorem]{Proposition}
\theoremstyle{definition}
\newtheorem{definition}[theorem]{Definition}
\theoremstyle{remark}
\newtheorem{remark}[theorem]{Remark}
\numberwithin{equation}{section}
\DeclareMathOperator{\spec}{spec}%
\DeclareMathOperator{\cov}{cov}%
 
\renewcommand{\theenumi}{\textup{\roman{enumi}}}
\newcommand{\defi}{\stackrel{\textup{\tiny def}}{=}}
\DeclareMathOperator{\fug}{\text{\Large $\pi$}}

\begin{document}
\title[Deformations of Fuchsian Systems and the Schlesinger System]%
{Deformations of  Fuchsian Systems\\ of Linear Differential Equations \\
 and the Schlesinger System.}
\author[V.Katsnelson]{Victor Katsnelson}
\address{Department of Mathematics,\newline
\hspace*{\parindent}the Weizmann Institute of Science,\newline%
\hspace*{\parindent}Rehovot\\ 76100 \\
Israel}%
\email{victor.katsnelson@weizmann.ac.il}
\thanks{The research of Victor Katsnelson was supported
by the Minerva Foundation.}
\author[D. Volok]{Dan Volok}
\address{Department of Mathematics,\newline
\hspace*{\parindent}the Weizmann Institute of Science,\newline%
\hspace*{\parindent}Rehovot\\ 76100 \\
Israel}%
\email{dan.volok@weizmann.ac.il} 
\date{December 5, 2004}
\dedicatory{To the  centenary of the Schlesinger system}
\begin{abstract}{\small%
We consider holomorphic deformations of Fuchsian systems
parameterized by the pole loci. It is well known that, in the case
when the residue matrices are non-resonant, such a deformation is
isomonodromic if and only if the residue matrices satisfy  the
Schlesinger system with respect to the parameter. Without the
non-resonance condition this result fails: there exist
non-Schlesinger isomonodromic deformations. In the present article
we introduce the class of the so-called isoprincipal deformations of
Fuchsian systems. Every isoprincipal deformation is also an
isomonodromic one.  In general, the class of the isomonodromic
deformations is much richer than the class of the isoprincipal
deformations, but in the non-resonant case these classes coincide.
We prove that a deformation is isoprincipal if and only if the
residue matrices satisfy the Schlesinger system. This theorem holds
in the general case, without any assumptions on the spectra of the
residue matrices of the deformation.  An explicit example
illustrating isomonodromic deformations, which are neither
isoprincipal nor meromorphic with respect to the parameter,  is also
given.
}%
\end{abstract}
\maketitle

\section*{NOTATION}
\begin{itemize}
\item
 $\mathbb C$  stands for the complex plane.

\item
 $\mathbb{CP}^1$  stands for the extended
complex plane ($=$ the Riemann sphere):
$$\mathbb{CP}^1=\mathbb{C}\cup\infty.$$

\item
 ${\mathbb C}^n$  stands for the
$n$-dimensional complex space.

\item
 In the coordinate notation, a point
$\boldsymbol{t}\in{\mathbb C}^n$ will be written as
$\boldsymbol{t}=(t_1,\ldots,t_n).$

\item
$\mathbb{C}^n_*$ is the set of points $\boldsymbol{t}\in{\mathbb
C}^n,$ whose coordinates $t_1,\dots,t_n$  are pairwise different:
$${\mathbb{C}}^n_{\ast}={\mathbb C}^n\setminus\bigcup_{\substack{1\leq
i,j\leq n\\ i\not=j}} \{\boldsymbol{t}:t_i=t_j\}.$$

\item $\mathfrak{M}_k$ stands for
the set of all $k\times k$ matrices with complex entries.

\item
 $[\cdot,\cdot]$ denotes the commutator: for
 $A,B\in\mathfrak{M}_{k},\ [A,B]=AB-BA$.

\item$I$  stands for the identity
matrix of an appropriate dimension.

\end{itemize}

\setcounter{section}{-1}
\section{Introduction}
The systematic study of linear differential equations in the complex
plane with coefficients dependent on parameters has been started by
Lazarus Fuchs in the late eighties of the nineteenth
century\footnote{L.Fuchs died in 1902.}
\cite{FuL1}, \cite{FuL2}, \cite{FuL3}, \cite{FuL4}. In particular,
L.Fuchs investigated the equations whose  monodromy  does not depend
on such parameters. These investigations were continued in the
beginning of the twentieth century by
L.Schlesinger%
\footnote{A student of L.Fuchs.}, \cite{Sch2}, \cite{Sch3},
\cite{Sch4},   R.Fuchs\footnote{The son of L.Fuchs.}, \cite{FuR1},
\cite{FuR2}, and R. Garnier, \cite{Gar}.

L.Schlesinger's research was closely related to the Hilbert 21st
problem (a.k.a. the Riemann - Hilbert monodromy problem), which
requires to construct a Fuchsian system with prescribed monodromy
(for the explanation of terminology see Section \ref{Fu} of the
present article). In the paper \cite{Sch2}, which appeared  exactly
one hundred years ago -- in 1905, L.Schlesinger proposed the idea
that it would be very fruitful to study the {\em deformations} of
Fuchsian systems
\begin{equation}\label{eqf}
\dfrac{dY}{dx}=\sum_{1\leq j\leq
n}\dfrac{Q_j(\boldsymbol{t})}{x-t_j}\cdot Y,\end{equation} where the
residues $Q_j$ depend holomorphically on  the pole loci
$\boldsymbol{t}=(t_1, \ldots , t_n)$, and investigate the dependence
of the solution $Y$ on  $\boldsymbol{t}$, as well as on $x$.

Emphasizing this idea, L.Schlesinger explained that he was guided by
the analogy with the theory of algebraic functions, where he had
studied algebraic functions as functions of both the "main variable"
and the loci of {\em ramification points} considered as parameters
(see \cite[pp.  287\,-\,288]{Sch1}).

Also in the paper \cite{Sch2}, the system of PDEs
\begin{equation}%
\label{Sch0} \left\{
\begin{aligned}
\dfrac{\partial Q_{j}}{\partial t_{k}}&=
\dfrac{[Q_{j},Q_{k}]}{t_{j}-t_{k}}, &\quad
 1\leq
j,k\leq n,\ k\not=j, \\[1.0ex]
\dfrac{\partial Q_{j}}{\partial t_{j}}&= - \sum_{\substack{1\leq
k\leq n
\\
k\not=j}} \dfrac{\left[Q_{j},Q_{k}\right]} {t_{j}-t_{k}}, & \quad
1\leq j\leq n,
\end{aligned}\right.
\end{equation}
which is now known as {\em the Schlesinger system}, was introduced
and the statement  that \textsl{the holomorphic deformation
\eqref{eqf} is isomonodromic if and only if its coefficients
$Q_j(\boldsymbol{t})$ satisfy the system \eqref{Sch0}} was
formulated. (See page 294 of \cite{Sch2}, four bottom lines of this
page.) This formulation was repeated in the book \cite{Sch3}, pp.
328-329, and later in the paper \cite{Sch4}, p. 106. (References to
the earlier paper \cite{Sch2} are relatively rare. Usually, one
refers  to the more recent paper \cite{Sch4}.)

Over the years the Schlesinger system and the isomonodromic
deformations of Fuchsian systems were extensively studied;  we would
like to mention in particular the papers of T. Miwa \cite{Miwa} and
of B. Malgrange \cite{Ma1}, where it was proved that the Schlesinger
system enjoys the Painlev\'e property (its solutions  are
meromorphic functions in the universal covering space over
$\mathbb{C}^n_*$), and the book by A.R. Its and V.{Yu}.\:Novokshenov
\cite{ItNo}, where  the connections between the isomonodromic
deformations and the transcendents of Painlev\'e were revealed.

However, in the 1990s the famous negative solution to the Riemann -
Hilbert monodromy problem due to A.A. Bolibrukh (see  \cite{Bol1},
\cite{Bol2}) gave a  strong motivation for the {\em revision} of the
classical results, concerning the isomonodromic deformations of the
Fuchsian systems.

For example, it should be noted that in the original works
\cite{Sch2} -- \cite{Sch4} of L. Schlesinger no assumptions
concerning the non-resonance of the matrices $Q_j$ are made. {\em In
such generality the above-cited statement of L. Schlesinger fails.}
If a holomorphic deformation \eqref{eqf} of Fuchsian systems is such
that the residues $Q_j$ satisfy  the Schlesinger system
\eqref{Sch0}, then this deformation is isomonodromic, but the
converse statement is not true, in general.

It was also A.A. Bolibrukh who constructed the first explicit
example of the {\em non-Schlesinger  isomonodromic} deformation. In
this example the monodromy is non-trivial and the residues
$Q_j(\boldsymbol{t})$ are rational functions of $\boldsymbol{t}$
(see \cite{Bol3} and \cite{Bol4} -- in both papers the same example
appears as Example 2; see also\,\footnote{Unfortunately, to our best
knowledge this review has not yet been translated into English.}
Section 3 of the review paper \cite{Bol5}, where this example
appears as Example 3).

 At the same time it was shown independently  in \cite{Kats1}  that {\em almost
every isomonodormic deformation of  Fuchsian systems with generic
rational solutions}\,\footnote{In particular, with trivial
monodromy.} {\em is non-Schlesinger} (for more details see Remark
\ref{ashkelon} in Section \ref{Examp} of this article).

Thus {\em the isomonodromic property of the deformation \eqref{eqf}
implies the Schlesinger system for the residues
$Q_j(\boldsymbol{t})$ under the non-resonance condition, but not in
general.}

Unfortunately,  careless treatment of the non-resonance condition is
very common in the history of problems related to the monodromy of
Fuchsian systems. It can also be found in some works of
 V. Volterra and of G.D. Birkhoff (see
\cite[Chapter XV, \S 9]{Gant} for details). This tradition continues
to certain extent in the above-mentioned paper \cite{Miwa} of T.
Miwa on the Painlev\'e property of isomonodromic deformations: in
this paper the non-resonance condition appears as the equation
(2.22), but is omitted both in the introduction and in the
formulation of the main result.

 Without the assumption of non-resonance  the main result
 of \cite{Miwa} does {\em not} hold: there exist isomonodromic deformations of the form
 \eqref{eqf}, where the residues $Q_j(\boldsymbol{t})$
 are {\em not} meromorphic  in
the universal covering space over $\mathbb{C}^n_*$. (The appropriate
example is presented in Section \ref{Examp} of this article. We note
that this phenomenon does not occur in the above-mentioned example
of the {non-Schlesinger  isomonodromic} deformation due to A.A.
Bolibrukh: in that example the residues $Q_j(\boldsymbol{t})$ are
rational functions of $\boldsymbol{t}$.)

The main goal of the present work is to answer the following
question: {\em how to describe the class of holomorphic deformations
\eqref{eqf} with the property that the residues
$Q_j(\boldsymbol{t})$ satisfy the Schlesinger system \eqref{Sch0},
when one omits the non-resonance assumption?}

The presentation of our results is organized as follows.

In the first section after this introduction we recall the basic
notions concerning the Fuchsian system and introduce a certain
canonical multiplicative decomposition of the fundamental solution
in a neighborhood of its singular point $t_j$. This is the so-called
{\em regular-principal factorization}: the fundamental solution is
represented as the product of a {\em regular factor,} holomorphic
and invertible at the point $t_j,$ and a {\em principal factor,}
holomorphic (multi-valued) and invertible everywhere {\em except} at
$t_j.$ This {\em principal factor} is the multiplicative analogue of
the {\em principal part} in the Laurent decomposition: it contains
the information about the nature of the singularity.

In Section \ref{DefFDE} we  introduce the main notion of the present
article (Definition \ref{DefIsoprDefo}): the so-called {\em
isoprincipal}\,\footnote{Iso- (from
\(\stackrel{.\!.}{\iota}\)\(\sigma o\varsigma\)
 \,-\,equal\,-\, in Old Greek)
 is a combining form.} families of Fuchsian systems. These are the
holomorphic families \eqref{eqf} with the property that all the
principal factors of a suitably normalized fundamental solution
$Y(x,\boldsymbol{t})$ are, in a certain sense, preserved. We show
that every isoprincipal family is also isomonodromic and that the
converse is true, when the non-resonance condition is in force.

In Section  \ref{SchlSyst} we formulate and prove our main result
(Theorem \ref{main}) that {\em the family \eqref{eqf} is
isoprincipal if and only if the residues $Q_j(\boldsymbol{t})$
satisfy  the Schlesinger system \eqref{Sch0}.} This result holds in
the general case, without the assumption of non-resonance.

In the next section we discuss the isoprincipal deformation of a
given Fuchsian system. Using our Theorem \ref{main}, we also outline
how to establish the Painlev\'e property of the Schlesinger system
and indicate possible generalizations.

Finally, in Section \ref{Examp} we illustrate the general theory
with explicit examples of the isoprincipal and  isomonodromic
deformations. In particular, we give an example of the isomonodromic
deformation which is {\em not} related to the Schlesinger system and
does {\em not} possess the Panlev\'e property. This example is based
on the theory of the isoprincipal families  of Fuchsian systems with
generic rational solutions, developed in \cite{Kats1}, \cite{Kats2}
and \cite{KaVo2}.
\section{Fuchsian differential systems}\label{Fu}
\subsection{Fuchsian differential systems}
A Fuchsian differential system is a linear system of ordinary
differential equations of the form
\begin{equation}
\label{FDE} \dfrac{dY}{dx}=\left(\sum\limits_{1\leq j\leq
n}\dfrac{Q_j}{x-t_j}\right)\,Y,
\end{equation}
where $Q_j$, $1\leq j\leq n,$ are square matrices of the same
dimension, say $Q_j\in\mathfrak{M}_k$, and $t_1,\ldots,t_n$ are
pairwise distinct points of the complex  plane $\mathbb{C}$. The
variable $x$ "lives" in the punctured Riemann sphere
$\mathbb{CP}^1\setminus\big\{t_1, \ldots , t_n\big\}$, the "unknown"
$Y$ is an \(\mathfrak{M}_k\)-valued matrix function of $x$. Under
the condition
\begin{equation}
\label{regconwt} \sum\limits_{1\leq j\leq n}Q_j=0
\end{equation}
the point $x_0=\infty$ is a {\em regular} point for the system
\eqref{FDE}. If this condition is satisfied (which we always assume
in the sequel), then in a neighborhood of the point $x_0=\infty$
there exists a fundamental solution $Y=Y(x)$ of \eqref{FDE}
satisfying the initial condition
\begin{equation}
\label{IC} Y(x)\bigm|_{x=\infty}=I.
\end{equation}
This solution $Y$ can be analytically continued   into the
multi-connected domain $\mathbb{CP}^1\setminus\{t_1, \ldots ,
t_n\}$.
 However, for $x\in\mathbb{CP}^1\setminus\{t_1, \ldots , t_n\}$ the value
of $Y$ at the point $x$ depends, in general, on the path $\alpha$
from $x_0=\infty$ to $x$, along which the analytic continuation is
performed:
$$Y=Y(x,\alpha).$$
 More precisely, $Y$ depends not on the path $\alpha$ itself, but on its homotopy
class  in $\mathbb{CP}^1\setminus\{t_1, \ldots , t_n\}$. Thus  $Y$
is a multi-valued holomorphic function in the punctured Riemann
sphere $\mathbb{CP}^1\setminus\{t_1, \ldots , t_n\}$ or, better to
say, $Y$ is a singled-valued holomorphic function on the {\em
universal covering surface} of the punctured Riemann sphere $
\mathbb{CP}^1\setminus\{t_1, \ldots , t_n\} $ with the distinguished
point $\infty$.

\subsection{Universal covering spaces}
Recall (see \cite[Chapter 1, Sections 3\,-\,5]{Fo} if need be) that
the {\em universal covering space} $\cov(\mathcal{X};x_0)$ of an
arcwise connected topological space $\mathcal{X}$ with the
distinguished point $x_0\in \mathcal{X}$ is the set of pairs
$(x,\alpha),$ where $x$ is a point in $\mathcal{X}$ and $\alpha$ is
a homotopy class of continuous mappings
$$\alpha:\{s\in\mathbb{R}\,:\,0\leq s\leq 1\}\mapsto \mathcal{X},\quad
\alpha(0)=x_0,\,\alpha(1)=x.$$ Such a mapping $\alpha$ is called a
{\em path  in} $\mathcal{X}$ {\em from} $x_0$ {\em to} $x.$ A path
in $\mathcal{X}$ from $x_0$ to $x_0$ is called a {\em loop with the
distinguished point} $x_0$.

 The product $\beta\cdot\alpha$ of
two paths $\alpha, \beta$ in $\mathcal{X}$, where $\alpha$ is a path
from $a$ to $b$ and $\beta$ is a path from $b$ to $c$, is defined as
the path from $a$ to $c$, obtained by going first along  $\alpha$
from $a$ to $b$ and then along $\beta$ from $b$ to $c$:
$$(\beta\cdot\alpha)(s)=\left\{\begin{aligned}\alpha(2s),&\quad 0\leq s\leq
\dfrac{1}{2},\\
\beta(2s-1),&\quad \dfrac{1}{2}\leq s\leq 1.
\end{aligned}
\right.
$$
 With respect to this product,
 the homotopy classes of loops in $\mathcal{X}$ with the distinguished point $x_0$ form the
 so-called {\em fundamental group} $\fug(\mathcal{X};x_0)$ of the  space $\mathcal{X}$ with
the distinguished point $x_0$.  The fundamental group
$\fug(\mathcal{X};x_0)$ acts on the universal covering space
$\cov(\mathcal{X};x_0)$ (on the right) as the group of {\em deck
transformations}
\begin{equation}\label{deck}
(x,\alpha)\mapsto(x,\alpha\gamma),\quad
(x,\alpha)\in\cov(\mathcal{X};x_0),\,\gamma\in\fug(\mathcal{X};x_0).
\end{equation}

\subsection{Monodromy}
Let $Y=Y(x,\alpha)$ be the solution of \eqref{FDE} -- \eqref{IC}
defined on  the {universal covering surface}
$\cov(\mathbb{CP}^1\setminus\{t_1, \ldots , t_n\};\,\infty). $ For
each loop $\gamma\in\fug(\mathbb{CP}^1\setminus\{t_1, \ldots ,
t_n\};\,\infty)$ let us consider the function $Y_\gamma$ defined on
$ \cov(\mathbb{CP}^1\setminus\{t_1, \ldots , t_n\};\,\infty)$ by
\begin{equation}
\label{mon1}
 Y_\gamma(x,\alpha)\defi Y(x,\alpha\gamma).\end{equation}
 The expression \eqref{mon1} means that
the value of $Y_\gamma$ at the point $(x,\alpha)$ of the universal
covering surface $ \cov(\mathbb{CP}^1\setminus\{t_1, \ldots ,
t_n\};\,\infty)$ is obtained by the analytic continuation of the
solution $Y$ of \eqref{FDE} -- \eqref{IC}: first along the loop
$\gamma$ from the distinguished point $x_0=\infty$ to itself,  then
along the path $\alpha$ from $x_0$ to $x$.

 Thus  $Y_\gamma$
 is also a  fundamental solution of the linear system \eqref{FDE} and, therefore,  there exists a unique
  invertible
constant matrix $M_\gamma\in\mathfrak{M}_k$ such that
\begin{equation}\label{mon2}
Y_\gamma(x,\alpha)\equiv Y(x,\alpha)\cdot M_\gamma,\quad
(x,\alpha)\in\cov(\mathbb{CP}^1\setminus\{t_1, \ldots ,
t_n\};\,\infty).
\end{equation}

\begin{definition}\label{DefMonMatr}
Let $Y=Y(x,\alpha)$ be the solution of \eqref{FDE} -- \eqref{IC} on
the universal covering surface $ \cov(\mathbb{CP}^1\setminus\{t_1,
\ldots , t_n\};\infty)$ and let
$\gamma\in\fug(\mathbb{CP}^1\setminus\{t_1, \ldots ,
t_n\};\,\infty)$.

 The constant (with respect to $x$) matrix $M_\gamma\in\mathfrak{M}_k$, which appears in  the identity
 \eqref{mon2},
 is said to be the {\em monodromy matrix\,} of the solution $Y$, corresponding
to the loop $\gamma$.
\end{definition}

Note that for a pair of loops
$\gamma_1,\gamma_2\in\fug(\mathbb{CP}^1\setminus\{t_1, \ldots ,
t_n\})$ and the cooresponding monodromy matrices
$M_{\gamma_1},M_{\gamma_2}$ of the solution $Y$ it holds that
$$Y(x,\alpha\gamma_1\gamma_2)= \left(Y\cdot M_{\gamma_2}\right)({x},\alpha\gamma_1)
=Y(x,\alpha\gamma_1)\cdot M_{\gamma_2}=Y(x,\alpha)\cdot
M_{\gamma_1}M_{\gamma_2}.$$ Therefore, the monodromy matrices of $Y$
satisfy the following multiplicative identity:
\begin{equation}
\label{linrep} M_{\gamma_1\gamma_2}=M_{\gamma_1}M_{\gamma_2}\quad
\forall\gamma_1,\gamma_2\in\fug(\mathbb{CP}^1\setminus\{t_1, \ldots
, t_n\}).
\end{equation}
This means that the mapping $\gamma\mapsto M_\gamma$ is a linear
representation of the fundamental group
$\fug(\mathbb{CP}^1\setminus\{t_1, \ldots , t_n\};\infty).$

\begin{definition}
Let $Y$ be the solution of \eqref{FDE} -- \eqref{IC} on the
universal covering surface $ \cov(\mathbb{CP}^1\setminus\{t_1,
\ldots , t_n\};\infty).$ The linear representation of the
fundamental group $\fug(\mathbb{CP}^1\setminus\{t_1, \ldots ,
t_n\};\infty)$
\begin{equation}\label{monrep}
\gamma\mapsto M_\gamma,\quad
\gamma\in\fug(\mathbb{CP}^1\setminus\{t_1, \ldots ,
t_n\};\infty),\end{equation} where $M_\gamma$ denotes the monodromy
matrix of the solution $Y$, corresponding to the loop $\gamma$,  is
called the {\em monodromy representation} of the solution $Y$.
\end{definition}

\subsection{The regular-principal factorization for a
fundamental solution of a Fuchsian system: single-valued case} Each
of the points $t_j,$ $ 1\leq j\leq n,$ is a singularity of the
solution $Y$. This means that at least one of the two functions $Y$
and $Y^{-1}$ is not holomorphic at $t_j$. More information about the
nature of the singularity at $t_j$ can be obtained from a certain
multiplicative decomposition of the solution $Y$ near the point
$t_j$, which is called the {\em regular-principal factorization.}

In order to explain the idea of the regular-principal factorization,
let us assume for the moment that the solution $Y=Y(x)$ is {\em
single-valued} in the domain $\mathbb{CP}^1\setminus\{t_1, \ldots ,
t_n\}$ -- that is, the monodromy representation of $Y$ is trivial:
$$M_\gamma=I\quad\forall \gamma\in\fug(\mathbb{CP}^1\setminus\{t_1, \ldots ,
t_n\};\infty).$$
 For $1\leq j\leq n$ let $\mathcal{V}_j$ be an open
simply connected neighborhood of $t_j$ in $\mathbb{C}$, such that
$t_k\not\in\mathcal{V}_j$ for $k\not=j$. Then it follows, for
instance, from G.D. Birkhoff's results on factorization of matrix
functions holomorphic  in the annulus (see \cite[\S 7]{Birk1}) that
in the punctured neighborhood $\mathcal{V}_j\setminus \{t_j\}$ the
solution $Y(x)$ admits a factorization of the form
$$Y(x)=H_j(x)\cdot P_j(x),\quad x\in\mathcal{V}_j\setminus \{t_j\},$$
where the function $H_j(x)$ is holomorphic and invertible in the
entire (non-punctured) neighborhood $\mathcal{V}_j$ and   the
function $P_j(x)$ is holomorphic, single-valued and invertible in
the punctured plane $\mathbb{C}\setminus \{t_j\}$. The factors
$H_{j}(x)$ and $P_{j}(x)$ are
 said to be, respectively, the {\em regular factor} and the {\em principal
factor} of the solution $Y$ at its singular point $t_j$.

In the general case, when the monodromy representation of $Y$ may be
non-trivial, the regular-principal factorization of $Y$ is more
involved.

Indeed, on the one hand the solution $Y$ is normalized at the
distinguished point $x_0=\infty$. In order to consider $Y$ in a
neighborhood of $t_j$, we have to choose a homotopy class of paths,
connecting the distinguished point $x_0=\infty$ with this
neighborhood of $t_j$, and such a choice is not unique.

On the other hand, even in the single-valued case $x_0=\infty$ is,
in general, a singular point of the principal factor $P_{j}$. In the
general case $P_j$ will have to be considered as a function on a
universal covering surface of the punctured plane
$\mathbb{C}\setminus\{t_j\}$. Therefore, $P_j$ needs to be
normalized at  some distinguished point in
$\mathbb{C}\setminus\{t_j\}$ (obviously different from $x_0=\infty$)
and continued analytically from there.

Thus, before we can present the regular-principal factorization of
$Y$ in the general case, we need some preparation.

\subsection{Branches of the solution of a Fuchsian system in a neighborhood of the singular point}
We propose the following terminology:
\begin{definition}\label{defclift}
Let $Y$ be the solution of \eqref{FDE} -- \eqref{IC} on the
universal covering surface $ \cov(\mathbb{CP}^1\setminus\{t_1,
\ldots , t_n\};\,\infty).$

   For $1\leq j\leq n$  assume that:
  \begin{enumerate}
  \item
   $\mathcal{V}$ is a domain in $ \mathbb{CP}^1\setminus\{t_1,
\ldots , t_n\}$;
\item $p$ is a point in the domain
$\mathcal{V}$;
 \item $\alpha$ is a
path in $\mathbb{CP}^1\setminus\{t_1, \ldots , t_n\}$ from the
distinguished point $x_0=\infty$ to the point $p$.
\end{enumerate}

 Define the function $Y_{\alpha}$ on the universal covering
surface $\cov(\mathcal{V};\,p)$ by the analytic continuation of the
solution $Y$ first along the path $\alpha$, then inside the domain
$\mathcal{V}$:
\begin{equation}\label{defbranch}Y_{\alpha}(x,\beta)\defi
Y(x,\beta\cdot\alpha),\end{equation}
 where
$x\in\mathcal{V}$ and $\beta$ is a path in $\mathcal{V}$ from $p$
to $x$.

Then the function $Y_{\alpha}$, holomorphic in
$\cov(\mathcal{V};\,p)$, is said to be {\em the branch of the
solution}\, $Y$ {\em in the domain} $\mathcal{V}$, {\em
corresponding to the path} $\alpha$.
\end{definition}

In what follows we shall be mostly dealing with the branches of the
solution $Y$ in domains of the form $\mathcal{V}_j\setminus\{t_j\}$,
where  $\mathcal{V}_j$ is a simply connected neighborhood of the
singular point $t_j$, such that $t_k\not\in\mathcal{V}_j$ for
$k\not=j$. In this case
  the universal covering surface
$\cov(\mathcal{V}_j\setminus\{t_j\};p_j)$ has a simple structure:
the fundamental group  $\cov(\mathcal{V}_j\setminus\{t_j\};p_j)$ is
cyclic, generated by the loop  with the distinguished point $p_j$
which makes one positive circuit of  $t_j$ in
$\mathcal{V}_j\setminus\{t_j\}$.

\begin{definition}\label{loops}
 For $1\leq j\leq n$  assume that:
  \begin{enumerate}
  \item
   $\mathcal{V}_j\subset\mathbb{C}$ is an open simply
connected neighborhood of the singular point $t_j$, such that
$t_k\not\in\mathcal{V}_j$ for $k\not=j$; \item $p_j$ is a point in
the punctured neighborhood $\mathcal{V}_j\setminus\{t_j\}$;
 \item $\alpha_j$ is a
path in $\mathbb{CP}^1\setminus\{t_1, \ldots , t_n\}$ from the
distinguished point $x_0=\infty$ to the point $p_j$.
\end{enumerate}

Let $\beta_j$ be  the loop  in the punctured neighborhood
$\mathcal{V}_j\setminus\{t_j\}$ with the distinguished point $p_j$
which makes one positive circuit of $t_j$,  and let $\gamma_j$ be
the loop in the punctured sphere
$\mathbb{CP}^1\setminus\{t_1,\dots,t_n\}$ with the distinguished
point $x_0=\infty $, defined by
\begin{equation}
\label{bigloop} \gamma_j\defi\alpha_j^{-1}\cdot\beta_j\cdot\alpha_j
\end{equation}
(that is, the loop $\gamma_j$   goes from $x_0=\infty$ to $p_j$
along $\alpha_j$, then makes one positive circuit  of $t_j$ along
the  small loop $\beta_j$, then
 goes again along $\alpha_j$, but in the opposite direction: from
$p_j$ to $x_0=\infty$).

Then: \begin{itemize} \item the loop $\beta_j$ is said to be {\em
the small loop around} $t_j$ {\em in the punctured neighborhood}
$\mathcal{V}_j\setminus\{t_j\}$; \item  the loop $\gamma_j$ is said
to be {\em the big loop around} $t_j$, {\em corresponding to the
path} $\alpha_j$.
\end{itemize}
\end{definition}

\begin{remark}\label{genmon}
Note that for a suitable choice of the paths
$\alpha_1,\dots,\alpha_n$ the corresponding big loops
$\gamma_1,\dots,\gamma_n$ generate the fundamental group
$\fug(\mathbb{CP}^1\setminus\{t_1,\dots,t_n\};\,\infty)$. These
generators are not free:  choosing $\alpha_1,\dots,\alpha_n$
carefully we can ensure, for example, that
$\gamma_1\cdots\gamma_n=1.$
\end{remark}

We observe that the surface
$\cov(\mathcal{V}_j\setminus\{t_j\};p_j)$ in Definition
\ref{defclift} is naturally embedded into the universal covering
surface $\cov(\mathbb{C}\setminus\{t_j\};p_j)$, which is isomorphic
to the {\em Riemann surface of the logarithm} $\ln\zeta$.

Although the basic properties of the function $\ln\zeta$ are very
well-known, we shall discuss them in some detail, because they  are
important for our future considerations.

\subsection{The Riemann surface of $\ln\zeta$}
For each fixed $\zeta\in\mathbb{C}\setminus\{0\}$ the equation
$$e^\lambda=\zeta$$ has a countable set of solutions $\lambda=\lambda(\zeta)\defi\ln\zeta$.
These solutions can be parameterized as
\begin{equation}\label{deflog}\ln\zeta=\ln|\zeta|+i\arg\zeta,\text{ where }
\ln|\zeta|\in\mathbb{R}\end{equation}
 and $\arg\zeta$ is an
equivalence class of real numbers  (the {\em values} of $\arg\zeta$)
modulo addition by $2\pi$.

Now we explain how  the function $\ln\zeta$ can be defined as a
single-valued holomorphic function on the universal covering surface
$\cov(\mathbb{C}\setminus\{0\};1)$ of the punctured plane
$\mathbb{C}\setminus\{0\}$ with the distinguished point $\zeta_0=1$.

 Let us choose a point $\zeta\in\mathbb{C}\setminus\{0\}$ and a
path $\vartheta$ in $\mathbb{C}\setminus\{0\}$ from $\zeta_0=1$ to
$\zeta$. Then there exists  a unique $\theta\in\mathbb{R}$ such
that, up to homotopy  in $\mathbb{C}\setminus\{0\}$, the path
$\vartheta$ can be parameterized as follows:
\begin{equation}\label{defarg}\vartheta(s)=e^{(\ln|\zeta|+i\theta)s},
\quad  0\leq s\leq 1\end{equation} (here $\ln|\zeta|$ is the
real-valued logarithm).
  The real
number $\theta$ is said to be {\em the value of}  $\arg\zeta$ {\em
corresponding to the path} $\vartheta$. In this manner we establish
a $1$-to-$1$ correspondence between the values of $\arg\zeta$ and
the homotopy classes of paths in $\mathbb{C}\setminus\{0\}$ from
$\zeta_0=1$ to $\zeta$.

Thus  the function $\arg\zeta$ is defined as a single-valued
continuous function on $\cov(\mathbb{C}\setminus\{0\};1)$.
Accordingly, the function $\ln\zeta$ is defined by \eqref{deflog} as
a single-valued holomorphic function on the universal covering
surface $\cov(\mathbb{C}\setminus\{0\};1)$. In the sequel we shall
often refer to the surface $\cov(\mathbb{C}\setminus\{0\};1)$ as
{\em the Riemann surface of} $\ln\zeta$.

The fundamental group $\fug(\mathbb{C}\setminus\{0\};1)$ is cyclic,
generated by the loop  with the distinguished point $\zeta_0=1$
which makes one turn counterclockwise around the origin. The
corresponding deck transformation of
$\cov(\mathbb{C}\setminus\{0\};1)$ is denoted by \begin{equation}
\label{dtrrsl} \zeta\mapsto \zeta\cdot e^{2\pi i},\end{equation} so
that the following monodromy relations hold:
\begin{equation}\label{monlog}\arg(\zeta e^{2\pi i})=\arg\zeta+2\pi,\quad
\ln(\zeta e^{2\pi i})=\ln\zeta +2\pi i.\end{equation}

\subsection{Transplants of functions defined on the Riemann surface of $\ln\zeta$}
Let us   consider the universal covering surface
$\cov(\mathbb{C}\setminus\{t_j\};p_j)$, where $t_j$ is some point in
the complex plane $\mathbb{C}$ and $p_j$ is a distinguished point in
the punctured plane $\mathbb{C}\setminus\{t_j\}$. Let us choose
 some value $\theta_j$ of $\arg(p_j-t_j)$ and let
$\vartheta_j$ denote the corresponding path (see \eqref{defarg}) in
$\mathbb{C}\setminus\{0\}$ from $\zeta_0=1$ to $\zeta=p_j-t_j$. Let
us define a mapping from the universal covering surface
$\cov(\mathbb{C}\setminus\{t_j\};p_j)$ into the Riemann surface of
$\ln\zeta$ as follows.

To each point $(x,\alpha)\in\cov(\mathbb{C}\setminus\{t_j\};p_j),$
where $\alpha$ is a path in $\mathbb{C}\setminus\{t_j\}$ from $p_j$
to $x$, we associate the point
$(x-t_j,\alpha_{t_j}\cdot\vartheta_j)\in\cov(\mathbb{C}\setminus\{0\};1)$,
where the path $\alpha_{t_j}$, which leads in
$\mathbb{C}\setminus\{0\}$ from $p_j-t_j$ to $x-t_j$, is obtained by
the parallel translation of the path $\alpha$:
$$\alpha_{t_j}(s)=\alpha(s)-t_j,\quad 0\leq s\leq
1.$$

This mapping  is   an isomorphism between
$\cov(\mathbb{C}\setminus\{t_j\};p_j)$ and the Riemann surface of
$\ln\zeta$. It will be denoted by
\begin{equation}\label{defisotra}
x\xrightarrow[\arg(p_j-t_j)=\theta_j]{} x-t_j,\quad
x\in\cov(\mathbb{C}\setminus\{t_j\};p_j),\
x-t_j\in\cov(\mathbb{C}\setminus\{0\};1),
\end{equation}
and the inverse mapping will be denoted by
\begin{equation}\label{defisotrainv}
x\xrightarrow[\arg(p_j-t_j)=\theta_j]{} x+t_j,\quad
x\in\cov(\mathbb{C}\setminus\{0\};1),\
x+t_j\in\cov(\mathbb{C}\setminus\{t_j\};p_j).
\end{equation}
Accordingly, we shall denote the deck transformation of
$\cov(\mathbb{C}\setminus\{t_j\};p_j)$, corresponding to the loop
with the distinguished point $p_j$ which makes one positive circuit
around $t_j$ in $\mathbb{C}\setminus\{t_j\}$, by
\begin{equation}\label{smalltr}
x\mapsto t_j+(x-t_j)\cdot e^{2\pi i},\quad
x\in\cov(\mathbb{C}\setminus\{t_j\};p_j).
\end{equation}

With this notation we can consider the function $\ln(x-t_j)$ as a
function of $x$,  holomorphic on the universal covering surface $
\cov(\mathbb{C}\setminus\{t_j\};p_j)$ and such that (see
\eqref{monlog})
\begin{equation}\label{monlogt}
\ln\left((x-t_j)\cdot e^{2\pi i}\right)=\ln(x-t_j)+2\pi i.
\end{equation}

 More generally,
we propose the following

\begin{definition}\label{DefTrans}
Let a function $E(\zeta)$ be defined on the Riemann surface of
$\ln\zeta$. Let $t_j\in\mathbb{C}$, let
$p_j\in\mathbb{C}\setminus\{t_j\}$ and let us choose a value
$\theta_j$ of $\arg(p_j-t_j)$. For each
$x\in\cov(\mathbb{C}\setminus\{t_j\};p_j)$ let $x-t_j$ denote the
image of $x$ in the Riemann surface of $\ln\zeta$ under the
isomorphism \eqref{defisotra}.

The function $E(x-t_j)$, defined as a function of $x$ on the
universal covering surface $ \cov(\mathbb{C}\setminus\{t_j\};p_j)$,
is said to be the {\em transplant of the function} $E(\zeta)$ {\em
into} $ \cov(\mathbb{C}\setminus\{t_j\};p_j)$, {\em corresponding to
the value} $\theta_j$ {\em of} $\arg(p_j-t_j)$.
\end{definition}

\subsection{The regular-principal factorization for a fundamental solution of a Fuchsian system: general case}
\begin{theorem}
\label{RePrDec} Let $Y$   be the solution
  of the Fuchsian system \eqref{FDE} -- \eqref{regconwt}, satisfying the initial condition \eqref{IC}.
   For $j=1,\dots,n$   assume that:\begin{enumerate}
  \item
   $\mathcal{V}_j\subset\mathbb{C}$ is an open simply
connected neighborhood of the singular point $t_j$, such that
$t_k\not\in\mathcal{V}_j$ for $k\not=j$; \item $p_j$ is a point in
the punctured neighborhood $\mathcal{V}_j\setminus\{t_j\}$;
 \item $\alpha_j$ is a
path in $\mathbb{CP}^1\setminus\{t_1, \ldots , t_n\}$ from the
distinguished point $x_0=\infty$ to the point $p_j$.
\end{enumerate}

 Then for each $j$, $1\leq j\leq n,$ the branch $Y_{\alpha_j}$ of the solution
$Y$ in the punctured neighborhood $\mathcal{V}_j\setminus\{t_j\}$,
corresponding to the path $\alpha_j$,
 admits  the factorization
\begin{equation}
\label{RPD}%
 Y_{\alpha_j}(x)=H_{\alpha_j}(x)\cdot P_{\alpha_j}(x), \quad
 x\in\cov(\mathcal{V}_j\setminus\{t_j\};\,p_j),
\end{equation}
where the factors possess the following properties:
\begin{itemize}
\item[(R)]
 the matrix function $H_{\alpha_j}(x)$ is holomorphic and invertible in the
 entire\,\footnote{\label{SiVa}In particular, the function $H_{\alpha_j}(x)$ is
single-valued in $\mathcal{V}_j$.} (non-punctured) neighborhood
$\mathcal{V}_j$;
\item[(P)]
 the matrix function $P_{\alpha_j}( x)$ is  holomorphic and invertible on the universal covering surface
 $\cov(\mathbb{C}\setminus\{t_j\};p_j)$.
 \end{itemize}
\end{theorem}

\begin{definition}
\label{DefFact}%
The factorization \eqref{RPD}, where the factors  $H_{\alpha_j}$ and
$P_{\alpha_j}$ possess the properties (R) and (P), is said to be
{\em the regular-principal factorization} of the branch
$Y_{\alpha_j}$ of the solution $Y$  in a punctured  neighborhood of
the singular point $t_j.$ The factors  $H_{\alpha_j}$ and
$P_{\alpha_j}$ are
 said to be, respectively, {\em the regular factor} and {\em the principal
factor} of the branch $Y_{\alpha_j}$.
\end{definition}

\begin{proof}[Proof of Theorem \ref{RePrDec}]
For $j=1,\dots,n$ let $\gamma_j$ be the big loop around $t_j$,
corresponding to the path $\alpha_j$ (see Definition \ref{loops}),
and let $M_{\gamma_j}$ be the corresponding monodromy matrix of $Y$.

Then, according to  Definitions \ref{DefMonMatr} and \ref{defclift},
the monodromy matrix $M_{\gamma_j}$ is given by

\begin{equation}\label{mongen}
M_{\gamma_j}=Y_{\alpha_j}^{-1}(x)Y_{\alpha_j}(t_j+(x-t_j)\cdot
e^{2\pi i}),\quad \forall
x\in\cov(\mathcal{V}_j\setminus\{t_j\};\,p_j),
\end{equation}
where $Y_{\alpha_j}$ is the  branch of $Y$ in
$\mathcal{V}_j\setminus\{t_j\}$, corresponding to the path
$\alpha_j$.

Since the matrix $M_{\gamma_j}$ is invertible,
 there exists a matrix  denoted by $\ln M_{\gamma_j}$, such that\,\footnote{Here
  we refer to \cite[Chapter VIII, Section 8]{Gant}. Such a matrix  $\ln
 M_{\gamma_j}$ is not unique, of course, but for our purposes any
 choice of  $\ln
 M_{\gamma_j}$ will do.}
 $$e^{\ln M_{\gamma_j}}=M_{\gamma_j}.$$

Let us choose a transplant $\ln(x-t_j)$ of the function $\ln\zeta$
into $\cov(\mathbb{C}\setminus\{t_j\};p_j)$. Then, according to
\eqref{monlogt},
 the matrix
function
$$ (x-t_j)^{\dfrac{1}{2\pi i}\ln M_{\gamma_j}}\defi e^{\ln(
x-t_j)\dfrac{1}{2\pi i}\ln M_{\gamma_j}},$$ which is holomorphic and
invertible on $\cov(\mathbb{C}\setminus\{t_j\};p_j)$,  satisfies the
relation
$$((x-t_j)\cdot e^{2\pi i})^{\dfrac{1}{2\pi i}\ln M_{\gamma_j}}=(x-t_j)^{\dfrac{1}{2\pi i}\ln
M_{\gamma_j}}\cdot M_{\gamma_j},\quad
x\in\cov(\mathbb{C}\setminus\{t_j\};p_j).
$$
Hence, in view of \eqref{mongen}, the branch $Y_{\alpha_j}$  of the
solution $Y$ in $\mathcal{V}_j\setminus\{t_j\}$ has the form
$$Y_{\alpha_j}( x)=\Phi(x)\cdot ( x-t_j)^{\dfrac{1}{2\pi i}\ln
M_{\gamma_j}},\quad x\in\cov(\mathbb{C}\setminus\{t_j\};p_j),$$
where $\Phi(x)$ is a matrix function,  holomorphic, invertible and
{\em single-valued} in the punctured neighborhood
$\mathcal{V}_{j}\setminus \{t_j\}$.

Now, according to \cite[\S 7]{Birk1}, we can factorize the function
$\Phi(x)$ as
$$\Phi(x)=\Phi_+(x)\cdot\Phi_-(x),\quad x\in\mathcal{V}_j\setminus \{t_j\}$$
where $\Phi_+(x)$ and $\Phi_-(x)$ are matrix functions,
single-valued, holomorphic and invertible in, respectively, the
entire (non-punctured) neighborhood $\mathcal{V}_j$  and the
punctured plane $\mathbb{C}\setminus \{t_j\}$. We set
 \begin{align*}
 H_{\alpha_j}(x)&\defi\Phi_+(x),\\
 P_{\alpha_j}( x)&\defi\Phi_-(x)\cdot( x-t_j)^{\dfrac{1}{2\pi i}\ln
M_{\gamma_j}} \end{align*}   and obtain the desired factorization
\eqref{RPD} with the properties (R) and (P). This completes the
proof. \end{proof}

\begin{remark}\label{regau}
 Of course, the principal and regular factors of the branch $Y_{\alpha_j}$ of the solution $Y$ in a
punctured neighborhood of the singularity $t_j$ are determined only
up to the transformation
\begin{equation}%
\label{gauprin}%
P_{\alpha_j}(x)\to T(x)\cdot P_{\alpha_j}( x), \quad
H_{\alpha_j}(x)\to T(x)^{-1}\cdot H_{\alpha_j}(x),
\end{equation}
 where $T(x)$ is an  invertible entire
matrix function. However, once the choice of, say, the regular
factor $H_{\alpha_j}$ is fixed,
 the principal factor $P_{\alpha_j}$ is
uniquely determined.

 Moreover,  if we choose a different path from
$x_0$ to $p_j$, say $\alpha_j^\prime$, then, in view of Definitions
\ref{DefMonMatr} and \ref{defclift}, the branches $Y_{\alpha_j}$ and
$Y_{\alpha_j^\prime}$ of the solution $Y$ are related by
$$Y_{\alpha_j^\prime}(x)=Y_{\alpha_j}(x)\cdot
M_{\alpha_j^{-1}\alpha_j^\prime},\quad
x\in\cov(\mathbb{C}\setminus\{t_j\};p_j),$$ where
$M_{\alpha_j^{-1}\alpha_j^\prime}$ is the  monodromy matrix of $Y,$
corresponding to the loop
$\alpha_j^{-1}\alpha_j^\prime\in\fug(\mathbb{CP}^1\setminus\{t_1,
\ldots , t_n\};\infty)$. Hence  the branch $Y_{\alpha_j^\prime}$
admits  the regular-principal factorization
$$Y_{\alpha_j^\prime}( x)=H_{\alpha_j^\prime}(x)\cdot P_{\alpha_j^\prime}( x), \quad
 x\in\cov(\mathbb{C}\setminus\{t_j\};p_j),$$
 where
\begin{align*}
H_{\alpha_j}(x)&=H_{\alpha_j^\prime}(x),&\quad x&\in\mathcal{V}_j, \\
P_{\alpha_j^\prime}( x)&=P_{\alpha_j}(x) \cdot
M_{\alpha_j^{-1}\alpha_j^\prime},&\quad
{x}&\in\cov(\mathbb{C}\setminus\{t_j\};p_j).\end{align*}

  Thus the regular factor  at the singular point $t_j$ can be chosen {\em
independently} of the choice of the path $\alpha_j$ and  will be
denoted simply by $H_j(x)$.
\end{remark}

\begin{remark}\label{wlog1}
\label{PosMD} Up to now, we have made no use of the fact that  the
system \eqref{FDE} is Fuchsian (that is, each singularity is a {\em
simple pole} for the coefficients of the system). In particular, the
regular - principal factorization \eqref{RPD} of the fundamental
solution of a linear differential system also takes place in a
neighborhood of the isolated singularity,  where the coefficients of
the system  have a higher order pole  or even an essential singular
point.

However, in the special case when the system is {\em Fuchsian}, the
general form of the fundamental solution in a neighborhood of its
singular point is quite well known (see, for instance, \cite[Chapter
XV, \S 10]{Gant}) and thus much more precise statements concerning
the principal factors of the solution of the system \eqref{FDE} can
be made.

If the matrix $Q_j$ is non-resonant\,\footnote{\label{DNRM} A square
matrix $Q$ is said to be \textsl{non-resonant} if distinct
eigenvalues of $Q$ do not differ by integers or, in other words, if
the spectra of the matrices $Q+nI$ and $Q$ are disjoint for every
$n\in\mathbb{Z}\setminus 0$.}, then the principal factor
$P_{\alpha_j}$ can be chosen in the form
\begin{equation}\label{prinnonres}
P_{\alpha_j}( x)=( x-t_j)^{A_{\alpha_j}},\end{equation}
 where
$A_{\alpha_j}$ is a matrix, similar to the matrix $Q_j$:
 \begin{equation}\label{similar}
 A_{\alpha_j}=C_{\alpha_j}^{-1}Q_jC_{\alpha_j}, \end{equation}
  where
$C_{\alpha_j}$ is an invertible matrix. In the general case (without
the assumption that the matrix $Q_j$ is non-resonant) the principal
factor  $P_{\alpha_j}$ can be chosen in the form
\begin{equation}
\label{SPF}%
 P_{\alpha_j}( x)=(x-t_j)^{Z_{\alpha_j}}\cdot ( x-t_j)^{ A_{\alpha_j}},
\end{equation}
where  $Z_{\alpha_j}$ is a diagonalizable matrix with integer
eigenvalues
 $l_1, \dots , l_k$ and
   $  A_{\alpha_j}$ is a non-resonant matrix, whose
  eigenvalues $\hat{\lambda}_p$ are related
 to the eigenvalues
$\lambda_p$ of the matrix $Q_j$ by the equations
\begin{equation}\label{specla} \hat{\lambda}_p=\lambda_p-l_p,\quad 1\leq p\leq
k.\end{equation}
 The matrices $Z_{\alpha_j}, A_{\alpha_j}$ also
possess certain additional properties,  but we shall not go into
further details, because in the sequel our considerations  will be
mostly based on the existence of the regular - principal
factorization \eqref{RPD} described in Theorem \ref{RePrDec} rather
than on the specific form of the factors.
\end{remark}

\section{Holomorphic Families of Fuchsian Differential Systems} \label{DefFDE}
\subsection{Families of Fuchsian systems, parameterized by the pole
loci} In the present paper we consider a {\em family} of linear
 differential systems of the form
\begin{equation}
\label{FLDEt}%
\dfrac{dY}{dx}=\bigg(\sum\limits_{1\leq j\leq
n}\dfrac{Q_j(\boldsymbol{t})}{x-t_j}\bigg) Y.
\end{equation}
The variable $x$ "lives" in  the punctured Riemann sphere
$\mathbb{CP}^1\setminus\{t_1, \ldots , t_n\}$, where $t_1, \ldots ,
t_n$ are pairwise distinct points of the complex plane $\mathbb{C}$.
However, now  $t_1, \ldots , t_n$ are not fixed but serve as the
parameters of the family. The string $\boldsymbol{t}=(t_1, \ldots ,
t_n)$ is considered as a point of $\mathbb{C}^n_{\ast}$ and the
residue matrices $Q_j$  are  assumed to  depend  on the the
parameter $\boldsymbol{t}$.  The "unknown" $Y$ is a square matrix
function  depending both on the "main variable" $x$ and on the
parameter $\boldsymbol{t}$.

\begin{definition}
\label{DHD}%
Assume that the matrix functions $Q_j(\boldsymbol{t}),$  $1\leq
j\leq n,$ are defined  and holomorphic for
$\boldsymbol{t}\in\boldsymbol{\mathcal{D}}$, where
$\boldsymbol{\mathcal{D}}$ is a domain in $\mathbb{C}^n_{\ast}.$
Then the family \eqref{FLDEt} is said to be a {\em holomorphic
family of Fuchsian systems, parameterized by the pole loci.}
\end{definition}

In what follows we assume that the condition
\begin{equation}%
\label{IRSP}%
\sum\limits_{1\leq j\leq n}Q_j(\boldsymbol{t})\equiv 0,\quad
\boldsymbol{t}\in\boldsymbol{\mathcal{D}},
\end{equation}%
holds. Thus for every fixed
$\boldsymbol{t}\in\boldsymbol{\mathcal{D}}$ the point $x=\infty$ is
a regular point for the system \eqref{FLDEt}, considered as a
differential system with respect to $x$.

For each fixed $\boldsymbol{t}\in\boldsymbol{\mathcal{D}}$, the
fundamental solution $Y=Y(x,\boldsymbol{t})$ of the differential
system \eqref{FLDEt} with the initial condition\footnote{In view of
\eqref{IRSP}, for each fixed
$\boldsymbol{t}\in\boldsymbol{\mathcal{D}}$ the point $x=\infty$ is
regular for the system \eqref{FLDEt} and hence the initial condition
\eqref{ICt} can be posed.}
\begin{equation}
\label{ICt} Y( x,\boldsymbol{t})\bigm|_{x=\infty}=I
\end{equation}%
  is defined and holomorphic as a function of $x$ on the
universal covering surface $\cov(\mathbb{CP}^1\setminus\{t_1, \ldots
, t_n\};\infty)$.

 In the present section our goal is to compare the properties of
$Y({x},\boldsymbol{t})$, such as the monodromy representation or the
principal factors, for different $\boldsymbol{t}$. More precisely,
we would like to understand what does it mean that "the monodromy
representation or the principal factors of $Y({x},\boldsymbol{t})$
{are the same for different} $\boldsymbol{t}$"? In the previous
section these notions were defined in terms of homotopy classes of
paths  on the punctured Riemann sphere $\mathbb{CP}^1\setminus\{t_1,
\dots , t_n\}$. However, for different $\boldsymbol{t}=(t_1, \dots ,
t_n)$ the domains $\mathbb{CP}^1\setminus\{t_1, \dots , t_n\}$ are
different. Thus one ought to explain how to consider "the same paths
for different $\boldsymbol{t}$". This can be done because we can
confine ourselves to {\em local} considerations.

\subsection{Cylindrical neighborhoods of points in $\mathbb{C}^n$}

\begin{definition}\label{defcyl}
Let $\mathcal{W}_j,$ $1\leq j\leq n$, be subsets of the complex
plane $\mathbb{C}$.

The Cartesian product
\begin{equation*}
\boldsymbol{\mathcal{W}}=\mathcal{W}_1\times \cdots \times
\mathcal{W}_n\subset\mathbb{C}^n
\end{equation*}
is said to be  {\em the cylindrical set with the  bases}
$\mathcal{W}_j,$ $1\leq j\leq n$.
\end{definition}

\begin{definition}
Let   $\boldsymbol{t}=(t_1, \ldots , t_n)$ be a point in
$\mathbb{C}^n$. For $j=1, \ldots , n$ let $\mathcal{W}_j$ be an
 open neighborhood of $t_j$ in $\mathbb{C}$, such that:
 \begin{enumerate}
 \item the set $\mathcal{W}_j$ is simply connected;
 \item the set $\mathbb{C}\setminus\overline{\mathcal{W}_j}$ is
 connected.
 \end{enumerate}
  Denote by
$ \boldsymbol{\mathcal{W}}$  the cylindrical set with the bases
$\mathcal{W}_j$, $1\leq j\leq n$.

The set $\boldsymbol{\mathcal{W}}$ is said to be an {\em open
cylindrical neighborhood of the point} $\boldsymbol{t}$ in
$\mathbb{C}^n$.
\end{definition}

\begin{definition}
Let subsets $\boldsymbol{\mathcal{S}},\boldsymbol{\mathcal{O}}$ of
$\mathbb{C}^n$ be such that: \begin{enumerate} \item the set
$\boldsymbol{\mathcal{O}}$ is open; \item the closure
$\overline{\boldsymbol{\mathcal{S}}}$ is compact; \item
$\overline{\boldsymbol{\mathcal{S}}}\subset\boldsymbol{\mathcal{O}}$.
\end{enumerate}

Then we say that the set $\boldsymbol{\mathcal{S}}$ is {\em
compactly included} in $\boldsymbol{\mathcal{O}}$ and denote this
relation by
$$
\boldsymbol{\mathcal{S}}\Subset\boldsymbol{\mathcal{O}}.
$$
\end{definition}

\subsection{Isomonodromic families of Fuchsian systems}

 Let $\boldsymbol{t}^0$ be a point in the domain\footnote{ Recall
that $\boldsymbol{\mathcal{D}}$ is the domain where the residue
matrices $Q_j(\boldsymbol{t})$ from \eqref{FLDEt} are defined and
holomorphic. Since
 $\boldsymbol{\mathcal{D}}\subseteq\mathbb{C}^n_{\ast}$,
the coordinates $t_1^0,t_2^0,\dots,t_n^0$ of every point
$\boldsymbol{t^0}\in\boldsymbol{\mathcal{D}}$ are pairwise
distinct.} $\boldsymbol{\mathcal{D}}$ and let
$\boldsymbol{\mathcal{W}}$ be a cylindrical neighborhood of\,
$\boldsymbol{t}^0$, such that
$\boldsymbol{\mathcal{W}}\Subset\boldsymbol{\mathcal{D}}$.  Then
 \begin{equation}%
 \label{DS}%
 \overline{\mathcal{W}_p}\cap\overline{\mathcal{W}_q}=\emptyset.
\ \ 1\leq p,q\leq n, p\not=q,\end{equation} hence\footnote{Here we
refer to a relatively delicate result from general topology:
\textsl{if $K_1,K_2$ are two disjoint compact subsets of
$\mathbb{R}^m$ and each of
two sets %
$\mathbb{R}^m\setminus \!K_1$,  $\mathbb{R}^m\setminus\!K_2$ %
is connected, then the set $\mathbb{R}^m\setminus \!\{K_1\cup K_2\}$
is connected, as well.} (See for example \cite{HW}, Corollary of
Theorem VI.10.) Actually we do not need  this result in full
generality. For our goal it is enough to consider only  the case of
$m=2$ and some very special compact sets $K\subset\mathbb{R}^2$,
such as finite unions of disks, etc. In this particular case, the
above stated result is elementary.} the set
$\mathbb{C}\setminus\bigcup_k\overline{\mathcal{W}}_k$ is connected.

For a fixed $\boldsymbol{t}\in\boldsymbol{\mathcal{W}}$ each
homotopy class of
 loops in the {\em punctured} sphere
$\mathbb{CP}^1\setminus\{t_1,\dots,t_n\}$ with the distinguished
point $x_0=\infty$ has a representative which is  a loop in the {\em
perforated} sphere
$\mathbb{CP}^1\setminus\bigcup_k\overline{\mathcal{W}}_k$.

On the other hand, if $\gamma$ is  a loop in the perforated sphere
$\mathbb{CP}^1\setminus\bigcup_k\overline{\mathcal{W}}_k$ with the
distinguished point $x_0=\infty$,   then for each fixed
$\boldsymbol{t}\in\boldsymbol{\mathcal{W}}$ the loop $\gamma$
 serves as a path of the
analytic continuation  with respect to $x$ for the solution
$Y(x,\boldsymbol{t})$ of \eqref{FLDEt} -- \eqref{ICt} and one can
consider the corresponding monodromy matrix\,\footnote{See
Definition \ref{DefMonMatr}.} $M_\gamma(\boldsymbol{t}).$

Although the loop $\gamma$ does not depend on $\boldsymbol{t}$,
 the corresponding monodromy matrix
$M_\gamma(\boldsymbol{t})$ does, in general. We distinguish the
following special case:
\begin{definition}\label{DefIsomon}
Let \eqref{FLDEt} be a holomorphic family of Fuchsian systems, where
the residue matrices $Q_j(\boldsymbol{t})$ are holomorphic in a
domain $\boldsymbol{\mathcal{D}}\subseteq\mathbb{C}^n_{\ast}$ and
satisfy \eqref{IRSP}. For each
$\boldsymbol{t}\in\boldsymbol{\mathcal{D}}$, let
 $Y({x},\boldsymbol{t})$ be the solution of \eqref{FLDEt} --
 \eqref{ICt}.

The family \eqref{FLDEt} is said to be an {\em isomonodromic family
of Fuchsian systems with the distinguished point} $x_0=\infty$ if
for every $\boldsymbol{t}^0\in\boldsymbol{\mathcal{D}}$  there
exists  a cylindrical open neighborhood
$\boldsymbol{\mathcal{W}}\Subset\boldsymbol{\mathcal{D}}$ of\,
$\boldsymbol{t}^0$, such  that
 the
following  holds.
\begin{itemize}
\item[]
For every loop
 $\gamma$ in the perforated sphere $\mathbb{CP}^1\setminus\bigcup_k
\overline{\mathcal{W}_k}$ with the distinguished point $x_0=\infty$
and every pair of points
$\boldsymbol{t^\prime},\boldsymbol{t^{\prime\prime}}\in\boldsymbol{\mathcal{W}}$
 the monodromy matrices $M_{\gamma}(\boldsymbol{t^\prime}), M_{\gamma}(\boldsymbol{t^{\prime\prime}})$ of the solutions
$Y( x,\boldsymbol{t^\prime})$, $Y(
x,\boldsymbol{t^{\prime\prime}})$, which correspond
 to this loop $\gamma$, are equal:
\begin{equation}
\label{isomon}
M_\gamma(\boldsymbol{t^\prime})=M_\gamma(\boldsymbol{t^{\prime\prime}})\end{equation}
\end{itemize}
\end{definition}

\begin{remark}
Note that if the family \eqref{FLDEt} is isomondromic  with the
distinguished point $x_0=\infty$ and
$\boldsymbol{t}^0\in\boldsymbol{\mathcal{D}}$, then  the monodromy
matrices of $Y$ are constant with respect to $\boldsymbol{t}$  in
every cylindrical open neighborhood $\boldsymbol{\mathcal{W}}$ of\,
$\boldsymbol{t}^0$, such  that
$\boldsymbol{\mathcal{W}}\Subset\boldsymbol{\mathcal{D}}$.
\end{remark}

\subsection{Isoprincipal families: the informal
definition}\label{infodef} Now we introduce the notion of the {\em
isoprincipal\,} family of Fuchsian systems, which is the central
notion in the present article.

Again, we assume that \eqref{FLDEt} is a holomorphic family of
Fuchsian systems, where the residue matrices $Q_j(\boldsymbol{t})$
are holomorphic in a domain
$\boldsymbol{\mathcal{D}}\subseteq\mathbb{C}^n_{\ast}$ and satisfy
\eqref{IRSP}, and  consider the solution
 $Y({x},\boldsymbol{t})$ of \eqref{FLDEt} -- \eqref{ICt}.

According to Proposition \ref{RePrDec}, for each fixed
$\boldsymbol{t}\in\boldsymbol{\mathcal{D}}$ a branch $Y_{\alpha_j}$
of the solution $Y$ in a neighborhood of $t_j$ admits the
regular-principal factorization, but now both the regular factor
$H_{j}$ and the principal factor $P_{\alpha_j}$ may depend on
$\boldsymbol{t}$:
\begin{equation}
\label{RPDt}%
 Y_{\alpha_j}( x,\boldsymbol{t})=H_{j}(x,\boldsymbol{t})\cdot
P_{\alpha_j}(x,\boldsymbol{t}).
\end{equation}
For example, if the principal factor is chosen in the form
\eqref{SPF}, then $Z_{\alpha_j}=Z_{\alpha_j}(\boldsymbol{t}),$ $
A_{\alpha_j}= A_{\alpha_j}(\boldsymbol{t})$ and
\begin{equation}
\label{SPFt}%
 P_{\alpha_j}( x,\boldsymbol{t})=
 (x-t_j)^{Z_{\alpha_j}(\boldsymbol{t})}\cdot
( x-t_j)^{ A_{\alpha_j}(\boldsymbol{t})}.
\end{equation}

Roughly speaking, the family \eqref{FLDEt} is {\em isoprincipal} if
for every $j$, $1\leq j\leq n$,  the matrices $Z_{\alpha_j}$ and $
A_{\alpha_j}$ in \eqref{SPFt} do not depend on $\boldsymbol{t}$:
$Z_{\alpha_j}(\boldsymbol{t})\equiv Z_{\alpha_j},
A_{\alpha_j}(\boldsymbol{t})\equiv  A_{\alpha_j}$, and
\begin{equation}
\label{SIPFt}%
 P_{\alpha_j}( x,\boldsymbol{t})=
 (x-t_j)^{Z_{\alpha_j}}\cdot ( x-t_j)^{ A_{\alpha_j}}.
\end{equation}
If the matrices $Z_{\alpha_j}$ and $A_{\alpha_j}$ in \eqref{SPFt} do
not depend on $\boldsymbol{t}$, then the principal factor
$P_{\alpha_j}( x,\boldsymbol{t})$ of the form \eqref{SIPFt}
possesses the following property: {\em it depends only on the
difference} $ x-t_j$.

For our goals, the specific form \eqref{SIPFt} of the principal
factors is of no importance. We just need each principal factor $
P_{\alpha_j}(x,\boldsymbol{t})$ to depend only on the difference $
x-t_j.$

 This means that there exist functions
$E_{\alpha_j},$    such that
\begin{equation}\label{prinfactrae}
P_{\alpha_j}( x,\boldsymbol{t})=E_{\alpha_j}( x-t_j),
\end{equation}
or, in the language of differential equations,
\begin{subequations}\label{firpla}
\begin{align}\dfrac{\partial P_{\alpha_j}}{\partial t_\ell}&=0,\quad
\ell\not=j,\\
\dfrac{\partial P_{\alpha_j}}{\partial x}&=-\dfrac{\partial
P_{\alpha_j}}{\partial t_j}=\left.\dfrac{d E_{\alpha_j}}{d
\zeta}\right|_{\zeta=x-t_j}.\end{align}
\end{subequations}

The formal definition of the  isoprincipal family of Fuchsian
systems (see Definition \ref{DefIsoprDefo} below) is more involved,
since for each $\boldsymbol{t}$ the branch $Y_{\alpha_j}(
x,\boldsymbol{t})$ and the principal factor $P_{\alpha_j}(
x,\boldsymbol{t})$ depend on the choice of a path $\alpha_j$ in the
punctured sphere $\mathbb{CP}^1\setminus\{t_1,\dots,t_n\}$, which
connects the distinguished point $x_0=\infty$ with a neighborhood of
$x=t_j$.

Moreover,  for each $\boldsymbol{t}$ the principal factor
$P_{\alpha_j}( x,\boldsymbol{t})$ should be a
transplant\,\footnote{See Definition \ref{DefTrans}.} of a function
$E_{\alpha_j}(\zeta)$, holomorphic on the Riemann surface of
$\ln\zeta$, into a universal covering surface over
$\mathbb{C}\setminus\{t_j\}$ and   these transplants should be
defined coherently with respect to $\boldsymbol{t}$.

\subsection{Isoprincipal families: the formal definition}
\begin{definition}\label{defcor}
Let a function $E(\zeta)$ be defined on the Riemann surface of
$\ln\zeta$. Let $\mathcal{W}_j$ be a simply connected domain,
compactly included in $\mathbb{C}$. Let
$p_j\in\mathbb{C}\setminus\overline{\mathcal{W}_j}$ and let us
choose a
 branch of $\arg(p_j-t_j)$, {\em continuous} with respect
to $t_j$ in $\mathcal{W}_j$ (since $\mathcal{W}_j$ is simply
connected and  $\mathcal{W}_j\Subset\mathbb{C}\setminus\{p_j\}$,
such a choice can be made).  For every $t_j\in\mathcal{W}_j$ let us
specify the value of $\arg(p_j-t_j)$ in this manner and consider the
corresponding transplant  $E(x-t_j)$ of the function $E(\zeta)$ {
into} $ \cov(\mathbb{C}\setminus\{t_j\};p_j).$

Then the family of   transplants
$\{E(x-t_j)\}_{t_j\in\mathcal{W}_j}$ is said to be {\em coherent
with respect to} ${t}_j$ {\em in} $\mathcal{W}_j$.
\end{definition}

\begin{definition}
Let $\boldsymbol{t}^0$ be a point in a domain
$\boldsymbol{\mathcal{D}}\subseteq\mathbb{C}^n_{\ast}$
 and let $\boldsymbol{\mathcal{W}},\boldsymbol{\mathcal{V}}$
be a  pair of open cylindrical neighborhoods of the point
$\boldsymbol{t}^0$, such that
$$\boldsymbol{\mathcal{V}}\Subset\mathbb{C}^n_{\ast}\quad\text{ and
}\quad \boldsymbol{\mathcal{W}}\Subset
\left\{\boldsymbol{\mathcal{V}}\cap\boldsymbol{\mathcal{D}}\right\}.$$
 Then the pair
$\boldsymbol{\mathcal{W}},\boldsymbol{\mathcal{V}}$ is said to be a
{\em nested pair of open cylindrical neighborhoods of the point}
$\boldsymbol{t}^0$ {in} $\boldsymbol{\mathcal{D}}.$
\end{definition}

\begin{definition}[Formal definition of the isoprincipal family]
\label{DefIsoprDefo}%
 Let \eqref{FLDEt} be a holomorphic family  of
Fuchsian systems, where the residue matrices $Q_j(\boldsymbol{t})$
are holomorphic in a domain
$\boldsymbol{\mathcal{D}}\subseteq\mathbb{C}^n_{\ast}$ and satisfy
\eqref{IRSP}. For  $\boldsymbol{t}\in\boldsymbol{\mathcal{D}}$ let
 $Y({x},\boldsymbol{t})$ be the solution of \eqref{FLDEt} --
 \eqref{ICt}.

The family \eqref{FLDEt} is said to be an {\em isoprincipal family
of Fuchsian systems with the distinguished point} $x_0=\infty$ if
for every $\boldsymbol{t}^0\in\boldsymbol{\mathcal{D}}$ there exists
a nested pair of open cylindrical neighborhoods of
$\boldsymbol{t}^0$:
 $$\boldsymbol{\mathcal{V}}=
\mathcal{V}_1\times \cdots \times\mathcal{V}_n,\quad
\boldsymbol{\mathcal{W}}=\mathcal{W}_1\times \cdots
\times\mathcal{W}_n,\quad
\boldsymbol{\mathcal{W}}\Subset\boldsymbol{\mathcal{V}},$$
 such that
 the following  holds.

\begin{itemize}
\item[] For every path $\alpha_j$ in the perforated sphere  $\mathbb{CP}^1\setminus
\bigcup_k\overline{\mathcal{W}_k}$  from the distinguished point
$x_0=\infty$ to a point
$p_j\in\mathcal{V}_j\setminus\overline{\mathcal{W}_j}$, $1\leq j\leq
n$, there exists a coherent family of transplants
$\{E_{\alpha_j}(x-t_j)\}_{t_j\in\mathcal{W}_j}$  of a function
 $E_{\alpha_j}(\zeta)$, holomorphic and invertible on the Riemann surface of $\ln\zeta$, such that for each
$\boldsymbol{t}\in\boldsymbol{\mathcal{W}}$
 the  branch $Y_{\alpha_j}({x},\boldsymbol{t})$ of the solution
$Y({x},\boldsymbol{t})$ in the punctured domain
$\mathcal{V}_j\setminus\{t_j\}$
  admits
   the
representation
\begin{equation}%
\label{RePrFactT} Y_{\alpha_j}({x},\boldsymbol{t})=
H_{j}(x,\boldsymbol{t})\cdot E_{\alpha_j}(x-t_j),\quad
{x}\in\cov(\mathcal{V}_j\setminus\{t_j\};\,p_j),
\end{equation}%
where $H_{j}(x,\boldsymbol{t})$ is a function, holomorphic (with
respect to $x$) and invertible in the entire domain $\mathcal{V}_j.$
\end{itemize}
\end{definition}

\begin{remark}  In view of Definition \ref{DefFact}, Definition \eqref{DefIsoprDefo}
means that the family \eqref{FLDEt} is isoprincipal with the
distinguished point $x_0=\infty$ if every branch of the solution
$Y(x,\boldsymbol{t})$ in a neighborhood of each singular point
$x=t_j$ admits the regular-principal factorization
\eqref{RePrFactT}, where the principal factor is the {\em
appropriately shifted copy}  of a function $E_{\alpha_j}(\zeta)$,
which is holomorphic and invertible on the Riemann surface of
$\ln\zeta$ and does not depend on $\boldsymbol{t}$. The meaning of
the words "appropriately shifted copy" is made precise in Definition
\ref{defcor}: this is what we call  a coherent family of transplants
of $E_{\alpha_j}(\zeta)$.

Thus Definition \ref{DefIsoprDefo} is a formal interpretation of the
{\em informal} definition in Section \ref{infodef}.

We would also like to note that it suffices to consider only the
branches corresponding to a certain choice of the paths
$\alpha_1,\dots,\alpha_n$ -- the same choice as the one mentioned in
Remark \ref{genmont} below.
\end{remark}
\subsection{Every isoprincipal family is an isomonodromic one.}
\begin{theorem}
\label{IsoprImplIsom}%
Let \eqref{FLDEt} be a holomorphic family of Fuchsian systems, where
the residue matrices $Q_j(\boldsymbol{t})$ are holomorphic in a
domain $\boldsymbol{\mathcal{D}}\subseteq\mathbb{C}^n_{\ast}$ and
satisfy \eqref{IRSP}.

Assume that the family  \eqref{FLDEt} is isoprincipal with the
distinguished point $x_0=\infty$.

Then this family  is isomonodromic with the distinguished point
$x_0=\infty$.
\end{theorem}%

Before we turn to the proof of Theorem \ref{IsoprImplIsom}, let us
introduce the following "$\boldsymbol{t}$-dependent" counterpart of
Definition \ref{loops}:

\begin{definition}\label{loopst}
Let  $\boldsymbol{t}^0$ be a point in the domain
$\boldsymbol{\mathcal{D}}$ and let
$\boldsymbol{\mathcal{W}}\Subset\boldsymbol{\mathcal{V}}$ be a
nested pair
 of open
cylindrical neighborhoods of $\boldsymbol{t}^0$  in
$\boldsymbol{\mathcal{D}}$. For $1\leq j\leq n$ let $\alpha_j$ be a
path in the perforated sphere
$\mathbb{CP}^1\setminus\bigcup_k\overline{\mathcal{W}_k}$ from the
distinguished point $x_0=\infty$ to a point
$p_j\in\mathcal{V}_j\setminus\overline{\mathcal{W}_j}.$

Furthermore, assume that $\beta_j$ is the loop  in the annulus
$\mathcal{V}_j\setminus\overline{\mathcal{W}_j}$ with the
distinguished point $p_j$ which makes one positive circuit of the
set ${\mathcal{W}_j}$,  and let $\gamma_j$ be the loop in the
perforated sphere
$\mathbb{CP}^1\setminus\bigcup_k\overline{\mathcal{W}_k}$ with the
distinguished point $x_0=\infty $, defined by
\begin{equation}
\label{bigloopt}
\gamma_j\defi\alpha_j^{-1}\cdot\beta_j\cdot\alpha_j.
\end{equation}

Then: \begin{itemize} \item the loop $\beta_j$ is said to be {\em
the small loop around the set} ${\mathcal{W}_j}$ {\em in the
annulus} $\mathcal{V}_j\setminus\overline{\mathcal{W}_j}$;
\item the loop $\gamma_j$ is said to be {\em the big loop around the
set} ${\mathcal{W}_j}$, {\em corresponding to the path} $\alpha_j$.
\end{itemize}
\end{definition}

\begin{remark}\label{genmont} Similarly to the case of a fixed
$\boldsymbol{t}$ (see Remark \ref{genmon}),  for a suitable choice
of the paths $\alpha_1,\dots,\alpha_n$ the corresponding big loops
$\gamma_1,\dots,\gamma_n$ generate the fundamental group
$\fug(\mathbb{CP}^1\setminus\bigcup_k\overline{\mathcal{W}_k};\,\infty)$.
\end{remark}

\begin{proof}[Proof of Theorem \ref{IsoprImplIsom}]
 Let $\boldsymbol{t}^0$ be a point in $\boldsymbol{\mathcal{D}}$ and let
$\boldsymbol{\mathcal{W}}\Subset\boldsymbol{\mathcal{V}}$  be a
nested pair
 of open
cylindrical neighborhoods of $\boldsymbol{t}^0$ as in Definition
\ref{DefIsoprDefo}.

In view of Remark \ref{genmont}, it suffices to prove that if
$\alpha_j$ is a path in
$\mathbb{CP}^1\setminus\bigcup_k\overline{\mathcal{W}_k}$ from the
distinguished point $x_0=\infty$ to a point
$p_j\in\mathcal{V}_j\setminus\overline{\mathcal{W}_j}$ and
$\gamma_j$ is the corresponding big loop around ${\mathcal{W}_j},$
then the monodromy matrix $M_{\gamma_j}(\boldsymbol{t})$ of
$Y(x,\boldsymbol{t})$ does not depend on $\boldsymbol{t}$:
\begin{equation}\label{mondot}
M_{\gamma_j}(\boldsymbol{t})=\text{const}, \quad
\boldsymbol{t}\in\boldsymbol{\mathcal{W}}.
\end{equation}

 In view of \eqref{mongen}, for each  fixed $\boldsymbol{t}\in
\boldsymbol{\mathcal{W}}$  the monodromy matrix
$M_{\gamma_j}(\boldsymbol{t})$  is given by
$$M_{\gamma_j}(\boldsymbol{t})= Y_{\alpha_j}^{-1}({x},\boldsymbol{t})\cdot Y_{\alpha_j}
(t_j+(x-t_j) e^{2\pi i},\boldsymbol{t}),\quad
\forall{x}\in\cov(\mathcal{V}_j\setminus\{t_j\};\,p_j),$$ where
$Y_{\alpha_j}({x},\boldsymbol{t})$ is the branch of
$Y({x},\boldsymbol{t})$ in $\mathcal{V}_j\setminus\{t_j\}$,
corresponding to the path $\alpha_j$.

 Substituting the expression
\eqref{RPDt} for $Y_{\alpha_j}$ into the above identity and taking
into
 account that the  factor $H_{j}(x,\boldsymbol{t})$ is a
single-valued function of $x$, we obtain
\begin{multline*}
M_{\gamma_j}(\boldsymbol{t})= E_{\alpha_j}^{-1}({x}-t_j)\cdot
E_{\alpha_j}(({x}-t_j) e^{2\pi
i})\\=\left.\big(E_{\alpha_j}^{-1}(\zeta)\cdot E_{\alpha_j}(\zeta
e^{2\pi i})\big)\right|_{\zeta={x}-t_j},\quad
\forall{x}\in\cov(\mathcal{V}_{j}\setminus\{t_j\};\,p_j).\end{multline*}

Thus the function $E_{\alpha_j}^{-1}(\zeta)\cdot E_{\alpha_j}(\zeta
e^{2\pi i})$, holomorphic on the Riemann surface of $\ln\zeta$, is
constant with respect to $\zeta$ on a certain non-empty open subset
of this surface. Therefore, this function is identically constant on
the Riemann surface of $\ln\zeta$ and we write
\begin{equation}\label{monprinfac}
M_{\gamma_j}(\boldsymbol{t})=E_{\alpha_j}^{-1}(\zeta)\cdot
E_{\alpha_j}(\zeta e^{2\pi i})\quad
\forall\boldsymbol{t}\in\boldsymbol{\mathcal{W}},\,\zeta\in\cov(\mathbb{C}\setminus\{0\};1).\end{equation}
But  the right-hand side of the last identity  does not depend  on
$\boldsymbol{t}$,  hence we obtain \eqref{mondot}.
\end{proof}

\subsection{Every non-resonant isomonodromic family is an
isoprincipal one.}
 The  converse to Theorem \ref{IsoprImplIsom} is
only conditionally true: it holds under the assumption that all the
matrices $Q_j$ are {\em non-resonant} (see footnote
${}^{\ref{DNRM}}$). In general, however,
 an
isomonodromic family can be non-isoprincipal. The appropriate
counterexample will be presented in Section \ref{Examp} of this
paper.

\begin{lemma}\label{IsomImplIsospec}
 Let \eqref{FLDEt} be a
holomorphic family of Fuchsian systems, where the residue matrices
$Q_j(\boldsymbol{t})$ are holomorphic in a domain
$\boldsymbol{\mathcal{D}}\subseteq\mathbb{C}^n_{\ast}$ and satisfy
\eqref{IRSP}. Assume that this family  is isomonodromic with the
distinguished point $x_0=\infty$.

Then the family \eqref{FLDEt} is isospectral in the following sense:
for every pair of points $\boldsymbol{t^\prime},
\boldsymbol{t^{\prime\prime}} \in\boldsymbol{\mathcal{D}}$ and each
$j$, $1\leq j\leq n$, the spectra
$\spec{Q_j(\boldsymbol{t^\prime})}$ and
$\spec{Q_j(\boldsymbol{t^{\prime\prime}})}$ are equal\,\footnote{As
usual, the spectra are considered "with multliplicities".}:
\begin{equation}\label{isospec}
\spec{Q_j(\boldsymbol{t^\prime})}=\spec{Q_j(\boldsymbol{t^{\prime\prime}})},\quad\forall\boldsymbol{t^\prime},
\boldsymbol{t^{\prime\prime}} \in\boldsymbol{\mathcal{D}},\ 1\leq
j\leq n.\end{equation}
\end{lemma}
\begin{theorem}
\label{IsomImplIsopr}%
Let \eqref{FLDEt} be a holomorphic family of Fuchsian systems, where
the residue matrices $Q_j(\boldsymbol{t})$ are holomorphic in a
domain $\boldsymbol{\mathcal{D}}\subseteq\mathbb{C}^n_{\ast}$ and
satisfy \eqref{IRSP}.

Assume that this family satisfies the following  conditions:
\begin{enumerate}
\item
It is isomonodromic with the distinguished point $x_0=\infty$;
\item
At least for one point $\boldsymbol{t}\in\boldsymbol{\mathcal{D}}$,
each of the matrices $Q_j(\boldsymbol{t}),j=1, \ldots , n,$ is
non-resonant.
\end{enumerate}

Then the family \eqref{FLDEt} is isoprincipal with the distinguished
point $x_0=\infty$.
\end{theorem}%

\begin{remark}\label{wlog}
It should be noted that, unlike the rest of our considerations in
the present article, the proofs of Lemma \ref{IsomImplIsospec} and
Theorem \ref{IsomImplIsopr} utilize the explicit form of the
principal factors mentioned in Remark \ref{PosMD}.
\end{remark}

\begin{proof}[Proof of Lemma \ref{IsomImplIsospec}]
 Let $\boldsymbol{t}^0$ be a point in $\boldsymbol{\mathcal{D}}$ and let
$\boldsymbol{\mathcal{W}}\Subset\boldsymbol{\mathcal{V}}$  be a
nested pair
 of open
cylindrical neighborhoods of $\boldsymbol{t}^0$ in
$\boldsymbol{\mathcal{D}}$ as in Definition \ref{DefIsomon}. For
$j=1,\dots, n$ let us choose a path $\alpha_j$ in the perforated
sphere $\mathbb{CP}^1\setminus\bigcup_k\overline{\mathcal{W}_k}$
from the distinguished point $x_0=\infty$ to a point
$p_j\in\mathcal{V}_j\setminus\overline{\mathcal{W}_j}$. As usual,
for each $\boldsymbol{t}\in\boldsymbol{\mathcal{W}}$ we denote by
$Y_{\alpha_j}(x,\boldsymbol{t})$ the branch of the solution
$Y({x},\boldsymbol{t})$ of \eqref{FLDEt} -- \eqref{ICt} in the
punctured domain $\mathcal{V}_j\setminus\{t_j\}$, corresponding to
this path $\alpha_j$.

Then, in view of Remark \ref{PosMD} (see \eqref{SPF},
\eqref{specla}), the branch $Y_{\alpha_j}({x},\boldsymbol{t})$
 admits
 the regular-principal factorization
 $$
 Y_{\alpha_j}({x},\boldsymbol{t})=H_j(x,\boldsymbol{t})\cdot(x-t_j)^{Z_{\alpha_j}(\boldsymbol{t})}\cdot
( x-t_j)^{ A_{\alpha_j}(\boldsymbol{t})},\quad
\boldsymbol{t}\in\boldsymbol{\mathcal{W}},\
{x}\in\cov(\mathcal{V}_j\setminus\{t_j\};\,p_j),$$
 where
$$ \spec\left(e^{2\pi i  A_{\alpha_j}(\boldsymbol{t})}\right)=\spec\left(e^{2\pi i
Q_j(\boldsymbol{t})}\right).$$ Let $\gamma_j$ be the big loop around
$\mathcal{W}_j$, corresponding to the path $\alpha_j$. Then the
 monodromy matrix $M_{\gamma_j} (\boldsymbol{t})$ of
$Y({x},\boldsymbol{t})$, corresponding to the loop $\gamma_j$ is
given by
\begin{equation}
\label{mongenj} M_{\gamma_j}(\boldsymbol{t})=
Y_{\alpha_j}^{-1}({x},\boldsymbol{t})\cdot Y_{\alpha_j}(t_j+(x-t_j)
e^{2\pi i},\boldsymbol{t})=e^{2\pi i
A_{\alpha_j}(\boldsymbol{t})}.\end{equation}
 Therefore,
$$\spec(M_{\gamma_j}(\boldsymbol{t}))=\spec\left(e^{2\pi i Q_j(\boldsymbol{t})}\right),
\quad\boldsymbol{t}\in\boldsymbol{\mathcal{W}}.$$ Since the family
\eqref{FLDEt} is isomonodromic, we have
$$ M_{\gamma_j}(\boldsymbol{t})=M_{\gamma_j}(\boldsymbol{t^0}),\quad\boldsymbol{t}\in\boldsymbol{\mathcal{W}},$$
hence
$$\spec\left(e^{2\pi i
Q_j(\boldsymbol{t})}\right)=\spec\left(e^{2\pi i
Q_j(\boldsymbol{t^0})}\right),\quad\boldsymbol{t}\in\boldsymbol{\mathcal{W}}.$$
This means that the spectra
 $\spec{Q_j(\boldsymbol{t})}$ and  $\spec{Q_j(\boldsymbol{t^0})}$ coincide modulo integers.
 But the function  $Q_j(\boldsymbol{t})$ is continuous with respect to $\boldsymbol{t}$ in $\boldsymbol{\mathcal{W}},$
 hence
 $$\spec{Q_j(\boldsymbol{t})}=\spec{Q_j(\boldsymbol{t^0})},
 \quad\boldsymbol{t}\in\boldsymbol{\mathcal{W}}.$$
 Since the above identity holds for all $\boldsymbol{t}$
  in a neighborhood of every point $\boldsymbol{t^0}\in\boldsymbol{\mathcal{D}}$,
we obtain the identity \eqref{isospec}.
\end{proof}

\begin{proof}[Proof of Theorem \ref{IsomImplIsopr}]
Let $\boldsymbol{t}^0$ be a point in $\boldsymbol{\mathcal{D}}$ and
let $\boldsymbol{\mathcal{W}}\Subset\boldsymbol{\mathcal{V}}$  be a
nested pair
 of open
cylindrical neighborhoods of $\boldsymbol{t}^0$ in
$\boldsymbol{\mathcal{D}}$ as in Definition \ref{DefIsomon}.

As in the proof of Lemma \ref{IsomImplIsospec}, for $j=1,\dots, n$
let us choose a path $\alpha_j$ in the perforated sphere
$\mathbb{CP}^1\setminus\bigcup_k\overline{\mathcal{W}_k}$ from the
distinguished point $x_0=\infty$ to a point
$p_j\in\mathcal{V}_j\setminus\overline{\mathcal{W}_j}$ and consider
for each fixed $\boldsymbol{t}\in\boldsymbol{\mathcal{W}}$ the
branch $Y_{\alpha_j}({x},\boldsymbol{t})$ of the solution
$Y({x},\boldsymbol{t})$ of \eqref{FLDEt} -- \eqref{ICt} in the
punctured domain $\mathcal{V}_j\setminus\{t_j\}$, corresponding to
this path $\alpha_j$.

Since by Lemma \ref{IsomImplIsospec}  the family \eqref{FLDEt} is
isospectral,  the matrix $Q_j(\boldsymbol{t})$ is non-resonant for
every $\boldsymbol{t}\in\boldsymbol{\mathcal{D}}$. Hence, in view of
Remark \ref{PosMD} (see the expressions \eqref{prinnonres},
\eqref{similar}),  the branch $Y_{\alpha_j}$
 admits
 the regular - principal factorization
$$
Y({x},\boldsymbol{t})=H_j(x,\boldsymbol{t})\cdot ( x-t_j)^{
A_{\alpha_j}(\boldsymbol{t})},\quad
{x}\in\cov(\mathcal{V}_j\setminus\{t_j\};\,p_j),$$ where the matrix
$A_{\alpha_j}(\boldsymbol{t})$ is similar to the matrix
$Q_j(\boldsymbol{t})$. Therefore, \begin{equation}\label{firid}
\spec( A_{\alpha_j}(\boldsymbol{t}))=\spec(
A_{\alpha_j}(\boldsymbol{t^0})),\quad
\boldsymbol{t}\in\boldsymbol{\mathcal{W}}.\end{equation}

In view of Definition \ref{DefIsoprDefo}, it remains to prove that
the matrix $A_{\alpha_j}(\boldsymbol{t})$ does not actually depend
on $\boldsymbol{t}$.

 Let $\gamma_j$ be the big loop around
$\mathcal{W}_j$, corresponding to the path $\alpha_j$. Then the
monodromy matrix $M_{\gamma_j} (\boldsymbol{t})$ of
$Y({x},\boldsymbol{t})$, corresponding to the loop $\gamma_j$, is
given by (see \eqref{mongenj})
$$
M_{\gamma_j}(\boldsymbol{t})= e^{2\pi i
A_{\alpha_j}(\boldsymbol{t})},\quad\boldsymbol{t}\in\boldsymbol{\mathcal{W}}.$$
Since the family \eqref{FLDEt} is isomonodromic, we have
$$ M_{\gamma_j}(\boldsymbol{t})=M_{\gamma_j}(\boldsymbol{t^0}),\quad\boldsymbol{t}\in\boldsymbol{\mathcal{W}},$$
hence $$e^{2\pi i A_{\alpha_j}(\boldsymbol{t})}= e^{2\pi i
A_{\alpha_j}(\boldsymbol{t^0})}, \quad
\boldsymbol{t}\in\boldsymbol{\mathcal{W}}.$$ But, in view of
\eqref{firid},  the last identity implies
$$ A_{\alpha_j}(\boldsymbol{t})=A_{\alpha_j}(\boldsymbol{t^0})=\text{const},
\quad \boldsymbol{t}\in\boldsymbol{\mathcal{W}}.$$
\end{proof}

\section{Isoprincipal Families of Fuchsian Systems and the Schlesinger System}
\label{SchlSyst} \subsection{The  Schlesinger system}

 A natural question arises: how to
express the property of a family of Fuchsian systems
\begin{equation}
\label{holdefo}%
\dfrac{dY}{dx}=\bigg(\sum\limits_{1\leq j\leq
n}\dfrac{Q_j(\boldsymbol{t})}{x-t_j}\bigg) Y,
\end{equation}
to be  isoprincipal   in
terms of the residues matrix functions $Q_j(\boldsymbol{t})$? Here,
as always, we assume that the residue matrices $Q_j(\boldsymbol{t})$
are  holomorphic in a domain
$\boldsymbol{\mathcal{D}}\subseteq\mathbb{C}^n_{\ast}$ and satisfy
\begin{equation}%
\label{regcont}%
\sum\limits_{1\leq j\leq n}Q_j(\boldsymbol{t})\equiv 0,\quad
\boldsymbol{t}\in\boldsymbol{\mathcal{D}}.
\end{equation}%

It turns out that an answer to this question is given by the
so-called {\em Schlesinger system} of PDEs:

\begin{equation}%
\label{Sch} \left\{
\begin{aligned}
\dfrac{\partial Q_{i}}{\partial t_{j}}&=
\dfrac{[Q_{i},Q_{j}]}{t_{i}-t_{j}}, &\quad
 1\leq
i,j\leq n,\ i\not=j, \\[1.0ex]
\dfrac{\partial Q_{i}}{\partial t_{i}}&= - \sum_{\substack{1\leq
j\leq n
\\
j\not=i}} \dfrac{\left[Q_{i},Q_{j}\right]} {t_{i}-t_{j}}, & \quad
1\leq i\leq n.
\end{aligned}\right.
\end{equation}

The following theorem is the main result of the present article.

\begin{theorem}[The main result]\label{main}
Let \eqref{holdefo} be a holomorphic family of Fuchsian systems,
where the residue matrices $Q_j(\boldsymbol{t})$ are holomorphic in
a domain $\boldsymbol{\mathcal{D}}\subseteq\mathbb{C}^n_{\ast}$ and
satisfy \eqref{regcont}.

Then the family \eqref{holdefo} is isoprincipal with the
distinguished point $x_0=\infty$ if and only if the residue matrices
$Q_j(\boldsymbol{t})$ satisfy with respect to $\boldsymbol{t}$ the
Schlesinger system \eqref{Sch} in the domain
$\boldsymbol{\mathcal{D}}$.
\end{theorem}
\begin{remark}\label{firin}
Note that \eqref{Sch} implies
$$\dfrac{\partial}{\partial t_j}\sum_{1\leq i\leq n}Q_i=0,\quad 1\leq
j\leq n,$$ that is, $\sum_iQ_i(\boldsymbol{t})$ is a first integral
of the Schlesinger system.

In particular, if  functions $Q_j(\boldsymbol{t})$,  $1\leq j\leq
n,$ satisfy the Schlesinger system \eqref{Sch} in the domain
$\boldsymbol{\mathcal{D}}$ and at some point
$\boldsymbol{t^0}\in\boldsymbol{\mathcal{D}}$ it holds that
$$\sum_{1\leq j\leq n}Q_j(\boldsymbol{t^0})=0,$$
 then these functions $Q_j(\boldsymbol{t})$ satisfy the relation
 \eqref{regcont}.
\end{remark}

\begin{remark} We would like to stress that in Theorem \ref{main} no assumptions on the
spectra of the residue matrices $Q_j(\boldsymbol{t})$ are made. Thus Theorem \ref{main}
for the {\em isoprincipal} families of Fuchsian systems
 can be viewed as an amended version of L. Schlesinger's statement, concerning
the {\em isomonodromic} deformations (see \cite{Sch2} and the introduction of the present article).

In the case of the isomonodromic families of Fuchsian systems Theorem \ref{main}
implies the following:
\begin{enumerate}
\item
 {\em If
the residue matrices $Q_j(\boldsymbol{t})$ satisfy  the Schlesinger
system \eqref{Sch}, then the
family \eqref{holdefo} is} isoprincipal and hence by Theorem
\ref{IsoprImplIsom} also {\em isomonodromic.}

\item {\em If the family \eqref{holdefo} is isomonodromic and,}
 in addition, {\em all the residue matrices $Q_j(\boldsymbol{t})$
are non-resonant (at least at some point), then} by Theorem
\ref{IsomImplIsopr} the family \eqref{holdefo} is isoprincipal and
hence {\em the residue matrices
$Q_j(\boldsymbol{t})$ satisfy  the Schlesinger system \eqref{Sch}.}
\end{enumerate}

We remark that in the statement (ii) the assumption of non-resonance for the residues
$Q_j(\boldsymbol{t})$ cannot be omitted: in Section \ref{Examp} we shall present
an example of the isomonodromic family \eqref{holdefo}, where the
residue matrices $Q_j(\boldsymbol{t})$ are resonant and do {\em not}
satisfy the Schlesinger system \eqref{Sch} (thus contradicting the statement of L. Schlesinger).

 Nevertheless, our proof of the "only if" part of the Theorem \ref{main}
 largely follows the original proof of L. Schlesinger for the
 isomonodromic case
 (see also   \cite[Section 3.5]{IKSY}, where the modern adaptation of
  Schlesinger's proof is presented). In
 particular, the overdetermined linear system \eqref{auxisopr},
 which appears in Proposition \ref{prlinisopr} below and is
 crucial  in the derivation  of the Schlesinger system, can be found
 in \cite[Section II]{Sch2}.
\end{remark}

The proof of Theorem \ref{main} will be split into parts and
presented as a series of propositions, culminating with Propositions
\ref{IsoToSchles} and \ref{SchlesToIso}.

\subsection{The auxiliary system related to the isomonodromic
family of Fuchsian systems} In order to prove Theorem \ref{main}, we
have to study the partial derivatives of the solution
$Y({x},\boldsymbol{t})$, satisfying the initial condition
\begin{equation}
\label{incont} Y({x},\boldsymbol{t})\bigm|_{x=\infty}=I,
\end{equation}%
with respect to the parameters $t_1,\dots,t_n$.

 First of all, let us choose and fix a
point $\boldsymbol{t^0}\in\boldsymbol{\mathcal{D}}$ and let
$\boldsymbol{\mathcal{W}}\Subset\boldsymbol{\mathcal{D}}$ be a
cylindrical open neighborhood of $\boldsymbol{t^0}$.
 {\em Since the
coefficients of the system \eqref{holdefo} and the initial condition
\eqref{incont} depend holomorphically on $\boldsymbol{t}$, the
solution $Y(x,\boldsymbol{t})$ is holomorphic jointly in $x$ and
$\boldsymbol{t}$ in the Cartesian product}
$\cov(\mathbb{CP}^1\setminus\bigcup_k\overline{\mathcal{W}_k};\,\infty)\times\boldsymbol{\mathcal{W}}.$

In particular, the partial derivatives
$\dfrac{\partial{Y}}{\partial{t_j}}(x, \boldsymbol{t^0})$
 are defined
and holomorphic with respect to  $x$ in
$\cov(\mathbb{CP}^{1}\setminus\bigcup_k\overline{\mathcal{W}_k};\,\infty).$
Since  these definitions of $\dfrac{\partial{Y}}{\partial{t_j}}(x,
\boldsymbol{t^0})$ agree for various choices of
$\boldsymbol{\mathcal{W}}$ as long as $\boldsymbol{\mathcal{W}}$ is
sufficiently small, we conclude that {\em for each fixed
$\boldsymbol{t^0}\in\boldsymbol{\mathcal{D}}$ the partial
derivatives $\dfrac{\partial Y}{\partial t_j}(x,\boldsymbol{t^0})$
are defined and holomorphic as functions of $x$  on the the same
surface as the function $Y(x,\boldsymbol{t^0})$ itself -- the
universal covering surface
$\cov(\mathbb{CP}^1\setminus\{t_1^0\dots,t_n^0\};\,\infty)$.}

It turns out that in terms of these partial derivatives of $Y$ the
property of the family \eqref{holdefo} to be isomonodromic can be
expressed as follows:

\begin{proposition}\label{prlinisom}
Let \eqref{holdefo} be a holomorphic family of Fuchsian systems,
where the residue matrices $Q_j(\boldsymbol{t})$ are holomorphic in
a domain $\boldsymbol{\mathcal{D}}\subseteq\mathbb{C}^n_{\ast}$ and
satisfy \eqref{regcont}.

Then the family \eqref{holdefo} is isomonodromic  with the
distinguished point $x_0=\infty$ if and only if  the solution
$Y(x,\boldsymbol{t})$ of \eqref{holdefo}, \eqref{incont}
 satisfies  a linear system of
the form
\begin{equation}\label{auxisom}\left\{\begin{aligned}\dfrac{\partial Y}{\partial x}&=
\sum\limits_{1\leq j\leq n}\dfrac{Q_j(\boldsymbol{t})}{x-t_j}\cdot
Y,\\
 \dfrac{\partial Y}{\partial t_j}&=T_j(x,t)\cdot Y,\quad 1\leq j\leq
n,\end{aligned}\right.\end{equation} where for each
$\boldsymbol{t}\in\boldsymbol{\mathcal{D}}$ the functions
$T_j(x,t)$, $1\leq j\leq n,$ are  \textsf{single-valued} holomorphic with
respect to $x$ in $\mathbb{CP}^1 \setminus\{t_1,\dots,t_n\}$.
\end{proposition}

\begin{definition}
Let the family \eqref{holdefo}  of Fuchsian systems be isomonodromic
with the distinguished point $x_0=\infty$.

 The system
\eqref{auxisom} with the single-valued coefficients $T_j(x,t)$,
$1\leq j\leq n,$ which appears in Theorem \ref{prlinisom}, is said
to be {\em the  auxiliary linear system  related to the
isomonodromic family \eqref{holdefo} of Fuchsian systems.}
\end{definition}

\begin{proof}[Proof of Proposition \ref{prlinisom}]
The first equation of the system \eqref{auxisom} is just the
Fuchsian system \eqref{holdefo} itself, hence we only need to prove
that for each fixed $\boldsymbol{t}\in\boldsymbol{\mathcal{D}}$ the
logarithmic derivatives
$$T_j( x,\boldsymbol{t})\defi
\dfrac{\partial Y}{\partial t_j}({x},\boldsymbol{t})\cdot
Y^{-1}({x},\boldsymbol{t}),$$
 which {\em a priori}  are defined  as holomorphic functions
of $x$ on the universal covering surface
$\cov(\mathbb{CP}^1\setminus\{t_1,\dots,t_n\};\,\infty),$ are {\em
single-valued} in the punctured sphere
$\mathbb{CP}^1\setminus\{t_1,\dots,t_n\}.$

 Let us  choose a point
$\boldsymbol{t}^0\in\boldsymbol{\mathcal{D}}$ and a cylindrical open
neighborhood
$\boldsymbol{\mathcal{W}}\Subset\boldsymbol{\mathcal{D}}$  of\,
$\boldsymbol{t}^0$. Let
$\gamma\in\fug(\mathbb{CP}^1\setminus\bigcup_k\overline{\mathcal{W}_k};\,\infty)$
 and  let us denote by
$$x\mapsto x\cdot\gamma,\quad x\in\cov(\mathbb{CP}^1\setminus
\bigcup_k\overline{\mathcal{W}_k};\,\infty)$$ the  deck
transformation of the universal covering surface
$\cov(\mathbb{CP}^1\setminus\bigcup_k\overline{\mathcal{W}_k};\,\infty)
$, corresponding\,\footnote{See \eqref{deck}.} to this loop
$\gamma$.

According to Definition \ref{DefMonMatr}, the monodromy matrix
$M_\gamma(\boldsymbol{t})$ of the solution $Y({x},\boldsymbol{t})$,
which corresponds to the loop $\gamma$, is given by
$$M_{\gamma}(\boldsymbol{t})=Y^{-1}({x},\boldsymbol{t})\cdot Y({x}\gamma,\boldsymbol{t}),\quad
x\in\cov(\mathbb{CP}^1\setminus\bigcup_{1\leq k\leq
n}\overline{\mathcal{W}_k};\,\infty),\
\boldsymbol{t}\in\boldsymbol{\mathcal{W}}.
$$ The monodromy matrix $M_\gamma(\boldsymbol{t})$ does not depend on ${x}$  and hence is  a holomorphic
single-valued\,\footnote{Recall (see Definition \ref{defcyl}) that
all the bases $\mathcal{W}_k$ of the cylindrical neighborhood
$\boldsymbol{\mathcal{W}}$ are simply connected.} function of
$\boldsymbol{t}$ in $\boldsymbol{\mathcal{W}}$. Differentiating the
equality
$$ Y({x}\gamma,\boldsymbol{t})=Y({x},\boldsymbol{t})\cdot M_{\gamma}(\boldsymbol{t})$$
with respect to $t_j$, we obtain
$$\dfrac{\partial Y}{\partial t_j}({x}\gamma,\boldsymbol{t})=
\dfrac{\partial Y}{\partial t_j}({x},\boldsymbol{t})\cdot
M_{\gamma}(\boldsymbol{t})+
Y({x},\boldsymbol{t})\cdot\dfrac{\partial M_{\gamma}}{\partial
t_j}(\boldsymbol{t}).$$ Therefore, the logarithmic derivative
$$T_j( x,\boldsymbol{t})=
\dfrac{\partial Y}{\partial t_j}({x},\boldsymbol{t})\cdot
Y^{-1}({x},\boldsymbol{t})$$  satisfies the monodromy
 relation
 $$T_j( x\gamma,\boldsymbol{t})=T_j( x,\boldsymbol{t})+
 Y({x},\boldsymbol{t})\cdot\dfrac{\partial M_{\gamma}}{\partial t_j}(\boldsymbol{t})\cdot  Y^{-1}
 ({x}\gamma,\boldsymbol{t}).$$

 The last equality implies that the following two statements are
 equivalent\,\footnote{This equivalence is stronger than the statement of Proposition \ref{prlinisopr} in the
 sense that
 it holds for each individual loop $\gamma$ and each individual index $j$.}:
 \renewcommand{\theenumi}{\textup{\arabic{enumi}}}
 \begin{enumerate}
 \item
The monodromy matrix  $M_{\gamma}(\boldsymbol{t})$ does not depend
on $t_j$:
 $$
\dfrac{\partial M_{\gamma}}{\partial t_j}(\boldsymbol{t})\equiv
0,\quad\boldsymbol{t}\in\boldsymbol{\mathcal{W}}.$$
\item
It holds  that
 $$
 T_j( x\gamma,\boldsymbol{t})\equiv T_j( x,\boldsymbol{t}),\quad
x\in\cov(\mathbb{CP}^1\setminus\bigcup_{1\leq k\leq
n}\overline{\mathcal{W}_k};\,\infty),\
\boldsymbol{t}\in\boldsymbol{\mathcal{W}}.$$
\end{enumerate}
However, the statement 2) holds for every
$\gamma\in\fug(\mathbb{CP}^1\setminus\bigcup_k\overline{\mathcal{W}_k};\,\infty)$
if and only if for each $\boldsymbol{t}\in\boldsymbol{\mathcal{W}}$
the function $T_j(x,\boldsymbol{t})$ is  a {\em single-valued}
function of $x$ in $\mathbb{CP}^1 \setminus\{t_1,\dots,t_n\}$. In
view of Definition \ref{DefIsomon}, this completes the proof.
\end{proof}
\subsection{The auxiliary system related to the isoprincipal family
of Fuchsian systems} Now we turn to the case, when the family
\eqref{holdefo} is not only isomonodromic but, moreover,
isoprincipal. We claim that in this special case the  auxiliary
linear system \eqref{auxisom} can be written explicitly in terms of
the residues $Q_j(\boldsymbol{t})$.

\begin{proposition}\label{prlinisopr}
Let \eqref{holdefo} be a holomorphic family of Fuchsian systems,
where the residue matrices $Q_j(\boldsymbol{t})$ are holomorphic in
a domain $\boldsymbol{\mathcal{D}}\subseteq\mathbb{C}^n_{\ast}$ and
satisfy \eqref{regcont}. Assume that the family \eqref{holdefo} is
isoprincipal  with the distinguished point $x_0=\infty$.

Then the solution $Y( x,\boldsymbol{t})$ of  \eqref{holdefo},
\eqref{incont}  satisfies the following auxiliary system:
\begin{equation}\label{auxisopr}\left\{\begin{aligned}\dfrac{\partial Y}{\partial x}&=
\sum\limits_{1\leq j\leq n}\dfrac{Q_j(\boldsymbol{t})}{x-t_j}\cdot
Y,\\
\dfrac{\partial Y}{\partial
t_j}&=-\dfrac{Q_j(\boldsymbol{t})}{x-t_j}\cdot Y,\quad 1\leq j\leq
n.\end{aligned}\right.\end{equation}
\end{proposition}

In the proof of Proposition \ref{prlinisopr} we shall use the
following
\begin{lemma}\label{harlem}
Let $\mathcal{U},\mathcal{V}$ be simply-connected domains in the
complex plane $\mathbb{C}$, such that:
\begin{enumerate}
\item $\mathcal{U}\Subset\mathcal{V}$;
\item the set $\mathcal{V}\setminus\overline{\mathcal{U}}$
is connected.
\end{enumerate}

Let $\boldsymbol{\mathcal{W}}$ be a domain in $\mathbb{C}^n$ and let
$H(x,\boldsymbol{t})$ be a function of $x\in\mathcal{V}$ and
$\boldsymbol{t}\in\boldsymbol{\mathcal{W}}$, possessing the
following properties:
\renewcommand{\theenumi}{\alph{enumi}}
\begin{enumerate}
\item
the function $H(x,\boldsymbol{t})$ is holomorphic (jointly in $x$
and $\boldsymbol{t}$) in
$\{\mathcal{V}\setminus\overline{\mathcal{U}}\}\times\boldsymbol{\mathcal{W}}$;
\item for each fixed $\boldsymbol{t}\in\boldsymbol{\mathcal{W}}$
the function $H(x,\boldsymbol{t})$ is holomorphic with respect to
$x$ in the entire domain $\mathcal{V}$.
\end{enumerate}

Then the function $H(x,\boldsymbol{t})$ is holomorphic (jointly in
$x$ and $\boldsymbol{t}$) in the domain
$\mathcal{V}\times\boldsymbol{\mathcal{W}}$.
\end{lemma}
\begin{remark}
Lemma  \ref{harlem} is a   special case of the well-known Hartogs
lemma. However, in this simple case the conclusion  follows
immediately from the Cauchy integral formula.

Indeed, let $\beta$ be a smooth loop  in the annulus
$\mathcal{V}\setminus\overline{\mathcal{W}}$ which makes one
positive circuit of the set $\mathcal{W}$  and let
$\Delta\Subset\boldsymbol{\mathcal{W}}$ be a polydisk:
$$\Delta=\Delta_1\times\cdots\Delta_n.$$ Then  for each  $\boldsymbol{t}\in\Delta$ and
$x\in\overline{\mathcal{W}}$ it holds that
\begin{multline*}H_j(x,\boldsymbol{t})=\dfrac{1}{2\pi
i}\oint\limits_{\beta_j}\dfrac{H_j(\zeta,\boldsymbol{t})}{\zeta-x}d\zeta\\[1ex]
=\dfrac{1}{(2\pi
i)^{n+1}}\oint\limits_{\beta_j}\left(\oint\limits_{\partial\Delta_n}\cdots\oint\limits_{\partial\Delta_1}
\dfrac{H_j(\zeta,\tau_1,\dots\tau_n)}{(\tau_1-t_1)\cdots(\tau_n-t_n)}d\tau_1\cdots
d\tau_n\right) \dfrac{d\zeta}{\zeta-x},
\end{multline*}
where $\partial\Delta_k$ denotes the boundary of the disk
$\Delta_k$. Since the contours of integration lie in the domain
$\{\mathcal{V}\setminus\overline{\mathcal{W}}\}\times
\boldsymbol{\mathcal{W}}$, where $H_j(x,\boldsymbol{t})$ is jointly
holomorphic  in $x$ and $\boldsymbol{t}$ (in particular,
continuous),  the integral represents a function, jointly
holomorphic in $x$ and $\boldsymbol{t}$ for $x$ in a neighborhood of
$\overline{\mathcal{W}_j}$ and $\boldsymbol{t}\in\Delta$.
\end{remark}

\begin{proof}[Proof of Proposition \ref{prlinisopr}]
 Let $\boldsymbol{t^0}\in\boldsymbol{\mathcal{D}}$ and let $\boldsymbol{\mathcal{W}}\Subset\boldsymbol{\mathcal{V}}$
 be a nested pair of cylindrical open neighborhoods of\, $\boldsymbol{t^0}$.

  For $\ell=1,\dots,n$  let us consider the logarithmic derivative
 \begin{equation}
T_\ell( x,\boldsymbol{t})= \dfrac{\partial Y}{\partial
t_\ell}({x},\boldsymbol{t})\cdot
Y^{-1}({x},\boldsymbol{t}).\end{equation}
 Since by Theorem \ref{IsoprImplIsom}
the isoprincipal family \eqref{holdefo} is also isomonodromic,
Proposition \ref{prlinisom} implies  that for each fixed
$\boldsymbol{t}\in\boldsymbol{\mathcal{W}}$ the function
$T_\ell(x,\boldsymbol{t})$ is single-valued holomorphic  with
respect to $x$ in the punctured sphere $\mathbb{CP}^1 \setminus
\{t_1,\dots,t_n\}$. We  have to prove  that
\begin{equation}\label{desire}
T_\ell(x,\boldsymbol{t})=-\dfrac{Q_\ell(\boldsymbol{t})}{x-t_\ell}.
\end{equation}

To this end let us introduce the auxiliary function
\begin{equation}\label{lau}
F_\ell(x,\boldsymbol{t})\defi
T_\ell(x,\boldsymbol{t})+\dfrac{Q_\ell(\boldsymbol{t})}{x-t_\ell},\end{equation}
which is single-valued holomorphic with respect to $x$ in the
punctured sphere $\mathbb{CP}^1\setminus\{t_1,\dots,t_n\}.$ We are
going to show that for $j=1,\dots, n$ the following statement holds:
\begin{itemize}
\item[$\boldsymbol{(\star)}$] The point $x=t_j$ is
a
 removable singularity of the function $F_\ell(x,\boldsymbol{t})$.
 \end{itemize}

 Then, according to Liouville theorem, the function
 $F_\ell(x,\boldsymbol{t})$ is constant
 with respect to $x$. Since, as  follows from
 \eqref{incont},\eqref{lau}
\begin{equation}\label{logderinf}
T_\ell(
x,\boldsymbol{t})\bigm|_{x=\infty}=F_\ell(x,\boldsymbol{t})\bigm|_{x=\infty}=
0,
\end{equation}
we can then conclude that $$F_\ell(x,\boldsymbol{t})\equiv 0,$$
which leads to the desired result \eqref{desire}.

Now we turn to the proof of the statement ${(\star)}$.
 To begin with,
 let us
  choose and fix a path
  $\alpha_j$ in the perforated sphere
  $\mathbb{CP}^{1}\setminus\bigcup_k\overline{\mathcal{W}_k}$ from the distinguished point $x_0=\infty$ to
  a point $p_j\in\mathcal{V}_j\setminus\overline{\mathcal{W}_j}.$
    Then,
according to Definition \ref{DefIsoprDefo}, for each
$\boldsymbol{t}\in\boldsymbol{\mathcal{W}}$ the branch
$Y_{\alpha_j}(x,\boldsymbol{t})$ of the solution
$Y({x},\boldsymbol{t})$ in $\mathcal{V}_j\setminus\{t_j\}$ admits
the regular-principal factorization
\begin{equation}\label{regprinfact}
Y_{\alpha_j}({x},\boldsymbol{t})= H_{j}(x,\boldsymbol{t})\cdot
E_{\alpha_j}({x}-t_j), \quad
x\in\cov(\mathcal{V}_j\setminus\{t_j\};\,p_j),
\end{equation}
where the  family $\{E_{\alpha_j}({x}-t_j)\}_{t_j\in\mathcal{W}_j}$
 is a coherent family of transplants of a function
$E_{\alpha_j}(\zeta)$, holomorphic and invertible on the Riemann
surface of $\ln\zeta$, and the function $H_{j}(x,\boldsymbol{t})$ is
holomorphic with respect to $x$ and invertible in the entire
(non-punctured) domain $\mathcal{V}_j$.

Then the principal factor $E_{\alpha_j}({x}-t_j)$ is jointly
holomorphic in $x$ and $t_j$ in
$\cov(\mathbb{C}\setminus\overline{\mathcal{W}_j};\,p_j)\}\times\mathcal{W}_j,$
and hence the regular factor
$$H_{j}(x,\boldsymbol{t})=Y_{\alpha_j}^{-1}({x},\boldsymbol{t})\cdot E_{\alpha_j}({x}-t_j)$$
 is jointly
holomorphic (single-valued) in $x$ and $\boldsymbol{t}$ in
$\{\mathcal{V}_j\setminus\overline{\mathcal{W}_j}\}\times
\boldsymbol{\mathcal{W}}$. Since the function
$H_{j}(x,\boldsymbol{t})$ is also holomorphic with respect to $x$ in
the entire  domain $\mathcal{V}_j$,  Lemma \ref{harlem} implies that
it is jointly holomorphic  in $x$ and $\boldsymbol{t}$ in
$\mathcal{V}_j\times \boldsymbol{\mathcal{W}}$.

Thus
  we can differentiate the equality \ref{regprinfact} with
respect to $t_\ell$, $1\leq\ell\leq n$. First we consider the case
 $\ell\not=j$. Then, since
\begin{equation}\label{tribuu1}\dfrac{\partial E_{\alpha_j}({x}-t_j)}{\partial t_\ell}=0,\quad
\ell\not=j,\end{equation}
  we obtain for $x\in\mathcal{V}_j\setminus\{t_j\}$
$$T_\ell(x,\boldsymbol{t})=\dfrac{\partial Y_{\alpha_j}}{\partial t_\ell}({x},\boldsymbol{t})=
\dfrac{\partial H_{j}}{\partial t_\ell}(x,\boldsymbol{t})\cdot
E_{\alpha_j}({x}-t_j)\\=\dfrac{\partial H_{j}}{\partial
t_\ell}(x,\boldsymbol{t})\cdot H_j^{-1}(x,\boldsymbol{t})\cdot
Y_{\alpha_j}({x},\boldsymbol{t}).$$ Hence
$$F_\ell( x,\boldsymbol{t})= \dfrac{\partial
H_j}{\partial t_\ell}({x},\boldsymbol{t})\cdot
H_j^{-1}({x},\boldsymbol{t})+\dfrac{Q_\ell(\boldsymbol{t})}{x-t_\ell},\quad
x\in\mathcal{V}_j\setminus\{t_j\},\ \ell\not=j,$$ which proves the statement
$(\star)$ in the case
 $\ell\not=j$.

Next we differentiate the equality \eqref{regprinfact} with respect
to $t_j$. Since \begin{equation}\label{tribuu2}\dfrac{\partial
E_{\alpha_j}({x}-t_j)}{\partial x}=-\dfrac{\partial
E_{\alpha_j}({x}-t_j)}{\partial t_j}=\left.\dfrac{d
E_{\alpha_j}(\zeta)}{d \zeta}\right|_{\zeta=x-t_j},\end{equation} we
obtain
\begin{multline*}T_j( x,\boldsymbol{t})= \dfrac{\partial H_j}{\partial t_j}({x},\boldsymbol{t})\cdot
H_j^{-1}({x},\boldsymbol{t})-\\-H_j({x},\boldsymbol{t})\cdot
\dfrac{\partial E_{\alpha_j}({x}-t_j)}{\partial x}\cdot
E_{\alpha_j}^{-1}({x}-t_j)\cdot H_j^{-1}({x},\boldsymbol{t}),\quad
x\in\cov(\mathcal{V}_j\setminus\{t_j\};\,p_j).
\end{multline*}
On the other hand, differentiating the equality \ref{regprinfact}
with respect to $x$, we get
\begin{multline*}
\dfrac{\partial Y_{\alpha_j}}{\partial x}({x},\boldsymbol{t})\cdot
Y_{\alpha_j}^{-1}({x},\boldsymbol{t})= \dfrac{\partial H_j}{\partial
x}({x},\boldsymbol{t})\cdot
H_j^{-1}({x},\boldsymbol{t})+\\+H_j({x},\boldsymbol{t})\cdot
\dfrac{\partial E_{\alpha_j}({x}-t_j)}{\partial x}\cdot
E_{\alpha_j}^{-1}({x}-t_j)\cdot H_j^{-1}({x},\boldsymbol{t}),\quad
x\in\cov(\mathcal{V}_j\setminus\{t_j\};\,p_j),
\end{multline*}
hence
 \begin{multline*}F_j( x,\boldsymbol{t})=\left(\dfrac{\partial H_j}{\partial
x}({x},\boldsymbol{t})+\dfrac{\partial H_j}{\partial
t_j}({x},\boldsymbol{t})\right)H_j^{-1}({x},\boldsymbol{t})+\\+\dfrac{Q_j(\boldsymbol{t})}{x-t_j}
 -\dfrac{\partial Y_{\alpha_j}}{\partial
x}({x},\boldsymbol{t})\cdot
Y_{\alpha_j}^{-1}({x},\boldsymbol{t}),\quad x\in\cov(\mathcal{V}_j\setminus\{t_j\};\,p_j).
\end{multline*}
Taking into account that
  $Y_{\alpha_j}({x},\boldsymbol{t})$ is a branch of the solution
 $Y(x,\boldsymbol{t})$ of the Fuchsian system \eqref{holdefo},
 we obtain
$$
F_j( x,\boldsymbol{t})=\left(\dfrac{\partial H_j}{\partial
x}({x},\boldsymbol{t})+\dfrac{\partial H_j}{\partial
t_j}({x},\boldsymbol{t})\right)H_j^{-1}({x},\boldsymbol{t})-\sum_{\substack{1\leq
k\leq n\\ k\not =j}}\dfrac{Q_k(\boldsymbol{t})}{x-t_k},\quad
x\in\mathcal{V}_j\setminus\{t_j\}.$$
  Thus the function
$F_j( x,\boldsymbol{t})$ has at $x=t_j$ a removable singularity,
which proves the statement $(\star)$ in the case  $j=\ell$.
\end{proof}
\begin{remark} Note that the relations \eqref{tribuu1}, \eqref{tribuu2},
which are instrumental in the proof of Proposition \ref{prlinisopr},
are precisely the relations  \eqref{firpla} in the informal
{definition}  of the isoprincipal family of Fuchsian systems in
Section \ref{infodef}.
\end{remark}
\begin{remark}\label{regsysor}
We observe that the auxiliary linear system \eqref{auxisopr} leads
to a linear system for the regular factor $H_j(x,\boldsymbol{t})$ of
the regular-principal factorization \eqref{regprinfact}.

Indeed, in view of \eqref{tribuu1} and \eqref{tribuu2}, we can
differentiate the equality \eqref{regprinfact} to obtain
\begin{align*}
\dfrac{\partial H_j}{\partial t_\ell}({x},\boldsymbol{t})\cdot
H_j^{-1}({x},\boldsymbol{t})&=\dfrac{\partial Y_{\alpha_j}}{\partial
t_\ell}({x},\boldsymbol{t})\cdot
Y_{\alpha_j}^{-1}({x},\boldsymbol{t}),\quad j\not =\ell,\\
\left(\dfrac{\partial H_j}{\partial
t_j}({x},\boldsymbol{t})+\dfrac{\partial H_j}{\partial
x}({x},\boldsymbol{t})\right)\cdot
H_j^{-1}({x},\boldsymbol{t})&=\left(\dfrac{\partial
Y_{\alpha_j}}{\partial t_j}({x},\boldsymbol{t})+\dfrac{\partial
Y_{\alpha_j}}{\partial x}({x},\boldsymbol{t})\right)\cdot
Y_{\alpha_j}^{-1}({x},\boldsymbol{t}),
\end{align*}
and therefore
\begin{equation}\label{befcha}
\left\{\begin{aligned}\dfrac{\partial H_j}{\partial
t_\ell}&=-\dfrac{Q_\ell}{x-t_\ell}\cdot H_j,\quad \ell\not=j,\\
\dfrac{\partial H_j}{\partial t_j}+\dfrac{\partial H_j}{\partial
x}&=\sum\limits_{\substack{1\leq \ell\leq n\\ \ell\not =
j}}\dfrac{Q_\ell}{x-t_\ell}\cdot
H_j.\end{aligned}\right.\end{equation} Using the change of variables
\begin{align}
\zeta&=x-t_j,\\ \label{chavar} L_j(\zeta,\boldsymbol{t})&\defi
H_j(\zeta+t_j,\boldsymbol{t}),\end{align}
 one can rewrite the system \eqref{befcha} in the
following form:
\begin{equation}\label{afcha}
\left\{\begin{aligned}\dfrac{\partial L_j}{\partial
t_\ell}&=-\dfrac{Q_\ell}{\zeta+t_j-t_\ell}\cdot L_j,\quad \ell\not=j,\\
\dfrac{\partial L_j}{\partial t_j}&=\sum\limits_{\substack{1\leq
\ell\leq n\\ \ell\not = j}}\dfrac{Q_\ell}{\zeta+t_j-t_\ell}\cdot
L_j.\end{aligned}\right.\end{equation} The system \eqref{afcha} is
nothing more than the system \eqref{auxisopr} with the constraint
$$x-t_j=\zeta=\text{const}.$$
Note that, although  $x=t_j$ is a singularity of the Fuchsian system
\eqref{holdefo} and the auxiliary system \eqref{auxisopr},  the
right-hand side of   the system \eqref{afcha} is holomorphic with
respect to $\zeta$  at $\zeta=0$ (compare with Remark {1.2} in \cite{Ma1}).  This occurs because the function
$L_j(\zeta,\boldsymbol{t})$, defined in \eqref{chavar},  is
holomorphic with respect to $\zeta$ and invertible at $\zeta=0.$
\end{remark}
\subsection{The Frobenius theorem}

The auxiliary system \eqref{auxisopr}, related to the isoprincipal
family \eqref{holdefo} is an {\em overdetermined} system of PDEs.
The criterion for the existence of solution of such a system is
known (see \cite[Section 2.11]{Na}) as {\em the Frobenius theorem:}

\begin{theorem}[The Frobenius theorem for Pfaffian\,\footnote{A
Pfaffian system is a system of the form \eqref{oversys}.} systems]\label{genfrob}
 Let $\boldsymbol{\Omega_p}$ and
$\boldsymbol{\Omega_q}$ be  domains in, respectively, $\mathbb{C}^p$
and $\mathbb{C}^q$.  Consider the following system of PDEs:
\begin{equation}
 \label{oversys}
 \dfrac{\partial {\boldsymbol{\lambda}}}{\partial
 \mu_j}=\boldsymbol{\phi_{j}}(\boldsymbol{\lambda},\boldsymbol{\mu}),\quad
  1\leq j\leq q,
 \end{equation}
 where $\boldsymbol{\lambda}(\boldsymbol{\mu})$ is an unknown
$\mathbb{C}^p$-valued function of the variable
$\boldsymbol{\mu}=(\mu_1,\dots,\mu_q)\in\mathbb{C}^q$ and
$\boldsymbol{\phi_{j}}(\boldsymbol{\lambda},\boldsymbol{\mu})=
(\phi_{1,j}(\boldsymbol{\lambda},\boldsymbol{\mu}),\dots,
\phi_{p,j}(\boldsymbol{\lambda},\boldsymbol{\mu}))$,  $1\leq j\leq
q,$ are given $\mathbb{C}^p$-valued functions, holomorphic with
respect to $\boldsymbol{\lambda},\boldsymbol{\mu}$ in the domain
$\boldsymbol{\Omega_p}\times\boldsymbol{\Omega_q}.$

 Then the
following statements (i) and (ii) are equivalent:
\begin{enumerate}\item For every pair of points $\boldsymbol{\lambda^0}\in\boldsymbol{\Omega_p},$
$\boldsymbol{\mu^0}\in\boldsymbol{\Omega_q}$ there exists  a
solution $\boldsymbol{\lambda}(\boldsymbol{\mu})$ of the system
\eqref{oversys},  holomorphic in a neighborhood
 of  $\boldsymbol{\mu^0}$  and
 satisfying the initial condition\,\footnote{Note that, according to the
 uniqueness theorem for {\em ordinary} differential equations, such
 a solution $\boldsymbol{\lambda}(\boldsymbol{\mu})$ is necessarily unique.}
 $$\boldsymbol{\lambda}(\boldsymbol{\mu^0})=\boldsymbol{\lambda^0}.$$
\item The $\mathbb{C}$-valued functions
$\phi_{i,j}(\boldsymbol{\lambda},\boldsymbol{\mu})$ satisfy
 the equations
 \begin{multline}\label{compacon}
\dfrac{\partial
\phi_{i,j}(\boldsymbol{\lambda},\boldsymbol{\mu})}{\partial\mu_k}+\sum_{1\leq\ell\leq
p}\dfrac{\partial
\phi_{i,j}(\boldsymbol{\lambda},\boldsymbol{\mu})}{\partial\lambda_\ell}\cdot
\phi_{\ell,k}(\boldsymbol{\lambda},\boldsymbol{\mu})=\\
=\dfrac{\partial
\phi_{i,k}(\boldsymbol{\lambda},\boldsymbol{\mu})}{\partial\mu_j}+\sum_{1\leq\ell\leq
p}\dfrac{\partial
\phi_{i,k}(\boldsymbol{\lambda},\boldsymbol{\mu})}{\partial\lambda_\ell}\cdot
\phi_{\ell,j}(\boldsymbol{\lambda},\boldsymbol{\mu}),\\
1\leq i\leq p,\ 1\leq j,k\leq q,
 \end{multline}
 in the domain $\boldsymbol{\Omega_p}\times\boldsymbol{\Omega_q}$.
 \end{enumerate}
\end{theorem}
\begin{definition}
 The condititon
\eqref{compacon}, formulated in the Frobenius theorem \ref{genfrob},
 is said to be {\em the compatibility condition} for the overdetermined system of PDEs \eqref{oversys}.
The overdetermined system of PDEs  \eqref{oversys} which satisfies
the compatibility condition \eqref{compacon} is said to be {\em
compatible.}
\end{definition}

\begin{remark}\label{simpfrob}
Note that the compatibility condition \eqref{compacon} can be obtained by
substituting the equations \eqref{oversys} into the identity

$$
\dfrac{\partial^2 {\lambda_i}}{\partial
 \mu_j\partial\mu_k}=\dfrac{\partial^2 {\lambda_i}}{\partial
 \mu_k\partial\mu_j},\quad 1\leq i\leq
p,\ 1\leq j,k\leq q.$$

Thus the statement (ii) of the Frobenius theorem \eqref{genfrob}
follows from the statement (i)  immediately.
\end{remark}

In what follows, we shall often deal with  {\em linear}
overdetermined systems of PDEs, {\em depending on a parameter.}
Because of the {\em global} existence theorem for such systems,
 the following stronger version of the
Frobenius theorem \ref{genfrob} holds  in this case:

\begin{theorem}[The Frobenius theorem for linear
systems with a parameter]\label{linfrob}
 Let $\boldsymbol{\Omega_q}$ and
$\boldsymbol{\Omega_r}$ be  domains in, respectively, $\mathbb{C}^q$
and $\mathbb{C}^r$. Consider the following linear systems of PDEs:
\begin{equation}\label{linfrobsys}\dfrac{\partial \boldsymbol{\lambda}}{\partial
 \mu_j}=
 \Phi_{j}(\boldsymbol{\mu},\boldsymbol{\tau})\cdot \boldsymbol{\lambda},\quad
  1\leq j\leq q,
 \end{equation}
where $\boldsymbol{\mu}\in\boldsymbol{\Omega_q}$ is the "main"
variable, $\boldsymbol{\tau}\in\boldsymbol{\Omega_r}$ is a parameter
and $\Phi_{j}(\boldsymbol{\mu},\boldsymbol{\tau})$, $1\leq j\leq q,$
are  $\mathfrak{M}_p$-valued functions, holomorphic with respect to
$\boldsymbol{\mu}$ in the domain $\boldsymbol{\Omega_q}$ for each
fixed $\boldsymbol{\tau}\in\boldsymbol{\Omega_r}$.

Then:
\begin{enumerate}
\item
For each fixed value of the parameter $\boldsymbol{\tau}$, say
$\boldsymbol{\tau}=\boldsymbol{\tau^0}\in\boldsymbol{\Omega_r}$, the
linear system
\begin{equation}\label{linfrobsys0}\dfrac{\partial \boldsymbol{\lambda}}{\partial
 \mu_j}=
 \Phi_{j}(\boldsymbol{\mu},\boldsymbol{\tau^0})\cdot \boldsymbol{\lambda},\quad
  1\leq j\leq q,
 \end{equation}
 has a fundamental solution\,\footnote{In other words, an
 $\mathfrak{M}_p$-valued solution $\boldsymbol{\lambda}(\boldsymbol{\mu},\boldsymbol{\tau^0})$,
 such that $\det\boldsymbol{\lambda}(\boldsymbol{\mu},\boldsymbol{\tau^0})\not=0$ for
 $\boldsymbol{\mu}\in\boldsymbol{\Omega_q}$.}
$\boldsymbol{\lambda}(\boldsymbol{\mu},\boldsymbol{\tau^0})$
 if and only if the functions
$\Phi_{j}(\boldsymbol{\mu},\boldsymbol{\tau^0})$ satisfy the
equations
 \begin{multline}\label{lincompacon0}
\dfrac{\partial
\Phi_{j}(\boldsymbol{\mu},\boldsymbol{\tau^0})}{\partial\mu_k}-
\dfrac{\partial
\Phi_{k}(\boldsymbol{\mu},\boldsymbol{\tau^0})}{\partial\mu_j}=
[\Phi_k(\boldsymbol{\mu},\boldsymbol{\tau^0}),\Phi_j(\boldsymbol{\mu},\boldsymbol{\tau^0})],\\
\boldsymbol{\mu}\in\boldsymbol{\Omega_q},\quad 1\leq j,k\leq q.
 \end{multline}
\item If the functions
$\Phi_{j}(\boldsymbol{\mu},\boldsymbol{\tau})$ are jointly
holomorphic with respect to $\boldsymbol{\mu},\boldsymbol{\tau}$ in
the domain $\boldsymbol{\Omega_q}\times\boldsymbol{\Omega_r}$ and
satisfy the equations
 \begin{multline}\label{lincompacon}
\dfrac{\partial
\Phi_{j}(\boldsymbol{\mu},\boldsymbol{\tau})}{\partial\mu_k}-
\dfrac{\partial
\Phi_{k}(\boldsymbol{\mu},\boldsymbol{\tau})}{\partial\mu_j}=
[\Phi_k(\boldsymbol{\mu},\boldsymbol{\tau}),\Phi_j(\boldsymbol{\mu},\boldsymbol{\tau})],\\
\boldsymbol{\mu}\in\boldsymbol{\Omega_q},\
\boldsymbol{\tau}\in\boldsymbol{\Omega_r},\quad 1\leq j,k\leq q,
 \end{multline}
 then for every point $\boldsymbol{\mu^0}\in\boldsymbol{\Omega_q}$
 and every $\mathfrak{M}_p$-valued function $\boldsymbol{\lambda^0}(\boldsymbol{\tau}),$
 holomorphic and invertible in $\boldsymbol{\Omega_r},$ there exists
 a unique fundamental solution $\boldsymbol{\lambda}(\boldsymbol{\mu},\boldsymbol{\tau})$
 of the system \eqref{linfrobsys},  jointly
holomorphic with respect to $\boldsymbol{\mu},\boldsymbol{\tau}$ in
the domain $\boldsymbol{\Omega_q}\times\boldsymbol{\Omega_r}$ and
satisfying the initial condition
\begin{equation}\label{inlinfrob}
\boldsymbol{\lambda}(\boldsymbol{\mu},\boldsymbol{\tau})\bigm|_{\boldsymbol{\mu}=\boldsymbol{\mu^0}}\,=
\boldsymbol{\lambda^0}(\boldsymbol{\tau}).\end{equation}
\end{enumerate}
\end{theorem}

\subsection{Proof of Theorem \ref{main}}
\begin{proposition}
\label{IsoToSchles} Let \eqref{holdefo} be a holomorphic family of
Fuchsian systems, where the residue matrices $Q_j(\boldsymbol{t})$
are holomorphic in a domain
$\boldsymbol{\mathcal{D}}\subseteq\mathbb{C}^n_{\ast}$ and satisfy
\eqref{regcont}. Assume that the family \eqref{holdefo} is
isoprincipal  with the distinguished point $x_0=\infty$. Then the
residues $Q_1(\boldsymbol{t}),\dots,Q_n(\boldsymbol{t})$ satisfy the
Schlesinger  system \eqref{Sch} in the domain
$\boldsymbol{\mathcal{D}}$.
\end{proposition}

\begin{proof}
According to Proposition \ref{prlinisopr},
 the solution $Y(x,\boldsymbol{t})$
of  \eqref{holdefo}, \eqref{incont}  satisfies the auxiliary linear
system \eqref{auxisopr}. Therefore, $Y(x,\boldsymbol{t})$ is the
fundamental solution of the overdetermined linear system
\eqref{auxisopr} with the initial condition $$Y(\infty,
\boldsymbol{t^0})=I,$$ where $\boldsymbol{t^0}$ is an arbitrary
fixed point in the domain $\boldsymbol{\mathcal{D}}$.

Hence, in view of Theorem \ref{linfrob}, the linear system
\eqref{auxisopr} (which is a special case of the system
\eqref{linfrobsys} with $q=n+1$, $r=0$,
$\boldsymbol{\lambda}(\boldsymbol{\mu})=Y(x,\boldsymbol{t})$) is
compatible. The compatibility condition \eqref{lincompacon} takes in
this case the form of the following two equations:
\begin{align}\label{fircomcon}
\dfrac{\partial}{\partial{x}}\left(\dfrac{Q_j}{x-t_j}\right)+\sum\limits_{1\leq
i\leq
n}\dfrac{\partial}{\partial{t_j}}\left(\dfrac{Q_i}{x-t_i}\right)&=
\sum\limits_{1\leq i\leq n}\dfrac{[Q_i,Q_j]}{(x-t_i)(x-t_j)},\quad
1\leq j\leq n,\\
\label{seccomcon}\dfrac{\partial}{\partial{t_i}}\left(\dfrac{Q_j}{x-t_j}\right)
-\dfrac{\partial}{\partial{t_j}}\left(\dfrac{Q_i}{x-t_i}\right)&=\dfrac{[Q_j,Q_i]}{(x-t_i)(x-t_j)},
\quad 1\leq i,j\leq n.
\end{align}

Since the residues $Q_j$ do not depend on $x$, from the equation
\eqref{fircomcon} we obtain
$$
\sum\limits_{1\leq i\leq n}\dfrac{\partial
Q_i}{\partial{t_j}}\cdot\dfrac{1}{x-t_i}= \sum\limits_{1\leq i\leq
n}\dfrac{[Q_i,Q_j]}{(x-t_i)(x-t_j)},\quad 1\leq j\leq n.$$
 Both sides
of the last equality are  rational functions of  $x$, which are
holomorphic  in the punctured sphere
$\mathbb{CP}^{1}\setminus\{t_1,\dots, t_n\}$ and have simple poles
at the points  $x=t_i$, $1\leq i\leq n$.  Equating the residues at
each pole $x=t_i$, $1\leq i\leq n$, we obtain the equations
\eqref{Sch}.

The equation \eqref{seccomcon} leads in a similar way to the first
of the equations \eqref{Sch} and, therefore, provides no additional
information.
\end{proof}

\begin{proposition}\label{SchlesToIso}
Let \eqref{holdefo} be a holomorphic family of Fuchsian systems,
where the residue matrices $Q_j(\boldsymbol{t})$ are holomorphic in
a domain $\boldsymbol{\mathcal{D}}\subseteq\mathbb{C}^n_{\ast}$ and
satisfy \eqref{regcont}. Assume that the residue matrices
$Q_j(\boldsymbol{t})$ satisfy  the Schlesinger system \eqref{Sch}.

Then the family \eqref{holdefo} is isoprincipal with the
distinguished point $x_0=\infty$.
\end{proposition}

\begin{proof}
Let us choose  $\boldsymbol{t^0}\in\boldsymbol{\mathcal{D}}$ and a
nested pair of open cylindrical neighborhoods
$\boldsymbol{\mathcal{W}}\Subset\boldsymbol{\mathcal{V}}$ of
$\boldsymbol{t^0}$ in $\boldsymbol{\mathcal{D}}.$ Let
$Y(x,\boldsymbol{t})$ be the solution  of \eqref{holdefo},
\eqref{incont} and for $1\leq j\leq n$  let $\alpha_j$ be a path
from $x_0=\infty$ to a point $p_j\in\mathcal{V}_j\setminus
\overline{\mathcal{W}_j}$.

Then, according to Theorem \ref{RePrDec} the branch
$Y_{\alpha_j}(x,\boldsymbol{t^0})$ of the solution
$Y(x,\boldsymbol{t^0})$ in $\mathcal{V}_j\setminus\{t^0\}$,
corresponding to the path $\alpha_j$, admits  the regular-principal
factorization
\begin{equation}\label{rpf0}
Y_{\alpha_j}(x,\boldsymbol{t^0})=H_j(x,\boldsymbol{t^0})\cdot
P_{\alpha_j}(x,\boldsymbol{t^0}),\quad
x\in\cov(\mathcal{V}_j\setminus\{\boldsymbol{t^0}\};\,p_j),\end{equation}
 where
the factors $H_j(x,\boldsymbol{t^0})$  and
$P_{\alpha_j}(x,\boldsymbol{t^0})$ are holomorphic with respect to
$x$ and invertible in, respectively, $\mathcal{V}_j$ and
$\cov(\mathbb{C}\setminus\{t_j^0\};\,p_j).$

>From here we proceed in four steps:

\paragraph{\bf Step 1}
 Firstly,
we construct a function $E_{\alpha_j}(\zeta)$, holomorphic and
invertible on the Riemann surace of $\ln\zeta$ and such that the
principal factor $P_{\alpha_j}(x,\boldsymbol{t^0})$ is a transplant
of $E_{\alpha_j}(\zeta)$ into
$\cov(\mathbb{C}\setminus\{t_j^0\};\,p_j)$:
$$P_{\alpha_j}(x,\boldsymbol{t^0})=E_{\alpha_j}(x-t_j^0).$$

To this end we choose some value $\theta_j^0$ of $\arg(p_j-t_j^0)$
and consider the corresponding  isomorphism from the Riemann surace
of $\ln\zeta$ onto $\cov(\mathbb{C}\setminus\{t_j^0\};p_j)$ (see
\eqref{defisotrainv}):
\begin{equation}\label{diti}
\zeta\xrightarrow[\arg(p_j-t_j^0)=\theta_j^0]{} \zeta+t_j^0,\quad
\zeta\in\cov(\mathbb{C}\setminus\{0\};1),\
\zeta+t_j^0\in\cov(\mathbb{C}\setminus\{t_j^0\};p_j).\end{equation}
 We
set \begin{equation}\label{defprinl} E_{\alpha_j}(\zeta)\defi
P_{\alpha_j}(\zeta+t_j^0,\boldsymbol{t^0}),\end{equation} where
$\zeta+t_j^0$ denotes the image in
$\cov(\mathbb{C}\setminus\{t_j^0\};p_j)$ of the point
$\zeta\in\cov(\mathbb{C}\setminus\{0\};1)$ under the isomorphism
\eqref{diti}.
\paragraph{\bf Step 2}
Secondly, we construct the regular factor $H_j(x,\boldsymbol{t})$
for the solution $Y(x,\boldsymbol{t})$ of \eqref{holdefo},
\eqref{incont}. In view of Remark \ref{regsysor}, we consider the
overdetermined linear system
\begin{equation}%
\label{regsys} \left\{
\begin{aligned}
\dfrac{\partial L_j(\zeta,\boldsymbol{t})}{\partial t_{i}}&=
-\dfrac{Q_i(\boldsymbol{t})}{\zeta+t_{j}-t_{i}}\cdot L_j(\zeta,\boldsymbol{t}), &\quad
 1\leq
i\leq n,\ i\not=j, \\[1.0ex]
\dfrac{\partial L_j(\zeta,\boldsymbol{t})}{\partial t_{j}}&= \sum\limits_{\substack{1\leq
i\leq n
\\
i\not=j}} \dfrac{Q_i(\boldsymbol{t})} {\zeta+t_{j}-t_{i}}\cdot L_j(\zeta,\boldsymbol{t}),
\end{aligned}\right.
\end{equation}
with the initial condition
\begin{equation}\label{inregl}
L_j(\zeta,\boldsymbol{t})\bigm|_{\boldsymbol{t}=\boldsymbol{t^0}}=H_j(\zeta+t_j^0,\boldsymbol{t^0}),
\end{equation}
where $H_j(x,\boldsymbol{t^0})$ is the regular factor in the
factorization \eqref{rpf0}.
 Here $\zeta\in\mathbb{C}$ is a small
parameter:
\begin{multline}\label{defeps}|\zeta|<\epsilon,\text{ where }
 \epsilon
>0 \text{ is such that }\\
x_1\in\mathcal{W}_j\implies \{x\,:\,|x-x_1|<\epsilon\}\subset\mathcal{V}_j,\quad
1\leq j\leq n.\end{multline}

We claim that the overdetermined linear system \eqref{regsys},
depending on the parameter $\zeta$, is compatible (see Calculation 1
below).

Therefore, the system \eqref{regsys} has a solution
$L_j(\zeta,\boldsymbol{t})$, satisfying the initial condition
\eqref{inregl}. This solution $L_j(\zeta,\boldsymbol{t})$ is jointly
holomorphic in $\zeta,\boldsymbol{t}$  for
$|\zeta|<\epsilon,\boldsymbol{t}\in\boldsymbol{\mathcal W}$  and
invertible.

We define the function $H_j(x,\boldsymbol{t})$ by
\begin{equation}\label{defhj} H_j(x,\boldsymbol{t})\defi
L_j(x-t_j,\boldsymbol{t}),\quad
\boldsymbol{t}\in\boldsymbol{\mathcal{W}},\ x\in\mathcal{V}_j
\end{equation}

\paragraph{\bf Step 3}
Thirdly, we consider the product
\begin{equation}\label{defzzz}Z_{\alpha_j}(\zeta,\boldsymbol{t})\defi
L_j(\zeta,\boldsymbol{t})\cdot E_{\alpha_j}(\zeta).\end{equation}
Note  that, in view of \eqref{inregl},
 for $\boldsymbol{t}=\boldsymbol{t^0}$ we have
 $$Z_{\alpha_j}(\zeta,\boldsymbol{t^0})=Y_{\alpha_j}(\zeta+t_j^0,\boldsymbol{t^0}),$$
 hence the function $Z_{\alpha_j}(\zeta,\boldsymbol{t^0})$ satisfies
 the linear system
\begin{equation}\label{inz}
\dfrac{d Z_{\alpha_j}(\zeta,\boldsymbol{t^0})} {d
\zeta} = \sum\limits_{1\leq i\leq n}
\dfrac{Q_i(\boldsymbol{t^0})} {\zeta+t_j^0-t_{i}^0}\cdot
Z_{\alpha_j}(\zeta,\boldsymbol{t^0}).\end{equation} Also, as
follows from \eqref{regsys}, we have
\begin{align}\label{derzj1}
\dfrac{\partial Z_{\alpha_j}(\zeta,\boldsymbol{t})}{\partial
t_i}&=-\dfrac{Q_i(\boldsymbol{t})}{\zeta+t_j-t_i}\cdot
Z_{\alpha_j}(\zeta,\boldsymbol{t}),\quad 1\leq i\leq n,\ i\not =j,\\
\label{derzj} \dfrac{\partial Z_{\alpha_j}(\zeta,\boldsymbol{t})}{\partial t_{j}}&=
\sum\limits_{\substack{1\leq i\leq n
\\
i\not=j}} \dfrac{Q_i(\boldsymbol{t})} {\zeta+t_j-t_{i}}\cdot Z_{\alpha_j}(\zeta,\boldsymbol{t}).
\end{align}

Furthermore,  the function $Z_{\alpha_j}(\zeta,\boldsymbol{t})$ can
be  shown (see Calculation 2  below) to satisfy with respect to
$\zeta$ the linear system
\begin{equation}\label{derzz}\dfrac{\partial Z_{\alpha_j}(\zeta,\boldsymbol{t})} {\partial \zeta}=
\sum\limits_{1\leq i\leq n} \dfrac{Q_i(\boldsymbol{t})}
{\zeta+t_j-t_{i}}\cdot Z_{\alpha_j}(\zeta,\boldsymbol{t}).\end{equation}

We note that for each fixed
$\boldsymbol{t}\in\boldsymbol{\mathcal{W}}$ the linear differential
system \eqref{derzz} with respect to $\zeta$ has no singularities
 in the punctured domain
 $$\mathcal{V}_{t_j}\defi\{\zeta:\zeta+t_j\in\mathcal{V}_j\}\setminus\{0\},$$
 hence  its fundamental solution
$Z_{\alpha_j}(\zeta,\boldsymbol{t})$ is holomorphic with respect to
$\zeta$ on a universal covering surface over the domain
$\mathcal{V}_{t_j}$. Therefore, the function
$$L_j(\zeta,\boldsymbol{t})=Z_{\alpha_j}(\zeta,\boldsymbol{t})\cdot E_{\alpha_j}^{-1}(\zeta)$$ is
 holomorphic with respect to $\zeta$ and invertible on this universal covering surface.
On the other hand, the function $L_j(\zeta,\boldsymbol{t})$ is also
{\em single-valued} holomorphic
  with respect to $\zeta$ and invertible in the open disk $\{\zeta\,:\,|\zeta|<\epsilon\}$ (see \eqref{defeps}).
  Hence the function $L_j(\zeta,\boldsymbol{t})$ is single-valued holomorphic
  with respect to $\zeta$ and invertible in the {\em non-punctured}
  domain $\mathcal{V}_{t_j}\cup\{0\}.$

  It follows that the function $H_j(x,\boldsymbol{t})$, defined in
\eqref{defhj}, is holomorphic (single-valued) with
 respect to $x$ and invertible in the entire domain $\mathcal{V}_j$.

\paragraph{\bf Step 4}
 Finally, we
consider the coherent family of transplants
$\{E_{\alpha_j}(x-t_j)\}_{t_j\in\mathcal{W}_j}$ (see Definition
\ref{defcor}), corresponding to the unique branch of $\arg(p_j-t_j)$
which is continuous with respect to $t_j$ in $\mathcal{W}_j$ and
takes the value $\theta^0$, chosen at Step 1, at the point $t_j=t_j^0$.

For each $\boldsymbol{t}\in\boldsymbol{\mathcal{W}}$ we {\em define}
the function  $Y_{\alpha_j}(x,\boldsymbol{t})$ by
\begin{equation}\label{defyyy}  Y_{\alpha_j}(x,\boldsymbol{t})\defi
H_j(x,\boldsymbol{t})\cdot
E_{\alpha_j}(x-t_j)=Z_{\alpha_j}(x-t_j,\boldsymbol{t}),\quad
x\in\cov(\mathcal{V}_j\setminus\{t_j\};\,p_j).
\end{equation}
Note that in view of \eqref{defprinl}, \eqref{inregl} this definition agrees
for $\boldsymbol{t}=\boldsymbol{t^0}$ with \eqref{rpf0} and the notation
 $Y_{\alpha_j}(x,\boldsymbol{t^0})$  for the branch of
 the solution $Y(x,\boldsymbol{t^0})$ in $\mathcal{V}_j\setminus\{t_j^0\}$,
 corresponding to the path $\alpha_j$. Now we show that
 {\em for every $\boldsymbol{t}\in\boldsymbol{\mathcal{W}}$
the function  $Y_{\alpha_j}(x,\boldsymbol{t})$, defined in \eqref{defyyy}, is the branch of
 the solution $Y(x,\boldsymbol{t})$ of \eqref{holdefo}, \eqref{incont}
 in the punctured domain $\mathcal{V}_j\setminus\{t_j\}$,
 corresponding to the path $\alpha_j$.}

First of all, in view of \eqref{derzz}, the function
$ Y_{\alpha_j}(x,\boldsymbol{t})$ satisfies with respect to $x$ the
system \eqref{holdefo}:
$$\dfrac{\partial  Y_{\alpha_j}}{\partial x}(x,\boldsymbol{t})=\sum_{1\leq i\leq
n}\dfrac{Q_i(\boldsymbol{t})}{x-t_i}\cdot  Y_{\alpha_j}(x,\boldsymbol{t}).$$
 Hence for each
$\boldsymbol{t}\in\boldsymbol{\mathcal{W}}$ the function
$Y_{\alpha_j}(x,\boldsymbol{t})$ can be analytically continued along the path $\alpha_j$
in the opposite direction:
from $p_j$ to $\infty$. The value of this continuation at $x=\infty$ will be denoted
by $\hat Y (\infty,\boldsymbol{t});$ in particular it holds that
\begin{equation}\label{ict00}
\hat{Y}(\infty,\boldsymbol{t^0})={Y}(\infty,\boldsymbol{t^0})=I.\end{equation}

Furthermore, the value $\hat Y (\infty,\boldsymbol{t})$
can be considered as the {\em initial} value at the distinguished point $x_0=\infty$ for
a fundamental solution $\hat{Y}(x,\boldsymbol{t})$ of the Fuchsian system \eqref{holdefo}, defined
on the universal covering surface
$\cov(\mathbb{CP}^1\setminus\{t_1,\dots,t_n\};\,\infty)$. The function
$ Y_{\alpha_j}(x,\boldsymbol{t})$ is the branch of this solution $\hat{Y}(x,\boldsymbol{t})$
in the punctured domain
$\mathcal{V}_j\setminus\{t_j\}$, corresponding to the path
$\alpha_j$.

Now we note that, in view of \eqref{derzj1} -- \eqref{derzj},
$$
\dfrac{\partial \hat{Y}}{\partial
t_i}(x,\boldsymbol{t})=-\dfrac{Q_i(\boldsymbol{t})}{x-t_i}\cdot \hat{Y}(x,\boldsymbol{t}),\quad 1\leq
i\leq n;
$$
in particular,
\begin{equation}
\left.\dfrac{\partial \hat{Y}}{\partial
t_i}(x,\boldsymbol{t})\right|_{x=\infty}=0,\quad 1\leq i\leq n.\label{dec00}\end{equation}
Combining \eqref{dec00} with \eqref{ict00}, we observe that the solution
$\hat{Y}(x,\boldsymbol{t})$ satisfies the initial condition
$$\hat{Y}(x,\boldsymbol{t})\bigm|_{x=\infty}=I,\quad\boldsymbol{t}\in\boldsymbol{\mathcal{W}},$$
and by the uniqueness theorem for linear differential systems
 coincides with the fundamental solution
$Y(x,\boldsymbol{t})$ of \eqref{holdefo}, \eqref{incont}.

Thus the function  $Y_{\alpha_j}(x,\boldsymbol{t})$, defined in \eqref{defyyy}, is the branch of
 the solution $Y(x,\boldsymbol{t})$ of \eqref{holdefo}, \eqref{incont}
 in the punctured domain $\mathcal{V}_j\setminus\{t_j\}$,
 corresponding to the path $\alpha_j$. The equality \eqref{defyyy} itself
 can now be considered as
the regular-principal factorization of the branch $Y_{\alpha_j}(x,\boldsymbol{t})$.
In view of Definition
\ref{DefIsoprDefo}, we conclude that {\em the family \eqref{holdefo} is
isoprincipal with the distinguished point $x_0=\infty$.}

In order to complete the proof, it remains to  present the calculations, omitted in the above
reasonings.

\paragraph{\bf Calculation 1}
We show that the overdetermined linear system \eqref{regsys},
depending on the parameter $\zeta$, is compatible.

In this case the compatibility condition \eqref{lincompacon} of
Theorem \eqref{linfrob} is represented  by the following two
equalities:
\begin{subequations}\label{compfol}
\begin{multline}\label{compfol1}
\dfrac{\partial}{\partial{t_j}}\left(\dfrac{Q_i}{\zeta+t_{j}-t_{i}}\right)+\sum\limits_{\substack{1\leq
k\leq n
\\
k\not=j}} \dfrac{\partial}{\partial{t_i}}\left(\dfrac{Q_k}
{\zeta+t_{j}-t_{k}}\right)=\\[1ex]= \sum\limits_{\substack{1\leq k\leq n
\\
k\not=j}}\dfrac{[Q_k,Q_i]}{(\zeta+t_{j}-t_{k})(\zeta+t_{j}-t_{i})},\quad
1\leq i\leq n,\ i\not =j,\end{multline} and
\begin{multline}\label{compfol2}
\dfrac{\partial}{\partial{t_k}}\left(\dfrac{Q_i}{\zeta+t_{j}-t_{i}}\right)-
\dfrac{\partial}{\partial{t_i}}\left(\dfrac{Q_k}{\zeta+t_{j}-t_{k}}\right)=\\[1ex]
=\dfrac{[Q_k,Q_i]}{(\zeta+t_{j}-t_{k})(\zeta+t_{j}-t_{i})},\quad
1\leq i,k\leq n,\ i,k\not =j.\end{multline}
\end{subequations}
The equality \eqref{compfol1} can be simplified as follows:
\begin{multline*}
\dfrac{\partial
Q_i}{\partial{t_j}}\cdot\dfrac{1}{\zeta+t_{j}-t_{i}}+\sum\limits_{\substack{1\leq
k\leq n
\\
k\not=j}} \dfrac{\partial Q_k}{\partial{t_i}}\cdot\dfrac{1}
{\zeta+t_{j}-t_{k}}=\\= \sum\limits_{\substack{1\leq k\leq n
\\
k\not=i,j}}\dfrac{[Q_k,Q_i]}{t_k-t_i}\left(\dfrac{1}{\zeta+t_{j}-t_{k}}-\dfrac{1}{\zeta+t_{j}-t_{i}}\right).
\end{multline*}
In view of \eqref{Sch}, the right-hand side of the last equality can
be rewritten as
\begin{multline*}\sum\limits_{\substack{1\leq k\leq n
\\
k\not=i,j}}\dfrac{\partial Q_k}{\partial
t_i}\left(\dfrac{1}{\zeta+t_{j}-t_{k}}-\dfrac{1}{\zeta+t_{j}-t_{i}}\right)=\\
=\dfrac{\partial Q_j}{\partial
t_i}\cdot\dfrac{1}{\zeta+t_{j}-t_{i}}+\sum\limits_{\substack{1\leq
k\leq n
\\
k\not=j}}\dfrac{\partial Q_k}{\partial
t_i}\cdot\dfrac{1}{\zeta+t_{j}-t_{k}},\end{multline*} where we have
used the fact  that, in view of Remark \ref{firin},
$$\sum\limits_{1\leq k\leq n}\dfrac{\partial Q_k}{\partial
t_i}=0,\quad 1\leq i\leq n.$$ Since, as follows from \eqref{Sch},
$$\dfrac{\partial Q_i}{\partial
t_j}=\dfrac{\partial Q_j}{\partial t_i},\quad 1\leq i,j\leq n,\
i\not=j,$$ we conclude that the equality \eqref{compfol1}, indeed,
holds. The equality \eqref{compfol2} can be verified analogously.

\paragraph{\bf Calculation 2}

We  show that the function $Z_{\alpha_j}(\zeta,\boldsymbol{t})$,
defined in \eqref{defzzz},
 satisfies
with respect to $\zeta$ the equation \eqref{derzz}.

This can be done as follows. We  consider the auxiliary function
\begin{equation}
X_{\alpha_j}(\zeta,\boldsymbol{t})\defi\dfrac{\partial
Z_{\alpha_j}}{\partial \zeta}(\zeta,\boldsymbol{t})-
\sum\limits_{1\leq i\leq
n}\dfrac{Q_i(\boldsymbol{t})}{\zeta+t_j-t_i}\cdot
Z_{\alpha_j}(\zeta,\boldsymbol{t}).
\end{equation}
It will be  shown below that $X_{\alpha_j}(\zeta,\boldsymbol{t})$
satisfies the linear system
\begin{equation}%
\label{regsys2} \left\{
\begin{aligned}
\dfrac{\partial X_{\alpha_j}}{\partial t_{i}}&=
-\dfrac{Q_i}{\zeta+t_{j}-t_{i}}\cdot X_{\alpha_j}, &\quad
 1\leq
i\leq n,\ i\not=j, \\[1.0ex]
\dfrac{\partial X_{\alpha_j}}{\partial t_{j}}&=
\sum\limits_{\substack{1\leq i\leq n
\\
i\not=j}} \dfrac{Q_i} {\zeta+t_{j}-t_{i}}\cdot X_{\alpha_j}.
\end{aligned}\right.
\end{equation}
(Note that this system is the same as the system \eqref{regsys} for
the function $L_j$.)

Since, in view of \eqref{inz}, the solution
$X_{\alpha_j}(\zeta,\boldsymbol{t})$ of the linear system
\eqref{regsys2} satisfies the initial condition
$$X_{\alpha_j}(\zeta,\boldsymbol{t})\bigm|_{\boldsymbol{t}=\boldsymbol{t^0}}=0,$$
 the uniqueness theorem for linear differential systems implies
$$X(\zeta,\boldsymbol{t})\equiv 0.$$ Therefore, the function $Z_{\alpha_j}(\zeta,\boldsymbol{t})$
 satisfies
 the equation \eqref{derzz}.

Now let us prove that $X_{\alpha_j}(\zeta,\boldsymbol{t})$
satisfies, indeed, the linear system \eqref{regsys2}.
 In view of \eqref{derzj}, we have
\begin{multline*}
\dfrac{\partial X_{\alpha_j}}{\partial t_j}= \sum_{\substack{1\leq
i\leq n
\\
i\not=j}}Q_i\cdot\dfrac{\partial}{\partial
\zeta}\left(\dfrac{Z_{\alpha_j}}{\zeta+t_j-t_i}\right)- \sum_{1\leq
i\leq n}\dfrac{\partial}{\partial
t_j}\left(\dfrac{Q_i}{\zeta+t_j-t_i}\right)\cdot Z_{\alpha_j}-\\-
\sum_{\substack{1\leq i,k\leq n\\ k\not =j }}
\dfrac{Q_iQ_k}{(\zeta+t_j-t_i)(\zeta+t_j-t_k)}\cdot
Z_{\alpha_j}\end{multline*}
\begin{multline*}\qquad=\sum_{\substack{1\leq i\leq n
\\
i\not=j}}\dfrac{Q_i}{\zeta+t_j-t_i}\cdot \dfrac{\partial
Z_{\alpha_j}}{\partial \zeta} - \sum_{1\leq i\leq n}\dfrac{\partial
Q_i}{\partial t_j}\cdot\dfrac{1}{\zeta+t_j-t_i}\cdot Z_{\alpha_j} -\\
- \sum_{\substack{1\leq i,k\leq n\\ i\not =j }}
\dfrac{Q_kQ_i}{(\zeta+t_j-t_i)(\zeta+t_j-t_k)}\cdot Z_{\alpha_j}.
\end{multline*}

Now we substitute \eqref{Sch} to obtain
\begin{multline*}
\dfrac{\partial X_{\alpha_j}}{\partial t_j}=\sum_{\substack{1\leq
i\leq n
\\
i\not=j}}\dfrac{Q_i}{\zeta+t_j-t_i}\cdot \dfrac{\partial
Z_{\alpha_j}}{\partial \zeta}-\sum_{\substack{1\leq i\leq n
\\
i\not=j}}\dfrac{Q_iQ_j}{\zeta(\zeta+t_j-t_i)}\cdot Z_{\alpha_j}-\\-
\sum_{\substack{1\leq i,k\leq n\\ i,k\not =j }}
\dfrac{Q_iQ_k}{(\zeta+t_j-t_i)(\zeta+t_j-t_k)}\cdot Z_{\alpha_j}\\
=\left(\sum_{\substack{1\leq i\leq n
\\
i\not=j}}\dfrac{Q_i}{\zeta+t_j-t_i}\right)\cdot\left(\dfrac{\partial
Z_{\alpha_j}}{\partial \zeta}-\sum_{1\leq k\leq n
}\dfrac{Q_k}{\zeta+t_j-t_k}\cdot
Z_{\alpha_j}\right)\\=\sum_{\substack{1\leq i\leq n
\\
i\not=j}}\dfrac{Q_i}{\zeta+t_j-t_i}\cdot X_{\alpha_j}.
\end{multline*}
The first equation of the system \eqref{regsys2} is obtained
analogously.
\end{proof}

\section{Isomondromic and isoprincipal deformations of Fuchsian systems}
\begin{definition}\label{deform}
Let a Fuchsian system \begin{equation}\label{Fu0}
 \dfrac{dY}{dx}=\left(\sum\limits_{1\leq j\leq
n}\dfrac{Q_j^0}{x-t_j^0}\right)\,Y,
\end{equation}
 where
$\boldsymbol{t^0}=(t_1^0,\dots, t_n^0)\in\mathbb{C}^n_*$,
$Q_1^0,\dots,Q_n^0\in\mathfrak{M}_k$ and
\begin{equation}
\label{regcon0} \sum\limits_{1\leq j\leq n}{Q_j^0}=0,\end{equation}
 be given.

Let a holomorphic  family of Fuchsian systems
\begin{equation}
\label{Fut}%
\dfrac{dY}{dx}=\bigg(\sum\limits_{1\leq j\leq
n}\dfrac{Q_j(\boldsymbol{t})}{x-t_j}\bigg) Y,
\end{equation}
where the residue matrices $Q_j(\boldsymbol{t})$ are holomorphic in
a neighborhood $\boldsymbol{\mathcal{D}}\subset\mathbb{C}^n_*$ of
$\boldsymbol{t^0}$, be such that:
\begin{align}
\label{regcontt} \sum\limits_{1\leq j\leq
n}{Q_j(\boldsymbol{t})}&\equiv 0,&\quad
\boldsymbol{t}&\in\boldsymbol{\mathcal{D}},\\
\label{incontt} Q_j(\boldsymbol{t^0})&=Q_j^0,&\quad 1\leq j&\leq n.
\end{align}

If the holomorphic  family of Fuchsian systems \eqref{Fut} is
isoprincipal (respectively,  isomonodromic) with the distinguished
point $x_0=\infty$, then it is said to be an {\em isoprincipal}
(respectively, {\em isomonodromic}) {\em deformation with the
distinguished point} $x_0=\infty$ of the Fuchsian system
\eqref{Fu0}.
\end{definition}

\begin{remark}
Note that, according to Theorem \ref{IsoprImplIsom}, an isoprincipal
deformation \eqref{Fut} of a given Fuchsian system \eqref{Fu0} is
also an isomonodromic one. As Theorem \ref{IsomImplIsopr} implies,
the converse is true under the condition that all the matrices
$Q_1^0,\dots, Q_n^0$ are non-resonant.
\end{remark}

Now we are going to address the question of the existence of an
isoprincipal deformation  of a given Fuchsian system. It follows
from Theorem \ref{main} and Remark \ref{firin} that the holomorphic
family \eqref{Fut} is an isoprincipal deformation with the
distinguished point $x_0=\infty$ of the Fuchsian system \eqref{Fu0}
if and only if the residues $Q_j(\boldsymbol{t})$ satisfy the
Schlesinger system
\begin{equation}%
\label{Sch2} \left\{
\begin{aligned}
\dfrac{\partial Q_{i}}{\partial t_{j}}&=
\dfrac{[Q_{i},Q_{j}]}{t_{i}-t_{j}}, &\quad
 1\leq
i,j\leq n,\ i\not=j, \\[1.0ex]
\dfrac{\partial Q_{i}}{\partial t_{i}}&= - \sum_{\substack{1\leq
j\leq n
\\
j\not=i}} \dfrac{\left[Q_{i},Q_{j}\right]} {t_{i}-t_{j}}, & \quad
1\leq i\leq n.
\end{aligned}\right.
\end{equation}
Thus the question is, whether the Cauchy problem for the Schlesinger
system with the initial condition \eqref{incontt} is solvable. An
answer to this question follows from the Frobenius theorem (Theorem
\ref{genfrob}):
\begin{proposition}\label{Schcomp}
Let $\boldsymbol{t^0}=(t_1^0,\dots, t_n^0)\in\mathbb{C}^n_*$,
$Q_1^0,\dots,Q_n^0\in\mathfrak{M}_k$ be given.

Then there exist a neighborhood
$\boldsymbol{\mathcal{D}}\subset\mathbb{C}^n_*$ of
$\boldsymbol{t^0}$ and unique  matrix functions
$Q_1(\boldsymbol{t}),\dots,Q_n(\boldsymbol{t})$, holomorphic in
$\boldsymbol{\mathcal{D}}$, which satisfy the Schlesinger system
\eqref{Sch2} and the initial condition \eqref{incontt}.
\end{proposition}

\begin{proof}
According to Theorem \ref{genfrob} and Remark \ref{simpfrob}, we
have to to verify the identity
$$\dfrac{\partial^2 Q_i}{\partial t_j\partial t_k}=\dfrac{\partial^2 Q_i}{\partial t_k\partial
t_j},\quad 1\leq i,j,k\leq n,$$ substituting the expressions
\eqref{Sch2} for the partial derivatives of $Q_i.$

In the case $i=k\not=j$ we have
\begin{multline*}
\dfrac{\partial^2 Q_i}{\partial t_i\partial
t_j}=\dfrac{\partial}{\partial{t_j}}\left(\sum_{\substack{1\leq
k\leq n
\\
k\not=i}} \dfrac{\left[Q_{j},Q_{i}\right]}
{t_{i}-t_{j}}\right)\\=\dfrac{[Q_i,[Q_i,Q_j]]-[Q_i,Q_j]}
{(t_{i}-t_{j})^2}+\sum_{\substack{1\leq k\leq n
\\
k\not=i}}\dfrac{[Q_i,[Q_j,Q_k]]+[Q_k,[Q_i,Q_j]]}{(t_i-t_j)(t_i-t_k)}\\
=\dfrac{[Q_i,[Q_i,Q_j]]-[Q_i,Q_j]}
{(t_{i}-t_{j})^2}-\sum_{\substack{1\leq k\leq n
\\
k\not=i}}\dfrac{[Q_j,[Q_k,Q_i]]}{(t_i-t_j)(t_i-t_k)}\\
=\dfrac{\partial}{\partial{t_i}}\left(\dfrac{\left[Q_{i},Q_{j}\right]}
{t_{i}-t_{j}}\right)=\dfrac{\partial^2 Q_i}{\partial t_j\partial
t_i},
\end{multline*}
where we have used the Jacobi identity
$$[A,[B,C]]+[B,[C,A]]+[C,[A,B]]=0,\quad \forall
A,B,C\in\mathfrak{M}_k.$$

In the case $k\not =i\not=j$ the computations are analogous and will
be omitted.
\end{proof}
\begin{remark}
Proposition \ref{Schcomp} was originally proved by L. Schlesinger in
\cite{Sch2}. Here we would also like to mention the paper \cite{Boa}
by P. Boalch, where some general considerations concerning the
compatibility of systems from a class, containing the Schlesinger
system, are presented.
\end{remark}

 As
an immediate consequence of Theorem \ref{main} and Proposition
\ref{Schcomp}, we obtain
\begin{theorem} \label{exisprdef} Let a Fuchsian system \eqref{Fu0} --
\eqref{regcon0} be given.

Then there exists a neighborhood
$\boldsymbol{\mathcal{D}}\subset\mathbb{C}^n_*$ of
$\boldsymbol{t^0}$ and a unique isoprincipal deformation \eqref{Fut}
of the Fuchsian system \eqref{Fu0} with the distinguished point
$x_0=\infty$, where the residue matrices $Q_j(\boldsymbol{t})$ are
holomorphic in the neighborhood $\boldsymbol{\mathcal{D}}$.
\end{theorem}

\begin{remark}
Theorem \ref{exisprdef} can also be proved in another way.

 Let
$\boldsymbol{\mathcal{W}}\Subset\boldsymbol{\mathcal{V}}$ be a
nested pair of open cylindrical neighborhoods of $\boldsymbol{t^0}$
in $\mathbb{C}^n_*$ and for $1\leq j\leq n$  let $\alpha_j$ be a
path from $x_0=\infty$ to a point $p_j\in\mathcal{V}_j\setminus
\overline{\mathcal{W}_j}$. Let $Y(x,\boldsymbol{t^0})$ be the
 solution of \eqref{Fu0} with the initial condition
$$ Y({x},\boldsymbol{t^0})\bigm|_{x=\infty}=I
$$
and let
$$Y_{\alpha_j}({x},\boldsymbol{t^0})=H_j(x,\boldsymbol{t^0})\cdot
E_{\alpha_j}(x-t_j^0),\quad
x\in\cov(\mathcal{V}_j\setminus\{t_j^0\};\,p_j)$$ be the
regular-principal factorization of the appropriate branch
$Y_{\alpha_j}({x},\boldsymbol{t^0})$ in the punctured domain
$\mathcal{V}_j\setminus\{t_j^0\}$ (here the function
$E_{\alpha_j}(\zeta)$, holomorphic and invertible in the Riemann
surface of $\ln\zeta$, is defined in the same way as in the proof of
Proposition \ref{SchlesToIso} -- see \eqref{defprinl}).

First we assume that the family \eqref{Fut}, where the residue
matrices $Q_j(\boldsymbol{t})$ are holomorphic in
$\boldsymbol{\mathcal{W}}$, is an isoprincipal deformation  of the
Fuchsian system \eqref{Fu0} with the distinguished point
$x_0=\infty$ and consider  the solution $Y(x,\boldsymbol{t})$   of
\eqref{Fut} with the initial condition
\begin{equation}\label{inconbo}
Y({x},\boldsymbol{t})\bigm|_{x=\infty}=I.
\end{equation}
Since $Y(x,\boldsymbol{t})$ is jointly holomorphic in $x$ and
$\boldsymbol{t}$ and invertible in the domain
$\cov(\mathbb{CP}^1\setminus\bigcup_k\overline{\mathcal{W}_k};\,\infty)\times
\boldsymbol{\mathcal{W}}$, so is the "ratio"
$$R(x,\boldsymbol{t})\defi Y(x,\boldsymbol{t})\cdot Y^{-1}(x,\boldsymbol{t^0}).$$
Moreover, taking into account the initial condition \eqref{inconbo}
and the fact that the family \eqref{Fut} is by Theorem
\ref{IsoprImplIsom} isomonodromic,  we reach the following
conclusion:
\begin{itemize}
\item[\bf (R)] The function $R(x,\boldsymbol{t})$ is holomorphic (single-valued)
with respect to $x,\boldsymbol{t}$ and invertible in the domain
$\left\{\mathbb{CP}^1\setminus\bigcup_k\overline{\mathcal{W}_k}\right\}\times\boldsymbol{\mathcal{W}}$
and it holds that
$$R({x},\boldsymbol{t})\bigm|_{x=\infty}\equiv I, \quad \boldsymbol{t}\in\boldsymbol{\mathcal{W}}.$$
\end{itemize}

Similarly, considering the coherent family of transplants
$\{E_{\alpha_j}(x-t_j)\}_{t_j\in\mathcal{W}_j},$ we observe that the
"ratio"
$$F_j(x,\boldsymbol{t})\defi E_{\alpha_j}(x-t_j)\cdot
E_{\alpha_j}^{-1}(x-t_j^0),\quad
x\in\cov(\mathbb{C}\setminus\overline{\mathcal{W}_j};\,p_j),\
\boldsymbol{t}\in\boldsymbol{\mathcal{W}}$$   is holomorphic
(single-valued) with respect to $x,\boldsymbol{t}$ and invertible in
$\left\{\mathbb{C}\setminus\overline{\mathcal{W}_j}\right\}\times\boldsymbol{\mathcal{W}}.$

Since the family \eqref{Fut} is assumed to be isoprincipal, there
exist functions $H_j(x,\boldsymbol{t})$ (the regular factors of
$Y(x,\boldsymbol{t})$), such that:
\begin{itemize}
\item[\bf (H)] Each function $H_j(x,\boldsymbol{t})$, $1\leq j\leq
n$,  is holomorphic (single-valued) with respect to
$x,\boldsymbol{t}$ and invertible in
$\mathcal{V}_j\times\boldsymbol{\mathcal{W}};$
\item[\bf (Pb)]
$$
F_j(x,\boldsymbol{t})\cdot
H_j^{-1}(x,\boldsymbol{t^0})=H_j^{-1}(x,\boldsymbol{t})\cdot
R(x,\boldsymbol{t}), \quad
\boldsymbol{t}\in\boldsymbol{\mathcal{W}},\ x\in
\mathcal{V}_j\setminus\overline{\mathcal{W}_j},\ 1\leq j\leq n.
$$
\end{itemize}
Note that the function $F_j(x,\boldsymbol{t})\cdot
H_j^{-1}(x,\boldsymbol{t^0})$ on the left-hand side of the equality
(Pb) is holomorphic (single-valued) with respect to
$x,\boldsymbol{t}$ and invertible in
$\left\{\mathbb{C}\setminus\overline{\mathcal{W}_j}\right\}\times\boldsymbol{\mathcal{W}}$.
It is  determined entirely in terms of the initial data
$\boldsymbol{t^0},Q_1^0,\dots,Q_n^0$. The equality (Pb) itself can
be viewed as a {\em factorization problem} for this function, where
one looks for the factors $H_j(x,\boldsymbol{t})$,
$R(x,\boldsymbol{t})$, possessing the analyticity properties (H) and
(R). If such factors can be found, then, reversing the reasonings
above, one arrives at the isoprincipal deformation of the Fuchsian
system \eqref{Fu0}.

The uniqueness of the solution, possessing  the properties (H) and
(R), for the factorization problem (Pb)   follows immediately from
the Liouville theorem. The existence of this solution for
$\boldsymbol{t}$ in a sufficiently small neighborhood
$\boldsymbol{\mathcal{W}}$ of $\boldsymbol{t^0}$ can be established
by elementary means, since for $\boldsymbol{t}=\boldsymbol{t^0}$ the
solution exists (it is trivial: $R=I$).

However, analyzing the factorization problem (Pb) more carefully and
systematically (see, for instance, \cite[\S 5]{BIK}) and taking into
account Theorem \ref{main}, one can reach the stronger conclusion
that {\em the functions} $Q_j(\boldsymbol{t}),$ {\em which satisfy
the Schlesinger system} \eqref{Sch2} {\em with the initial
condition} \eqref{incontt}, {are meromorphic in the
 universal covering space} $\cov(\mathbb{C}^n_*;,\,\boldsymbol{t^0}).$
 This result was obtained by B. Malgrange in \cite{Ma1} and (in the non-resonant case) by T. Miwa
 in \cite{Miwa}. Our proof, which involves the isoprincipal deformations
 and is outlined above, will be presented in more detail elsewhere.
\end{remark}

\begin{remark}
As was already mentioned (see Remarks \ref{wlog1} and \ref{wlog}),
most of the considerations of the present article need not be
restricted to the case of linear differential systems with only
Fuchsian singularities.

For instance, one can consider linear systems of the form
$$\dfrac{dY}{dx}=\left
(\sum_{j=1}^n\sum_{k=0}^{p_{j}} \dfrac{Q_{j,k}}{(x-
t_{j})^{k+1}}\right )Y,$$ where $$\sum_{j=1}^n Q_{j,0}=0.$$ In this
case the regular-principal factorization of the solution and the
notion of the isoprincipal family
$$\dfrac{dY}{dx}=\left
(\sum_{j=1}^n\sum_{k=0}^{p_{j}} \dfrac{Q_{j,k}(\boldsymbol{t})}{(x-
t_{j})^{k+1}}\right )Y$$ can be introduced  as in Definitions
\ref{DefFact}, \ref{DefIsoprDefo}.

Furthermore, the auxiliary linear system related to the isoprincipal
family takes the form (compare with \eqref{auxisom})
$$\left\{\begin{aligned}\dfrac{\partial Y}{\partial x}&=
\left (\sum\limits_{j=1}^n\sum_{k=0}^{p_{j}}
\dfrac{Q_{j,k}(\boldsymbol{t})}{(x-
t_{j})^{k+1}}\right )Y,\\
\dfrac{\partial Y}{\partial t_j}&=-\left(\sum_{k=0}^{p_{j}}
\dfrac{Q_{j,k}(\boldsymbol{t})}{(x- t_{j})^{k+1}}\right )Y,\quad
1\leq j\leq n.\end{aligned}\right.$$ The compatibility condition for
this overdetermined system is given by the system
$$\left\{
\begin{aligned}
\dfrac{\partial Q_{i,k}}{\partial t_{j}}&=
\sum\limits_{\substack{0\leq\ell\leq{p_{i}-k}\\0\leq q\leq
p_j}}(-1)^l\binom{l+q}{l} \dfrac{ [Q_{i,k+l},Q_{j,q} ]}{( t_{i}-
t_{j})^{q +l+1}} , \quad  0\leq k\leq p_i,\  1\leq
i,j\leq n,\ i\not=j,\\[1ex]
\dfrac{\partial Q_{i,k}}{\partial t_{i}}&= -
\sum\limits_{\substack{1\leq j\leq n
\\
j\not=i}}\sum\limits_{\substack{0\leq\ell\leq{p_{i}-k}\\
0\leq q\leq p_j}}(-1)^l\binom{l+q}{l} \dfrac{ [Q_{i,k+l},Q_{j,q}
]}{( t_{i}- t_{j})^{q +l+1}},  \quad 0\leq k\leq p_i,\ 1\leq i\leq
n,
\end{aligned}\right.$$
which itself is compatible and contains the Schlesinger system
\ref{Sch2} as a special case (corresponding to $p_1=\dots=p_n=0$).
\end{remark}
\section{An example}\label{Examp}
In order to illustrate the distinction between the isoprincipal and
the isomonodromic deformations of Fuchsian systems in the case when
the non-resonance condition is violated, we present the following
explicit example.

 Let us consider the linear differential system
\begin{equation}\label{sysex}
\dfrac{dY}{dx}=\begin{pmatrix}\dfrac{-1}{x(x-1)}&0\\0&\dfrac{-1}{(x-2)(x-3)}\end{pmatrix}Y.\end{equation}
 Note that the system \eqref{sysex} is of the form
 $$\dfrac{dY}{dx}=
\left(\dfrac{Q^0_0}{x}+\dfrac{Q^0_1}{x-1}+\dfrac{Q^0_2}{x-2}+\dfrac{Q^0_3}{x-3}\right)Y,$$
where
$$Q^0_0=-Q^0_1=\begin{pmatrix}1&0\\0&0\end{pmatrix},\quad
Q^0_2=-Q^0_3=\begin{pmatrix}0&0\\0&1\end{pmatrix},$$ and hence it is
Fuchsian and resonant. Moreover, $x=\infty$ is a regular point of
the system \eqref{sysex}, since
$$Q^0_0+Q^0_1+Q^0_2+Q^0_3=0.$$ The  solution $Y(x)$ satisfying the
initial condition $Y_{|\infty}=I$ is a rational matrix functions (in
particular, single-valued); it has the following form:
$$Y(x)=\begin{pmatrix}\dfrac{x}{x-1}&0\\0&\dfrac{x-2}{x-3}\end{pmatrix}\stackrel{\textup{\tiny
def}}{=}Y^0(x).$$ The principal factors of $Y^0$ can be
chosen\,\footnote{See Remark \ref{regau}.} in the form
\begin{alignat*}{2}P_0(x)&=\begin{pmatrix}x&0\\0&1\end{pmatrix},&\quad
P_1(x)&=\begin{pmatrix}\dfrac{1}{x-1}&0\\0&1\end{pmatrix},\\[1ex]
P_2(x)&=\begin{pmatrix}1&0\\0&x-2\end{pmatrix},&\quad
P_3(x)&=\begin{pmatrix}1&0\\0&\dfrac{1}{x-3}\end{pmatrix}.\end{alignat*}

 Now  we are going to construct the isoprincipal deformation of
the system \eqref{sysex}. For simplicity, we assume that the
singular points $x=1,2,3$ are fixed while the position of one
singularity (x=0) varies: $x=t$. Then we need to determine a
holomorphic family of  matrix functions $Y(x,t)$ such that
\begin{equation}\label{noconex}Y(x,t)\bigm|_{x=\infty}=I,\quad
Y(x,t)\bigm|_{t=0}=Y^0(x),\end{equation}
 $P_0(x-t)$
is the principal factor of $Y(x,t)$ at $x=t$  and for $k=1,2,3$
$P_k(x)$ is the principal factor of $Y(x,t)$ at $x=k$. Such a family
is unique and consists of rational matrix functions of the form:
$$Y(x,t)=\begin{pmatrix}\dfrac{x-t}{x-1}&0\\0&\dfrac{x-2}{x-3}\end{pmatrix}.$$
The function $Y(x,t)$ satisfies with respect to $x$ the Fuchsian
system
\begin{multline}\label{isoprdex}\dfrac{dY}{dx}=\begin{pmatrix}\dfrac{t-1}{(x-t)(x-1)}&0\\0&\dfrac{-1}{(x-2)(x-3)}
\end{pmatrix}Y\\[1ex]=
\left(\dfrac{Q_0^0}{x-t}+\dfrac{Q_1^0}{x-1}+\dfrac{Q_2^0}{x-2}+\dfrac{Q_3^0}{x-3}\right)Y,\end{multline}
and the constant functions
$$Q_j(t)\equiv Q_j^0,\quad j=0,1,2,3,$$
satisfy, of course, the Schlesinger system
$$\left\{\begin{aligned}\dfrac{d Q_j}{dt}&=\dfrac{[Q_0,Q_j]}{t-j},\quad
j=1,2,3,\\\dfrac{d Q_0}{dt}&=\sum_{j=1}^3\dfrac{[Q_j,Q_0]}{t-j},
\end{aligned}\right.$$
with the initial condition
$$Q_j(0)=Q_j^0,\quad j=0,1,2,3.$$

However, the deformation \eqref{isoprdex} is not the only possible
isomonodromic deformation of the system \eqref{sysex}. Indeed, let
us consider a family of rational functions
$$Y(x,t)=\begin{pmatrix}\dfrac{x-t}{x-1}&\dfrac{-2t(x-t)h(t)}{(x-1)(x-3)(t-3)}\\[1ex]
0&\dfrac{x-2}{x-3}\end{pmatrix},$$ where $h(t)$ is a function,
holomorphic in a neighborhood of $t=0$ (outside this neighborhood
$h(t)$ may have arbitrary singularities). Such a family also
satisfies the normalizing conditions \eqref{noconex}; moreover, we
have
\begin{align*}
Y(x,t)^{-1}&=\begin{pmatrix}\dfrac{x-1}{x-t}&\dfrac{2th(t)}{(x-2)(t-3)}\\[1ex]
0&\dfrac{x-3}{x-2}\end{pmatrix},\\[1ex] \dfrac{\partial
Y(x,t)}{\partial
x}Y(x,t)^{-1}&=\begin{pmatrix}\dfrac{t-1}{(x-t)(x-1)}&\dfrac{2t(x-t)h(t)}{(x-1)(x-2)(x-3)(t-3)}\\[1ex]
0&\dfrac{-1}{(x-2)(x-3)}\end{pmatrix}.\end{align*} Thus we obtain
the following deformation of the system \eqref{sysex}:
\begin{multline}\label{isomondex}\dfrac{dY}{dx}=\begin{pmatrix}\dfrac{t-1}{(x-t)(x-1)}&
\dfrac{2t(x-t)h(t)}{(x-1)(x-2)(x-3)(t-3)}\\[1ex]
0&\dfrac{-1}{(x-2)(x-3)}\end{pmatrix}Y\\[2ex]=
\left(\dfrac{Q_0(t)}{x-t}+\dfrac{Q_1(t)}{x-1}+\dfrac{Q_2(t)}{x-2}+\dfrac{Q_3(t)}{x-3}\right)Y,\end{multline}
where
\begin{gather*}
Q_0(t)\equiv Q_0^0,\quad
Q_1(t)=\begin{pmatrix}-1&\dfrac{-t(t-1)h(t)}{t-3}\\0&0\end{pmatrix}\\[2ex]
Q_2(t)=\begin{pmatrix}0&\dfrac{2t(t-2)h(t)}{t-3}\\0&1\end{pmatrix},\quad
Q_3(t)=\begin{pmatrix}0&-th(t)\\0&-1\end{pmatrix}.
\end{gather*}
Since the monodromy of $Y(x,t)$ for every $t$ is trivial, the
deformation \eqref{isomondex} is isomonodromic, but if
$h(t)\not\equiv 0$, then it is not isoprincipal and the coefficients
$Q_j(t)$ do not satisfy the Schlesinger system. Moreover, since we
only require $h(t)$ to be holomorphic in a neighborhood of $t=0$,
the behavior of the functions $Q_j(t)$ outside this neighborhood may
be arbitrary. In particular, these functions $Q_j(t)$  need not be
  meromorphic with respect to $t$ in $\mathbb{C}\setminus\{1,2,3\}$.

\begin{remark}\label{ashkelon} The example presented above is based on the theory
of holomorphic families \eqref{Fut} of Fuchsian systems, whose
solutions $Y(x,\boldsymbol{t})$ are generic rational matrix
functions of $x$ for every fixed $\boldsymbol{t}$. This theory was
developed in \cite{Kats1}, \cite{Kats2} and  \cite{KaVo2} (see also
the electornic version of the latter work \cite{KaVo2-e}.) For such
Fuchsian systems the number  of poles $n$ is even, so we write $2n$
instead of $n$. All such Fuchsian systems with this property can be
parameterized as follows. If $k$ is the dimension of the square
residue matrices $Q_j$, then to each pole $t_j$ a $k-1$-vector is
related as a "free" parameter. Therefore, to each $2n$-tuple
$\boldsymbol{t}=(t_1, \ldots ,t_2n)$ the total of $(k-1)\times 2n$
complex parameters is related. To every choice of these parameters
(satisfying a certain non-degeneracy condition) corresponds a
different system of the form \eqref{Fut}, whose solution
$Y(x,\boldsymbol{t})$ is a generic rational matrix function.
Considering $\boldsymbol{t}$ as variable, we assign $(k-1)\times 2n$
complex parameters to each
$\boldsymbol{t}\in\mathbb{C}^{2n}_{\ast}$.

Thus the families \eqref{Fut}, possessing the property that
$Y(x,\boldsymbol{t})$ is a generic rational matrix function of $x$
for each fixed $\boldsymbol{t}$, can be parameterized by
$(k-1)\times 2n$ complex valued functions of $\boldsymbol{t}$, and
these {\em functional} parameters are free. To obtain holomorphic
families, we have to require that these complex valued functions are
holomorphic. Of course, the monodromy of any rational matrix
function of $x$ is trivial. Hence we can parameterize the class of
the {\em isomonodromic} deformations \eqref{Fut} of Fuchsian systems
with generic rational solutions by $(k-1)\times 2n$ {functional}
parameters, and these parameters can be {\em arbitrary} (up to the
mentioned non-degeneracy condition) holomorphic functions of
$\boldsymbol{t}$.

It turns out that the deformation corresponding to a given choice of
the functional parameters is isoprincipal if and only if all these
parameters are {\em constant} functions. In particular, the class of
the isomonodromic deformations (considered without the non-resonance
condition) is much richer than its subclass of the isoprincipal
deformations.
\end{remark}


\begin{thebibliography}{MoMo2}
\small
\bibitem[Birk1]{Birk1}
\textsc{Birkhoff, G.D.} \textsl{The generalized {Riemann} problem
for linear differential equations and the allied problems for linear
difference and q-difference equations.} Proc. Amer. Acad.
Arts and Sci., \textbf{49} (1913), 521 - 568.\\
 \small{Reprinted in:} \cite{Birk2}, 259 - 306.
\bibitem[Birk2]{Birk2}
\textsc{Birkhoff, G.D.} \textsl{Collected Mathematical Papers.
Vol.I}. American Mathematical Society, New York, 1950. i-xi, 754 pp.
\bibitem[Boa]{Boa} \textsc{Boalch, P.} \textsl{Symplectic manifolds
and isomonodromic deformations.} Adv. in Math. \textbf{163} (2001),
137\,-\,205.
\bibitem[Bol1]{Bol1} {\fontfamily{UWCyr}\selectfont\cyracc
\textsc{Bolibruh, A. A.} \textit{Problema Rimana-Gil{\cprime}berta}.
Uspehi matem. nauk, \textbf{45}:2 (1990),  s. 3 - 47 .} Engl.
transl.: \textsc{Bolibruch, A.A.} \textsl{The Riemann-Hilbert
problem.} Russian Math. Surveys, \textbf{45}:2 (1990), pp. 1 - 47.
\bibitem[Bol2]{Bol2} {\fontfamily{UWCyr}\selectfont\cyracc
\textsc{Bolibrukh, A.} \textit{21-ya problema Gil{\cprime}berta dlya
line\u{\i}nykh fuksovykh sistem}. Nauka, 1994 (Tr. Mat. Inst.
Steklova, Ross. Akad. Nauk, \textbf{206}).} (In Russian). English
transl.: \textsc{Bolibrukh, A.} \textsl{The 21st Hilbert Problem for
Linear Fuchsian Systems.} Providence: Am. Math. Soc., 1995 (Proc.
Steklov Inst. Math., \textbf{206}).
\bibitem[Bol3]{Bol3}
\textsc{Bolibruch, A.} \textsl{On isomonodromic deformations of
Fuchsian systems.} Journ. of Dynamical and Control Systems, Vol.
\textbf{3}, No.4 (1997), 589\,-\,604.
\bibitem[Bol4]{Bol4}
{\fontfamily{UWCyr}\selectfont\cyracc \textsc{Bolibrukh, A.}
\textit{Ob izomonodromnykh sliyaniyakh fuksovykh osobennoste\u{\i}}.
Trudy Matem. Inst. Steklova, Ross. Akad. Nauk,} \textbf{221} (1998),
127\,-\,142. (In Russian). English transl.: \textsc{Bolibrukh, A.}
\textsl{On isomonodromic confluences of fuchsian singularities.}
Proc. Steklov Inst. Math., \textbf{221}(1998), 117\,-\,132.
\bibitem[Bol5]{Bol5} {\fontfamily{UWCyr}\selectfont\cyracc \textsc{Bolibrukh,\,A.}
\textit{Differentsial{\cprime}nye uravneniya s meromorfnymi
ko\-{\`e}f\-fi\-tsi\-en\,ta\,mi}. Str. 29\,-\,82 v: Sovremennye
Problemy Matematiki, Vypusk} \textbf{1},
{\fontfamily{UWCyr}\selectfont\cyracc  Matem. Inst. Steklova,
Ross. Akad. Nauk, Moskva,} 2003 (In Russian).\\ %
 \textsc{Bolibrukh, A.}
\textsl{Differential equations with meromorphic coefficients,}
pp.\.29\,-\,82 in: Modern Problems of Mathematics, \textbf{1},
Steklov Inst. Math., Russian Acad. Sci., Moskow, 2003.
\bibitem[BIK]{BIK} \textsc{Bolibruch,\,A.A., A.R.\,Its} and
{A.A.\,Kapaev}. \textsl{On the Riemann-Hilbert-Birkhoff inverse
monodromy problem and the Painlev\'e equations}.
{\fontfamily{UWCyr}\selectfont\cyracc Algebra i analiz, tom}
\textbf{16}:1 (2004), 121\,-\,162, and also St. Petersburg
Mathematical Journal, vol.\,\textbf{16}:1 (2005), pp.\,105\,-\,142.
\bibitem[Fo]{Fo} \textsc{Forster, O.} \textsl{Riemannschen
Fl\"achen.} Springer-Verlag,\\
Berlin\(\boldsymbol{\cdot}\)Heidelberg\(\boldsymbol{\cdot}\)New\,York,
1977, x+223 pp. (in German).\\ %
English Transl.: \textsc{Forster, O.} \textsl{Lectures on Riemann
Surfaces} (Graduate Texts in
Mathematics, \textbf{81}.) Springer-Verlag, %
New\,York\(\boldsymbol{\cdot}\)Berlin, 1981. viii+254 pp.\\
Russian Transl.: {\fontfamily{UWCyr}\selectfont\cyracc
\textsc{Forster, O.} \textit{Rimanovy Poverkhnosti}, Mir, Moskva,
1980.}
\bibitem[FuL]{FuL} \textsc{Fuchs, L.} \textsl{Gesammelte Mathematische Werke.
Band 3.} {\small (Herausgegeben by \textsc{R. Fuchs} und \textsc{L
Schlesinger})}. Mayer \& M\"uller, Berlin, 1909.
\bibitem[FuL1]{FuL1} \textsc{Fuchs, L.} \textsl{Zur Theorie der linearen
Differentialgleichungen.} Sitzungsberichte der K. preuss. Akademie
der Wissenschaften zu Berlin. Einleitung und No. 1 - 7, 1888, S.
1115 - 1126; No. 8 - 15, 1888, S.1273 - 1290; No. 16 - 21, 1889, S.
713 - 726; No. 22\.- 31, 1890, S. 21 - 38. \small{Reprinted in:}
\cite{FuL}, S. 1 - 68.
\bibitem[FuL2]{FuL2} \textsc{Fuchs, L.}  \textsl{\"Uber lineare
Differentialgleichungen, welche von Parametern unabh\"angige
substitutionsgruppen besitzen.} Sitzungsberichte der K. preuss.
Akademie der Wissenschaften zu Berlin. 1892, S. 157 - 176.
\small{Reprinted in:} \cite{FuL}, S. 117 - 139.
\bibitem[FuL3]{FuL3} \textsc{Fuchs, L.}  \textsl{\"Uber lineare
Differentialgleichungen, welche von Parametern unabh\"angige
Substitutionsgruppen besitzen.} Sitzungsberichte der K. preuss.
Akademie der Wissenschaften zu Berlin. Einleitung und No. 1 - 4,
1893, S. 975 - 988; No. 5 - 8, 1894,, S. 1117 - 1127.
\small{Reprinted in:} \cite{FuL}, S. 169 - 195.
\bibitem[FuL4]{FuL4} \textsc{Fuchs, L.}  \textsl{\"Uber die Abh\"angigkeit
der L\"osungen einer linearen differentialgleichung von den in den
Coefficienten auftretenden Parametren.} Sitzungsberichte der K.
preuss. Akademie der Wissenschaften zu Berlin. 1895, S. 905 - 920.
\small{Reprinted in:} \cite{FuL}, S. 201 - 217.
\bibitem[FuR1]{FuR1} \textsc{Fuchs, R.} \textsl{Sur quelquea \'equations
diff\'erentielles lin\'eares du second ordre.} %
Compt. Rend. de l'Acad\'emie des Sciences, Paris. \textbf{141}
(1905), pp. 555 - 558.
\bibitem[FuR2]{FuR2} \textsc{Fuchs, R.} \textsl{\"Uber lineare homogene
Differentialgleichungen zweiter Ordnung mit drei im Endlichen
gegebenen wesentlich
singul\"aren Stellen.} %
Mathematische Annalen, \textbf{63} (1907), pp. 301 - 321.
\bibitem[Gant]{Gant} {\fontfamily{UWCyr}\selectfont\cyracc
\textsc{Gantmaher, F. R.} \textit{Teoriya matric. 2-e izdanie}.
Nauka, 1966, 575 s.} (In Russian). English transl.:
\textsc{Gantmacher, F.R.}. \textsc{The theory of matrices.} Chelsea,
New York, 1959, 1960.
\bibitem[Gar]{Gar} \textsc{Garnier, R.} \textsl{Sur une classe d'\'equations
diff\'erentielles dont les int\'egrales g\'en\'erales ont leurs
points critiques fixes.} Compt. Rend. de l'Acad\'emie des Sciences,
Paris. \textbf{151} (1910), pp. 205 - 208.
\bibitem[HW]{HW} \textsc{Hurewich, W.} and
\textsc{H. Wallman}: \textsl{Dimension Theory.} Princeton University
Press, Princeton, 1941.
\bibitem[IKSY]{IKSY} \textsc{Iwasaki,\,K., H.\,Kimura,\,
S.\,Shimomura} and \textsc{M.\,Yoshida.} \textsl{From Gau{\ss} to
Painlev\'e.} ({\footnotesize\textsl{A Modern Theory of Special
Functions.}}) Aspects of Mathematics: E; vol. \textbf{16}.
Friedr.\,Vieweg\,\&\,Sohn, Braunschweig, 1991, xii+347
pp. 
\bibitem[ItNo]{ItNo} \textsc{Its,\,A.R.} and \textsc{V.Yu.\,Novokshenov}.\textsl{The Isomonodromic Deformation Method in
 the Theory of Painlev{\'e} equations}, Lect. Notes in Math, \textbf{1191},
 Springer-Verlag, Berlin\,-\,New\,York, 1986.
\bibitem[Kats1]{Kats1} \textsc{Katsnelson, V.} \textsl{Fuchsian differential
systems related to rational matrix fuctions in general position and
the joint system realization ,} pp. 117 - 143  in: Israel
Mathematical Conference Proceedings, Vol. \textbf{11} (1997),
\textsl{Proceedings of the Ashkelon Workshop on Complex Function
Theory (May 1996)}, \textsc{Zalcman, L.} - editor.
\bibitem[Kats2]{Kats2} \textsc{Katsnelson, V.} \textsl{Right and left joint
system representation of a rational matrix function in general
position (System representation theory for dummies) ,} pp. 337 - 400
in: \textsl{Operator theory, system theory and related topics. (The
Moshe Liv\v{s}ic anniversary volume. Proceedings of the Conference
on Operator Theory held at Ben-Gurion University of the Negev,
Beer-Sheva and in Rehovot, June 29 - July 4, 1997)}, \textsc{Alpay,
D.} and \textsc{V. Vinnikov} - editors. (Operator Theory:
Advances and Applications, vol. \textbf{123}), %
Birkh\"auser Verlag, Basel, 2001.
\bibitem[KaVo1]{KaVo1} \textsc{Katsnelson, V.} and
\textsc{D. Volok}: \textsl{Rational Solutions of the Schlesinger
System and Isoprincipal Deformations of Rational Matrix Functions I}
, pp.  291 - 348 in: \textsl{Current Trends in Operator Theory and
its Applications}, \textsc{Ball, J.A., Helton, J.W., Klaus, M.} and
\textsc{L. Rodman} -  editors. (Operator Theory:
Advances and Applications, vol. \textbf{149}), %
Birkh\"auser Verlag, Basel, 2004.
\bibitem[KaVo1-e]{KaVo1-e} \textsc{Katsnelson, V.} and
\textsc{D. Volok}: \textsl{Rational Solutions of the Schlesinger
System and Isoprincipal Deformations of Rational Matrix Functions
I},\\
arXiv.org e-print archive: \texttt{http://arxiv.org, math.
CA/0304108.}
\bibitem[KaVo2]{KaVo2} \textsc{Katsnelson, V.} and
\textsc{D. Volok}: \textsl{Rational Solutions of the Schlesinger
System and Isoprincipal Deformations of Rational Matrix Functions
II}, pp.  165 - 203 in: \textsl{Operator Theory, System Theory And
Scattering Theory: Multidimensional Generalizations}, \textsc{Alpay,
D.} and \textsc{V. Vinnikov} - editors. (Operator Theory:
Advances and Applications, vol. \textbf{157}), %
Birkh\"auser Verlag, Basel, 2005.
\bibitem[KaVo2-e]{KaVo2-e} \textsc{Katsnelson, V.} and
\textsc{D. Volok}: \textsl{Rational Solutions of the Schlesinger
System and Isoprincipal Deformations of Rational Matrix Functions
II},\\
arXiv.org e-print archive: \texttt{http://arxiv.org, math.
CA/0406489.}
\bibitem[Ma1]{Ma1} \textsc{Malgrange, B.}
\textsl{Sur les deformations isomonodromiques. I. Singularites
regulieres} (In French) [On isomonodromic deformations. I. Regular
singularities.] Pp. 401\,-\,426 in \cite{MaPh}.
\bibitem[Ma2]{Ma2} \textsc{Malgrange, B.}
\textsl{Sur les deformations isomonodromiques. II. Singularites
irregulieres} (In French) [On isomonodromic deformations. II.
Irregulare singularities.] Pp. 427\,-\,438 in \cite{MaPh}.
\bibitem[MaPh]{MaPh}\textsl{Math\'ematique et Physique.
(\footnotesize S\'eminaire de l'Ecole Normale Sup\'erieure
1979\,-\,1982).} (In French). [Mathematics and Physics.
(\footnotesize Seminar of Ecole Normale Sup\'erieure
1979\,-\,1982.)]\textsc{Boutet de Monvel,\,L., Douady,\,A., Verdier,
J.-L.}\,-\,eds. (Progress in Mathematics, \textbf{37}),
Birkh\"auser, %
Boston\(\boldsymbol{\cdot}\)Basel\(\boldsymbol{\cdot}\)Stuttgart, %
1983.

\bibitem[Miwa]{Miwa} \textsc{Miwa, T.} \textsl{Painleve' property of monodromy preserving
deformation equations and the analyticity of $\tau $ functions.}
Publ. Res. Inst. Math. Sci.  \textbf{17}  (1981), no. 2, 703 - 721.
\bibitem[Na]{Na}\textsc{Narasimhan, R.} \textsl{Analysis on real and complex manifolds.}
(Advanced Studies in Pure Mathematics, Vol. \textbf{1}).
Masson\,\&\,Cie, \'Editeur, Paris; North-Holland Publishing Co.,
Amsterdam; 1973. x+246 pp.\\ %
Russian Transl.: {\fontfamily{UWCyr}\selectfont\cyracc
\textsc{Narasimhan, R.} \textit{ Analiz na
de\u{\i}stvitel{\cprime}nyh i kompleksnyh mnogoobraziyah.} MIR,
Moskva 1971, 232 s.}
\bibitem[Sch1]{Sch1} \textsc{Schlesinger, L.}
\textsl{Sur la d{\'e}termination des fonctions alg{\`e}briques
uniformes sur une surface de  Riemann donn{\'e}e.} (French.) Annales
Scientifiques de l'{\'E}cole Normale Sup{\'e}rieure (S\'er.
\textbf{3}), t.\textbf{20} (1903), p. 331-347.
\bibitem[Sch2]{Sch2} \textsc{Schlesinger, L.} \textsl{\"Uber die L\"osungen
gewisser linearer Differentialgleichungen als Funktionen der
singul\"aren Punkte.} Journal f\"ur reine und angew. Math,
\textbf{129} (1905), pp. 287 - 294.
\bibitem[Sch3]{Sch3} \textsc{Schlesinger, L.} \textsl{Vorlesungen \"{u}ber
lineare Differentialgleichungen.} Leipzig und Berlin, 1908.
\bibitem[Sch4]{Sch4} \textsc{Schlesinger, L.} \textsl{\"Uber eine Klasse von
Differentialsystemen beliebiger Ordnung mit festen kritischen
Punkten.} Journal f\"ur reine und angew. Math, \textbf{141} (1912),
pp. 96 - 145.



\end{thebibliography}
\end{document}